\documentclass[11pt]{article}
\usepackage{graphicx}
\usepackage{amssymb}
\usepackage{eucal}
\usepackage[bottom,none]{draftcopy}
\newcommand\p{{\partial}}

\newcommand\supp{{\mbox{supp}\,}}                              %
\newcommand\BbbZ{{\mathbb Z}} 
\newcommand\BbbN{{\mathbb N}}                                                
\newcommand\BbbR{{\mathbb R}}                                                
\newcommand\BbbC{{\mathbb C}}  
\newcommand\BbbS{{\mathbb S}} 

\newtheorem{theorem}{Theorem}[subsection]
\newtheorem{corollary}{Corollary}[subsection]
\newtheorem{definition}{Definition}[subsection]
\newtheorem{remark}{Remark}[subsection]
\newtheorem{lemma}{Lemma}[subsection]
\newtheorem{proposition}{Proposition}[subsection]
\newtheorem{claim}{Claim}[subsection]

\makeatletter
\@addtoreset{equation}{section}

\makeatother

\title{Variable coefficient Schr\"odinger flows
for ultrahyperbolic operators}
\author{C.E. Kenig\thanks{Partially supported by  NSF and IBERDROLA program of Profesores
Visitantes} \\
University of Chicago, Chicago Il., USA 60637, cek@math.uchicago.edu 
\and
G. Ponce\thanks{Partially supported by  NSF and IBERDROLA program of Profesores Visitantes} \\
University of California, Santa Barbara, Ca, USA ponce@math.ucsb.edu
\and
C. Rolvung\\
Nykredit Markets \& Asset Management,\\
Kalvevod Bridge 1--3, DK-1780 Copenhagen V, Denmark,
cro@nykredit.dk
\and
 L. Vega\thanks{Partially supported by a MECyT grant}\\
 Universidad del Pais Vasco, Apdo. 644, 48080, Bilbao, Spain, mtpvegol@lg.ehu.es}
\date{}
\begin{document}
\maketitle
\newpage

\tableofcontents

\newpage

\section{INTRODUCTION }\label{sec1}
 In this paper we shall consider nonlinear
Schr\"odinger equations of the form
\begin{equation}
 \begin{array}{rcl}
\p_t u&=& i {\cal L}(x) u + \vec b_1(x)\cdot \nabla_x u + \vec b_2(x)\cdot \nabla_x \bar u \\
&+& c_1(x)u+c_2(x)\bar u + P(u,\bar u,\nabla_x u, \nabla_x\bar u),\label{eq1.1}
 \end{array}
\end{equation}
where $x\in\mathbb R^n$, $t>0$, $\displaystyle{\cal L}(x)=-\sum_{j,k=1}^n\p_{x_j}\left(a_{jk}(x)\p_{x_k}\right)$, 
$A(x)=\left(a_{jk}(x)\right)_{j,k=1,..,n}$
 is a real, symmetric and nondegenerate variable coefficient matrix, and $P$ is a
polynomial  with no linear or constant terms.

 Equations of the form described in (\ref{eq1.1}) with  $A(x)$ merely invertible as opposed
to positive definite arise in connection with water wave problems,  and in higher
dimensions as completely integrable models, see \cite{5}, 
\cite{14}, and \cite{32}.

 In this work we shall study the existence, uniqueness and regularity of the local solutions to the
initial value problem (IVP) associated to the equation (\ref{eq1.1}).
The class of equations is rather general and appropriate assumptions have to be imposed on the smoothness and decay of the
coefficients
$a_{jk}, \vec b_1, \vec b_2, c_1$ and $c_2$ and on the initial data $u(x,0)=u_0(x)$ as well as on the asymptotic
behavior of $a_{jk}(x)$ as $|x|\to \infty$. Also it will be necessary to measure
the regularity of
solutions in weighted Sobolev spaces of high indexes. The main result we obtain in this direction
is Theorem \ref{th6.2.2} -also see Remark \ref{re6.2.3}, in Section \ref{sec6}.
 
 One of the main difficulties in (\ref{eq1.1}) is that the nonlinear terms incur in the so called \lq\lq loss of derivatives". 
This can be avoided if $P$ 
is assumed to have a special symmetric form and $\vec b_1$ is real valued. 
In this case, the standard energy method gives local 
well-posedness of the corresponding IVP in $H^s(\mathbb R^n)$ for $s>n/2+1$ 
independently of the dispersive nature of (\ref{eq1.1}), 
see \cite{15},
\cite{23}. Another approach used to overcome this loss of derivatives is to restrict oneself to
working with 
${\cal L}=\Delta$, $\vec b_1=\vec b_2=0$
 in   suitable analytic function spaces, see \cite{10} and references therein.

In \cite{18}, C. E. Kenig, G. Ponce and L. Vega used linear dispersive smoothing effects of the associated linear equation
to show that the IVP for the equation (\ref{eq1.1}) with ${\cal L} =\Delta$, $\vec b_1= \vec b_2=0$, 
$c_1=c_2=0$ and general $P$ is locally
well posed in (possibly weighted) Sobolev spaces with high index for small initial data. For the case $n=1$, N. 
Hayashi and T. Ozawa \cite {11} removed the smallness condition by using an integrating factor
which  reduces the problem to a system 
where the energy method applies. H. Chihara \cite{2} removed the smallness assumptions 
in weighted Sobolev spaces \cite{18} in any dimension
$n$ by considering systems of two equations
which he diagonalized to essentially eliminate the conjugate first order terms. The remaining first
order terms are treated by a  method similar to the one used by S. Mizohata \cite{26} and S. Doi
\cite{7} to solve  linear  Schr\"odinger equations with lower order terms. It consists in applying
a pseudo-differential operator
$K$ to the equation. The commutator
$i[K\Delta-\Delta K]$  basically absorbs the first orderterm to overcome the loss of derivatives in a way related to the method of integrating factors. 
In Chihara's approach the ellipticity of ${\cal L}=\Delta$ 
is  key in the
diagonalization argument.
C. E. Kenig, G. Ponce and L. Vega \cite {22} obtained local well-posedness 
for the IVP for (\ref{eq1.1}) in the non-elliptic constant coefficient case
$$
{\cal L}=\sum_{j=1}^k\p_{x_j}^2-\sum_{j=k+1}^n \p_{x_j}^2,\;\;\;\;\;\;\;\,k=1,..,n-1,
$$
using the pseudo-differential operators of \cite{1} in the linear problem to avoid the 
diagonalization process. Furthermore, their results
are valid in Sobolev spaces with no weights if $P$ has no quadratic terms.
 
 It may be gathered from this short summary of background literature that a thorough understanding
of  linear Schr\"odinger equations is
important in the attempt to solve the nonlinear problem for (\ref{eq1.1}).

 Our approach in this work will be illustrated with the special case of (\ref{eq1.1})
\begin{equation}
\left\{
\begin{array}{l}
\p_t  u =i{\cal L}(x)u + \vec b_1(x)\cdot\nabla_x u +u\p_{x_1}u,\\
u(x,0)=u_0(x),
\end{array}
\right.\label{eq1.2}
\end{equation}
or equivalently
\begin{equation}
\left\{
\begin{array}{l}
\p_t  u =i{\cal L}(x)u + [\vec b_1(x)\cdot\nabla_x u +u_0(x)\p_{x_1}u] +(u-u_0)\p_{x_1}u,\\
u(x,0)=u_0(x).
\end{array}
\right. \label{eq1.3}
\end{equation}

The nonlinear part $(u-u_0)\p_{x_1}u$ of (\ref{eq1.3}) should be small for small $t$ because of the factor
$(u-u_0)$, but the factor $\p_{x_1}u$ still incurs loss of one derivative. The linear 
part of (\ref{eq1.3}) has a modified first order coefficient.
 
 It is therefore useful to study linear Schr\"odinger equations of the form
\begin{equation}
\left\{
\begin{array}{l}
\p_t  u =i{\cal L}(x)u + \vec b_1(x)\cdot\nabla_x u + \vec b_2(x)\cdot\nabla_x \bar u+c_1(x)u+c_2(x)\bar u + f(x,t),\\
u(x,0)=u_0(x).
\end{array}
\right.\label{eq1.4}
\end{equation}

 Solutions $u$ of (\ref{eq1.4})
gain one derivative compared to
$f$ and
$1/2$ derivative compared to
$u_0$,
 on the average in time and modulo spatial weights, under suitable assumptions. More
precisely,
 for
$s\in\BbbZ^+$ and $N\in \BbbZ^+$, $N>1$, the solution of (\ref{eq1.4}) satisfies
\begin{equation}
\begin{array}{l}
\displaystyle\int_0^T\int_{\BbbR^n}\left|J^{s+1/2}u(x,t)\right|\langle x\rangle^{-N}dxdt\\
\qquad\displaystyle\leq c\left((1+T)\|u_0\|_{H^s}^2 + \int_0^T\int_{\BbbR^n}\left|J^{s-1/2}f(x,t)\right|\langle x\rangle^{N}dxdt\right)
\end{array}\label{eq1.5}
\end{equation}

For the proof of these results
for the constant coefficient case ${\cal L}(x)=\Delta$,  $\vec b_1=\vec b_2=0$ 
see \cite{16}, \cite{24}, \cite{3}, \cite{29}, \cite{31},
\cite{17}, \cite{18}). The estimate (\ref{eq1.5}) allows to overcome the loss of one derivative 
introduced by  the nonlinear part of (\ref{eq1.3})
and, more generally, 
to solve (\ref{eq1.1}).

 Consider (\ref{eq1.4}) with ${\cal L}(x)=\Delta$, $\vec b_1=(i,0,..,0)$, $b_2=c_1=c_2=f=0$. 
The solution of this constant coefficient equation is
given via Fourier transform by
$ \hat u(\xi,t)=\exp(-t(i|\xi|^2+\xi_1))\hat u_0(\xi)$. The multiplier
 $\exp(-t(i|\xi|^2+\xi_1))$ is unbounded for $t\ne 0$, so (\ref{eq1.4}) is not wellposed in $L^2$ in this case.
In fact, the following condition, deduced by S. Mizohata \cite{26}, has been proven to be necessary for the well-posedness in
$L^2$ of (\ref{eq1.4}) with
${\cal L}(x)=\Delta$, $\vec b_2=0$
\begin{equation}
\sup_{x\in{\BbbR}^n,\,\omega\in{\BbbS}^{n-1}}\;\left|\int_0^{\infty}\mbox{Im} \,\vec b_1(x+r\omega)\cdot \omega
dr\right|<\infty.\label{eq1.6}
\end{equation}
So the decay assumptions on $\mbox{Im} \,\vec b_1$ are natural. An application of $J^s=(I-\Delta)^{s/2}$
 to (\ref{eq1.4}) gives a new equation with $\nabla_x a_{jk}(x)$ appearing in the first order coefficient. 
 Well posedness of (\ref{eq1.4}) is of interest for
any $s$, so decay assumptions on  $\nabla_x a_{jk}(x)$ seem also natural.

 To justify the decay assumptions on $\vec b_2(x)$, we use the result in \cite{20}. Consider
the IVP
\begin{equation}
\left\{
\begin{array}{l}
\p_t  u =i\p^2_{x_1x_2} u + \vec b_2(x)\cdot\nabla_x \bar u,\;\;\;\;x\in\BbbR^2,\;t>0,\\
u(x,0)=u_0(x)
\end{array}
\right.\label{eq1.7}
\end{equation}
with $\vec b_2=(i,0)$. It was shown in \cite{20} that the IVP (\ref{eq1.7}) is wellposed, however its
solutions gain only $1/4$ of derivative compared with $u_0$ instead of the expected $1/2$, i.e.
(\ref{eq1.5}) holds with $f=0$ and $1/4$ instead of $1/2$. Moreover, if we replace
$\p^2_{x_1x_2}$ by $\p^2_{x_1}+\p^2_{x_2}$ the expected gain of $1/2$ derivatives is obtained. So one has that in
the non-elliptic case, decay assumptions on $\vec b_2$ are also necessary to obtain (\ref{eq1.5}). The main
result we obtain concerning (\ref{eq1.4}) is Theorem \ref{th5.1.2} in Subsection \ref{subseq5.1}.

 When considering the variable coefficient equation in (\ref{eq1.1}) one should study of the bicharacteristic flow, i.e.
solutions
$(X(s;x,\xi),\Xi(s;x,\xi))$ of the system
\begin{equation}
\left\{
\begin{array}{l}
\displaystyle\frac{d\,}{ds}X_j(s;x,\xi)=2\sum_{k=1}^na_{jk}(X(s;x,\xi))\Xi_k(s;x,\xi),\\
\displaystyle\frac{d\,}{ds}\Xi_j(s;x,\xi)=-\sum_{k,l=1}^n\p_{x_j}a_{kl}(X(s;x,\xi))\Xi_k(s;x,\xi)\Xi_l(s;x,\xi),\\
\left(X(0;x,\xi\right),\Xi(0;x,\xi))=(x,\xi).
\end{array}
\right.\label{eq1.8}
\end{equation}

 In the constant coefficient case, ${\cal L}= -a_{jk}^0\p_{x_j}\p_{x_k}$, one has that the
bicharacteristic flow is
$$
\left(X(s;x,\xi),\Xi(s;x,\xi)\right)=(x+2sA^0\xi, \xi),\;\;\;\;\;\;\;A^0=(a_{jk}^0)_{j,k=1,..,n}.
$$

 For  ${\cal L}=-\Delta$ one gets $\left(X(s;x,\xi),\Xi(s;x,\xi)\right)=(x+2s\xi,\xi)$ and the condition in
(\ref{eq1.6}) can be seen as an
integrability
one along the bicharacteristics.
As we will see in this work (Sections 4-5), roughly speaking the operator $K$ whose symbol is
\begin{equation}
k(x,\xi)=\exp\,\left(\int_{-\infty}^0b(X(s;x,\xi),\Xi(s;x,\xi))ds\right),\;\;\;\mbox{with}\;\;b(x,\xi)=-\mbox{Im}\,
\vec b_1(x)\cdot
\xi,\label{eq1.9}
\end{equation}
 will play the role of the \lq\lq integrating factor" introduced in \cite{11}. 
Such constructions were also previously used
in the works \cite{4} and \cite{8}. The commutator  term $i[K{\cal L}-{\cal L} K]$, used to cancel the term $\vec
b_1(x)\cdot
\nabla$,  corresponds to differentiation 
of $K$ along
the bicharacteristic flow. Unfortunately the symbol in (\ref{eq1.9}) is not in a standard class. 
It satisfies
\begin{equation}
\left|\p_x^{\alpha}\p_{\xi}^{\beta} k(x,\xi)\right|\leq c_{\alpha 
\beta}\langle x\rangle^{|\beta|}\langle \xi\rangle^{-|\beta|}.\label{eq1.10}
\end{equation}
We observe that in  the particular constant coefficient elliptic case, ${\cal L}=\Delta$ we have
$$
\begin{array}{rcl}
k(x,\xi)&=&\displaystyle
 \exp\,\left(\int_{-\infty}^0 -\mbox{Im}\,\vec b_1(x+s\xi)\cdot\xi ds\right)\\
&=&\displaystyle
\exp\,\left(\int_{-\infty}^0 -\mbox{Im}\,\vec
b_1(x+s\hat{\xi})\cdot\hat{\xi}ds\right),\;\;\;\;\;\;\;\hat{\xi}=\xi/|\xi|,
\end{array}
$$ 
which is  related to Mizohata's condition in (\ref{eq1.6}).

 In the case where ${\cal L}(x)$ is elliptic, the class described in (\ref{eq1.10}) was introduced and
studied by W. Craig, T. Kappeler and W. 
Strauss in \cite{4}. However, it should be pointed out that 
in the non-elliptic case, i.e. ${\cal L}(x)$ is just
nondegenerate,
the geometric assumption (4.2) in \cite {4} is not satisfied by symbols of interest -see subsection \ref{subseq3.1}
in Section \ref{sec3}. 
 Moreover, we observe that the Hamiltonian
$\displaystyle h_2(x,\xi)=\sum_{j,k}a_{jk}(x)\xi_k\xi_j $ is preserved under the flow. 
Then one of the main differences of the flows
considered here with those associated to elliptic operators is that ellipticity gives the \it a priori \rm estimate 
$$
\nu^{-2}|\xi_0|^2\leq \left|\Xi(s;x_0,\xi_0)\right|^2\leq \nu^2|\xi_0|^2, 
$$
which guarantees that the solutions of the system (\ref{eq1.8}) are globally defined.
  
 We will assume that the bicharacteristic flow is non-trapping, i.e. for each
$(x_0,\xi_0)\in\BbbR^n\times(\BbbR^n-\{0\})$  and for each $\mu>0$ there exists $s_0>0$ such
that
\begin{equation}
|X(s_0;x_0,\xi_0)|\geq \mu.\label{eq1.11}
\end{equation}
The non-trapping condition appears naturally, since Ichinose \cite{13} showed that a
necessary condition for the well-posedness in $L^2$ of (\ref{eq1.1}) with ${\cal L}$ elliptic, $\vec
b_2\equiv0$, $c_1\equiv0,$ $ c_2\equiv0, $ and $P\equiv0,$ is that the analog of (\ref{eq1.6})
must hold in this case, with the integration taking place along the bicharacterisitics.
 The non-trapping condition also is
essential in the works of
\cite{4} and \cite{8}. In fact even when also $\vec b_1\equiv0$
and $f\equiv 0$, Doi showed in \cite{9} that it is necessary for (\ref{eq1.5}) to hold.

 In the ultra-hyperbolic case (i.e. with a merely non-degenerate matrix $A$), under appropriate decay assumptions and
asymptotic behavior as
$|x|\to
\infty$ 
 on the coefficients $a_{jk}(x)$, we shall prove that the bicharacteristic flow is
globally defined and \lq\lq uniformly non-trapping". Moreover in order  to keep the
structure of the conjugate first order terms,
 so that after applying the operator $K$ we
can obtain energy estimates, we need the symbol of $K$ to be even -see Definition \ref{def5.2.2} (iv) in
Section \ref{sec5}. Therefore we have to study carefully the bicharacteristic 
for backward and forward time. In particular when looking at the forward bicharacteristic
the more delicate part is when it is not outgoing -see Theorem \ref{th4.1.5} of Section \ref{sec4}. In that region we prove that  outside a bounded
ball, in the
$x$ variable, it behaves in dyadic annuli as  the free flow -see Theorem \ref{th4.1.5} in Section \ref{sec4} for a
precise statement. 

 As in \cite{18}, \cite{22}, the proof of the nonlinear results relies on two kinds
of linear estimates. 
The first one is concerned with the smoothing effect described in (\ref{eq1.5}) for solutions of the IVP (\ref{eq1.4})
\begin{equation}
 \begin{array}{l}
\displaystyle\int_0^{T}\int_{\BbbR^n}\left|J^{s+1/2}u(x,t)\right|^2\langle x\rangle^{-N}dxdt\\
\qquad\displaystyle\leq c(1+T)\sup_{0\leq t\leq T}\|u(t)\|^2_{H^s}
+c\int_0^{T}\int_{\BbbR^n}\left|f(x,t)\right|^2dxdt,
\end{array}\label{eq1.12}
 \end{equation}
and
$$
 \begin{array}{l}
\displaystyle\int_0^{T}\int_{\BbbR^n}\left|J^{s+1/2}u(x,t)\right|^2\langle x\rangle^{-N}dxdt
 \\
 \qquad\displaystyle\leq c(1+T)\sup_{0\leq t\leq T}\|u(t)\|^2_{H^s}+\int_0^{T}\int_{\BbbR^n}\left|J^{s-1/2}f(x,t)\right|^2\langle x\rangle^{N}dxdt
\end{array}
 $$
for $N>1$.

 The second kind is related with the local well-posedness in $L^2$ (and in $H^s$) of the IVP
(\ref{eq1.4}). To establish this result we  follow an indirect approach. First we truncate at infinity
the operator ${\cal L}(x)$ using $\theta\in C^{\infty}_0(\BbbR^n)$ 
with $\theta(x)=1,\;|x|\leq 1$, and $\theta(x)=0,\;
|x|\geq2$.
For $R>0$ we define$$
{\cal L}^R(x)=\theta(x/R){\cal L}(x)+(1-\theta(x/R)){\cal L}^0,
$$
where ${\cal L}^0=-a^0_{jk}\p_{x_j}\p_{x_k}$, $A^0=(a^0_{jk})_{j,k=1,..,n}$ is a (constant) matrix,
with the decay assumption 
$a_{jk}(x)-a_{jk}^0\in{\cal S}(\BbbR^n)$, $j,k=1,..,n$, (although we will work  in the ${\cal
S}(\BbbR^n)$ class, it will be clear from our proofs that the same results hold if we just
assume that the corresponding estimate holds for  a finite number of seminorms in (\ref{eq2.1.2}) Section
2). 
 Thus,
$$
{\cal L}(x)={\cal L}^R(x)+{\cal E}^R(x).
$$
 For $R$ large enough we consider the bicharacteristic flow $(X^R(s;x,\xi),\Xi^R(s;x,\xi))$ 
associated to the operator ${\cal L}^R(x)$ and  the 
corresponding integrating factor $K^R$, i.e. the operator with symbol as in (\ref{eq1.9}) but
evaluated in the 
 bicharacteristic flow $(X^R(s;x,\xi),\Xi^R(s;x,\xi))$. To obtain the $L^2$ local
well-posedness
of the IVP (\ref{eq1.4}) 
 we show that there exists $N_0$ depending only on the dimension such that for any   $M\in\BbbZ^+$ there exists $ R_0=R_0(M)$ such that for $R\geq R_0$
$$
\begin{array}{l}
\displaystyle\sup_{0\leq t\leq T}\left\|K^R u(t)\right\|^2_{L^2} \leq c R^{N_0}\|u(0)\|_{L^2}^2\\
\qquad\displaystyle+R^{-M}\int_0^T\int_{\BbbR^n}\left|J^{1/2}u\right|^2\langle x\rangle^{-N}dxdt+
cT R^{N_0} \sup_{0\leq t\leq T}\|u(t)\|^2_{L^2}.
\end{array}
$$
Next, we deduce  several estimates concerning the operator $K^R$. In particular, for
$E^R=I-\tilde K^R (K^R)^*$,
where the symbol of
$\tilde K^R$ differs from that of $K^R$ only in the sign of the exponent, and $
(K^R)^*$ is the adjoint of $K^R$, which allows us to treat $E^R u(t)$ as an error term -see Lemma \ref{le5.2.8} in Section \ref{sec5}. Collecting
these results  we get that
$$
\begin{array}{l}
\displaystyle\sup_{0\leq t\leq T}\left\| u(t\right\|^2_{L^2} \leq R^{N_0}\left\|u(0\right)\|_{L^2}^2\\
\qquad\displaystyle+R^{-M}\int_0^T\int_{\BbbR^n}\left|J^{1/2}u\right|^2\langle x\rangle^{-N}dxdt+
cT \sup_{0\leq t\leq T}\|u(t)\|^2_{L^2},
\end{array}
$$
which combined with (\ref{eq1.12}) yields the desired estimate, i.e. the local well-posedness in $L^2$ of the IVP (\ref{eq1.4}) 
for $T$ sufficiently small
$$
\sup_{0\leq t\leq T}\| u(t)\|^2_{L^2}
\leq c(T)\left(\|u_0\|^2_{L^2}+\int_0^T \|f(t)\|^2_{L^2}dt\right).
$$

  The smoothing effect and local well-posedness in $H^s$ of (\ref{eq1.4}) in the case where ${\cal L}$ has elliptic variable
coefficients will be proven in Section \ref{sec2} -see Lemma \ref{le 2.2.2} and Theorem \ref{th 2.3.1} of that section. This builds
on S. Doi's pseudo-differential method in \cite{7}, \cite{8}  and on H. Chihara's diagonalization
method for systems \cite{2} and uses only classical pseudo-differential operators.

 The diagonalization method cannot be used when ${\cal L}$ is ultrahyperbolic.
When ${\cal L}$ has constant coefficients it is possible to cancel the loss of derivatives using a
pseudo-differential  transformation  which falls under the scope of
Calder\'on-Vaillancourt's theorem -see \cite{22}, but this does not seem  to extend to the variable coefficient
case. Instead one is led to study a new class of  symbols and this is done in
Section \ref{sec3}. As we already mentioned the corresponding operators  in the elliptic case were studied
in
\cite{4}. 

A typical example of the symbols that need to be considered  (take $n=2$ for
simplicity) is 
$$a\bigl(x\cdot (\xi_2, \,-\xi_1)/|\xi|\bigl)\chi(\xi),
$$
where $\chi\in  C^\infty$, $\chi\equiv1$ for $ |\xi|\geq2$, $\chi\equiv0$ for $ |\xi|<1$, and
$a\in  C_0^\infty(\BbbR)$, in the elliptic case; and
$$a\bigl(x\cdot (\xi_2, \,\xi_1)/|\xi|\bigl)\chi(\xi),
$$
in the ultrahyperbolic one with $\chi$ and $a$ as before.
 
 As is explained in \cite{21}, the operators of the elliptic case are easily reduced to
classical pseudo-differential operators by expressing $a$ in terms of its Fourier transform and,
given  $\tau\in\BbbR,$ using the invertible change of variable $\xi\rightarrow\, \xi+\tau(\xi_2,
\,-\xi_1)/|\xi|$.  In the ultrahyperbolic setting this approach fails since the corresponding
mappings are not invertible, and hence the theory, in particular the $L^2$ boundedness, is more
delicate -see Theorem \ref{th3.2.1} in Section \ref{sec3}. The proofs of the rest of the results concerning the
calculus of the operators arising from these symbols-Theorem \ref{th3.3.4}, Theorem \ref{th3.3.5} and Theorem \ref{th3.3.6} of the
same section, are reduced after
some manipulations to the $L^2$ boundedness.

 In Section \ref{sec4} we study the bicharacteristic flow in the ultrahyperbolic case for ${\cal L}(x)$ and 
 its truncated version ${\cal L}^R(x)$. There
we shall deduce several estimates to be used in establishing the smoothing effect and the local well-posedness
of (\ref{eq1.4}) with ${\cal L}(x)$ non-elliptic which will be given in Section \ref{sec5}. This also relies on the calculus
of Section \ref{sec3}.

 Finally, the smoothing effect in (\ref{eq1.4}) is used to solve (\ref{eq1.1}) in Section \ref{sec6}. Solutions of (\ref{eq1.1}) are fixed points
of an integral mapping which is a contraction on a suitable function space in a small time 
interval, so Banach's contraction mapping principle
applies.

The results in the elliptic case, i.e. those in Section \ref{sec2}, are due to C. Rolvung, and appear in
his PhD dissertation \cite{27}. The results in the ultrahyperbolic case, for ${\cal L}$ a $
C_0^\infty$ perturbation of a constant coefficient operator ${\cal L}_0$ also appear in \cite{27}.
 
 \newpage 

\section{THE LINEAR ELLIPTIC EQUATION}\label{sec2}

The local well posedness and smoothing 
 effect for linear elliptic equations are considered in this 
 section. This builds on S. Doi`s method involving classical
 pseudo-differential operators and the sharp G\aa rding inequality \cite{8}
 as well as on a diagonalization as in \cite{2}.

\subsection{Pseudo-differential Operators}\label{subseq2.1}

  First we will recall some results from the theory of pseudo-differential
operators.

 The class $S^m=S^m_{1,0}$ of classical symbols of order $m \in \BbbR$  is defined by
\begin{equation}
S^m=\left\{p(x,\xi)\in C^{\infty}(\BbbR^n\times \BbbR^n)\,:\,|p|_{S^m}^{(j)} <\infty,\;j\in \BbbN\right\}\label{eq2.1.1}
\end{equation}
where
\begin{equation}
|p|_{S^m}^{(j)} = \sup \left\{\,\left\|\langle\xi\rangle^{-m+|\alpha|}\,\p_{\xi}^{\alpha}\p_x^{\beta}p(\cdot,\cdot)
\right\|_{L^{\infty}(\BbbR^n\times \BbbR^n)}
\,:\,|\alpha+\beta|\leq j\right\}\label{eq2.1.2}
\end{equation}
and $\langle\xi\rangle=(1+|\xi|^2)^{1/2}$.

 The pseudo-differential
operators $\Psi_p$ associated to the symbol $p\in S^m$ is defined by
\begin{equation}
\Psi_p f(x)= \int_{\BbbR^n}\,e^{ix\cdot\xi}p(x,\xi)\hat f(\xi)\frac{d\xi}{(2\pi)^{n/2}},\;\;\;\;\;\;\;
f\in{\cal S}(\BbbR^n).\label{eq2.1.3}
\end{equation}

 For example, a partial differential operator
$$
P=\sum_{|\alpha|\leq N}a_{\alpha}(x)\p_x^{\alpha},
$$
with $a_{\alpha}\in C^{\infty}_b(\BbbR^n)$ is a pseudo-differential operator $P=\Psi_p$ with
symbol
$$
p(x,\xi)=\sum_{|\alpha|\leq N}a_{\alpha}(x)(i\xi)^{\alpha}.
$$

 The fractional differentiation operator $J^s=\Psi_{\langle \xi\rangle^s}$ is also a 
pseudo-differential operator. The collection of symbol classes $S^m,\,m\in\BbbR$, is in some
cases closed under the division and square root operations. This is not the case for
polynomials
in $\xi$ and sometimes allows one to construct approximate inverses and square roots of
pseudo-differential operators.

 The following facts will be used throughout this work and the proofs can be found for example
in \cite{25}.

\begin{theorem}[Sobolev boundedness]\label{th 1.1}

Let $m\in\BbbR$, $p\in S^m$ and $s\in\BbbR$. Then
$\Psi_p$ extends to a bounded linear operator from
$H^{m+s}(\BbbR^n)$ to $H^{s}(\BbbR^n)$. Moreover, there
exist
$j=j(n;m;s)\in\BbbN$ and $c=c(n;m;s)$ such that
\begin{equation}
\|\Psi_p f\|_{H^s} \leq c\,|p|_{S^m}^{(j)}\,\| f\|_{H^{m+s}}.\label{eq2.1.4}
\end{equation}
Finally for $p\in {\cal S}^0$ and $\lambda
(|x |)=(1+|x|^2)^{-N/2}=\langle x\rangle^{-N}$, $N>1$ there exists $j=j(n,N)$ such that
\begin{equation}
 \int_{\BbbR^n}\, |\Psi_p f(x)|^2\,\lambda(|x|)\,dx
\leq c\,|p|_{S^0}^{(j)}\, \int_{\BbbR^n}\, | f(x)|^2\,\lambda(|x|)\,dx.\label{eq2.1.5}
\end{equation}
\end{theorem}
 
The proof of (\ref{eq2.1.5}) can be seen for example in \cite{22}.

\begin{theorem} [Symbolic calculus]\label{th 1.2}

Let $m_1,\,m_2\in \BbbR$, $p_1\in S^{m_1}$, $p_2\in S^{m_2}$. Then
 there exist $p_3\in S^{m_1+m_2-1}$, $p_4\in S^{m_1+m_2-2}$, and $p_5\in S^{m_1-1}$ such that
\begin{equation}
\begin{array}{rcl}
\Psi_{p_1}\Psi_{p_2}&=&\Psi_{p_1p_2}+\Psi_{p_3},\\
\Psi_{p_1}\Psi_{p_2}-\Psi_{p_2}\Psi_{p_1}&=&\Psi_{-i\{p_1,p_2\}}+\Psi_{p_4},\\
(\Psi_{p_1})^*&=&\Psi_{\bar p_1} + \Psi_{p_5}
\end{array}\label{eq2.1.6}
\end{equation}
where $\{p_1,p_2\}$ denotes the Poisson bracket, i.e.
$$
\{p_1,p_2\}=\sum_{j=1}^n\,\left(\p_{\xi_j}p_1\,\p_{x_j}p_2-\p_{x_j}p_1\,\p_{\xi_j}p_2\right),
$$
and such that for any $j\in\BbbN$ there exist $j'\in\BbbN$ and $c_1=c_1(n;m_1;m_2;j)$,
$c_2=c_2(n;m_1;j)$ such that
\begin{equation}
\begin{array}{rcl}
\displaystyle\left|p_3\right|_{S^{m_1+m_2-1}}^{(j)} +\left |p_4\right|_{S^{m_1+m_2-2}}^{(j)}&
\leq &\displaystyle  c_1\,\left|p_1\right|_{S^{m_1}}^{(j')}\,\left|p_2\right|_{S^{m_2}}^{(j')}\\ \\
\,\left|p_5\right|_{S^{m_1-1}}^{(j')} &\leq &c_2 \,\displaystyle\left|p_1\right|_{S^{m_1}}^{(j')}.
\end{array}\label{eq2.1.7}
\end{equation}
\end{theorem}

\begin{theorem}[Sharp G\aa rding inequality]\label{th 1.3}

Let $p\in S^{1}$ and suppose that there exists $R\geq 0$ such that $\mbox{Re} \,p(x,\xi)\geq 0\,$ for $\,\|\xi\|\geq R$. 
Then there exist $j=j(n)$ and $c=c(n;R)$ such that
\begin{equation}
\mbox{Re} \,\left\langle \Psi_{p}f;f\right\rangle_{L^2}\geq -c\,|p|_{S^1}^{(j)}\,\|f\|_{L^2}^2,\;\;\;\;\;
f\in{\cal S}(\BbbR^n).\label{eq2.1.8}
\end{equation}
\end{theorem}

 This result is due to L. H\"ormander \cite{12}.
 
\subsection{The Bicharacteristic Flow}\label{subseq2.2}
 
  The basic idea is to apply a pseudo-differential operator $K$ to the equation (\ref{eq1.4}) of Section
1
 in such a way that the commutator $i[K{\cal L}-{\cal L} K]$ cancels $\vec b_1(x)\cdot\nabla_x$.	
 It turns out that $i[K{\cal L}-{\cal L} K]$ corresponds to differentiation along the	
 bicharacteristic flow which will now be introduced.

	Let $A(x)=(a_{jk}(x))$ be a real, and symmetric $n\times n$ matrix of functions
 $a_{jk}\in C^{\infty}_b$. We will assume that	
\begin{equation}\label{eq2.2.1}	
\left|\nabla a_{jk}(x)\right|=o(|x|^{-1})\;\;\;\mbox{as}\;\;\;|x|\to\infty,\;\;\;j,k=1,..,n,	
\end{equation}	
and that $A(x)$ is positive definite, i.e.	
\begin{equation}\label{eq2.2.2}	
\exists\,\nu>0\;\;\;\;\;\forall\;x,\,\xi\in \BbbR^n\;\;\;\;\;	
\nu^{-1}|\xi|^2\leq\bigg |\sum_{j,k=1}^n\,a_{jk}(x)\xi_j\xi_k\bigg|\leq \nu|\xi|^2.	
\end{equation}
Let $h_2$ be the principal symbol of ${\cal L}=-\p_{x_j}a_{jk}(x)\p_{x_k}$, i.e.
\begin{equation}\label{eq2.2.3}
h_2(x,\xi)= \sum_{j,k=1}^n\,a_{jk}(x)\xi_j\xi_k.
\end{equation}
 The bicharacteristic flow is the flow of the Hamiltonian vector field
\begin{equation}\label{eq2.2.4}
H_{h_2}=\sum_{j=1}^n\,\left[\p_{\xi_j}h_2\cdot\p_{x_j}-\p_{x_j}h_2\cdot\p_{\xi_j}\right]
\end{equation}
and is denoted by $(X(s;x_0,\xi_0),\Xi(s;x_0,\xi_0))$, i.e.
\begin{equation}\label{2.2.5}
  \left\{
  \begin{array}{l}
\displaystyle \frac{d}{ds} X_j(s;x_0,\xi_0)=
2\sum_{k=1}^n a_{jk}\left(X(s;x_0,\xi_0)\right)\,\Xi_k(s;x_0,\xi_0),\\
\\
\displaystyle\frac{d}{ds} \Xi_j(s;x_0,\xi_0)=-
\sum_{k,l=1}^n \p_{x_j}a_{lk}\left(X(s;x_0,\xi_0)\right) 
\Xi_k(s;x_0,\xi_0) \Xi_l(s;x_0,\xi_0)
\end{array}
\right.
  \end{equation}
for $j=1,..,n$, with
\begin{equation}\label{eq2.2.6}
\left(X(0;x_0,\xi_0),\Xi(0;x_0,\xi_0)\right)=(x_0,\xi_0).
\end{equation}

 The bicharacteristic flow
exists in the time interval $s\in(-\delta,\delta)$ with 
$\delta=\delta(x_0,\xi_0)$, and $\delta(\cdot)$ depending continuously on $(x_0,\xi_0)$.

 The bicharacteristic flow preserves $h_2$, so ellipticity gives
\begin{equation}\label{eq2.2.7}
\nu^{-2}\left|\xi_0|^2\leq |\Xi(s;x_0,\xi_0\right|^2\leq \nu^2|\xi_0|^2,
\end{equation}
and hence $\delta=\infty$.

 It will be assumed  that the bicharacteristic flow is non-trapped
which means that the set 
$$
\left\{X(s;x_0,\xi_0)\,:\,s\in\BbbR\right\}
$$
is unbounded in $\BbbR^n$ for each $(x_0,\xi_0)\in\BbbR^n\times\BbbR^n-\{0\}$.
  
   Note that $h_2$ is homogeneous of degree $2$ in $\xi$ so that
\begin{equation}\label{eq2.2.8}
  \left\{
  \begin{array}{l}
X(s;x,t\xi)=X(ts;x,\xi),\\
\Xi(s;x,t\xi)=t\Xi(ts;x,\xi).
\end{array}
\right.
  \end{equation}

 The next result
shows that the Hamiltonian vector field is differentiation
 along the bicharacteristics.

\begin{lemma}\label{le 2.2.1}

 Let $\phi\in C^{\infty}(\BbbR^n\times\BbbR^n)$. Then
\begin{equation}\label{eq2.2.9}
(H_{h_2}\phi)(x,\xi)=\p_s \left[\phi(X(s;x,\xi),\Xi(s;x,\xi))\right]\left.\right|_{s=0}.
\end{equation}
\end{lemma}

 The following key lemma is due to S. Doi \cite{8} (Lemmas 2.3-2.5).

\begin{lemma}\label{le 2.2.2}

 Let $A(x)$ and its bicharacteristic flow satisfy the assumptions above.
Suppose $\lambda\in L^1\left([0,\infty))\cap C([0,\infty)\right)$ is strictly positive and nonincreasing.
Then there exist $c>0$ and a real symbol $p\in S^0$, both depending on $h_2$ and $\lambda$,
such that
\begin{equation}\label{eq2.210}
H_{h_2}p=\{h_2,p\}(x,\xi)\geq \lambda(|x|)\,|\xi|-c,\;\;\;\;\forall \,(x,\xi)\in \BbbR^n\times\BbbR^n.
\end{equation}
\end{lemma}
  
  An extension of this result to the case of invertible $A(x)$ will be given in Section \ref{sec5}, Lemma \ref{le5.1.1}.

\subsection{Linear Elliptic Smoothing Effects}\label{subseq2.3}

 In this subsection we consider the IVP associated to the  linear Schr\"odinger equation
\begin{equation}\label{eq2.3.1}
  \left\{
  \begin{array}{l}
\displaystyle\p_tu=-i\sum_{j,k=1}^n \p_{x_j}(a_{jk}(x)\p_{x_k}u)+\vec b_1(x)\cdot\nabla_x u+\vec 
b_2(x)\cdot\nabla_x \bar u
+c_1(x)u+c_2(x)\bar u + f(x,t),\\
\displaystyle u(x,0)=u_0(x),
\end{array}
\right.
  \end{equation}
where $A(x)= (a_{jk})_{j,k=1,..,n}$ satisfies (\ref{eq2.2.1})-(\ref{eq2.2.2}) and its bicharacteristic flow satisfies the
assumptions in the previous subsection, $\vec b_l=(b^1_l,..,b^n_l)\in\left(C^{\infty}_b\right)^n$,
$l=1,2$ and $ c_1, c_2\in C^{\infty}_b$.

 Combining the equation in (\ref{eq2.2.1}) and its complex conjugate we obtain a system in
$\vec w=(u,\,\bar u)^T$
\begin{equation}\label{eq2.3.2}
  \left\{
  \begin{array}{l}
\p_t \vec w=(iH+B+C) \vec w+\vec f,\\
\vec w(x,0)=\vec w_0(x),
\end{array}
\right.
  \end{equation}
where
$$
H=\begin{array}{l}
  \left(
  \begin{array}{ccc}
{\cal L}&0\\0&-{\cal L}
  \end{array}
  \right)
  ,\;\;
\;\;C=  \left(
  \begin{array}{ccc}
 c_{11}&c_{12}\\c_{21}&c_{22}
 \end{array}
  \right),
\end{array}
$$

$$
B=\left(
  \begin{array}{ccc} B_{11}&B_{12}\\B_{21}&B_{22} \end{array}
  \right),
$$
with 
$$
{\cal L}=-\sum_{j,k=1}^n\p_{x_j}\,(a_{jk}(x)\p_{x_k}),
$$
$$
B_{lm}u(x,t)=\sum_{j=1}^n b_{lmj}(x)\p_{x_j}u(x,t),\;\;\;b_{lmj}\in
C^{\infty}_b\;\;l,m=1,2,\;j=1,..,n,
$$
and
$$
c_{lm}=c_{lm}(x)\in C^{\infty}_b,\;\;l,m=1,2,\;\;\vec f(x,t)=(f(x,t),\bar f(x,t))^T.
$$

 The following well-posedness and smoothing results contain three parts depending 
on the regularity and the decay of the external force $f(x,t)$.

\begin {theorem}\label{th 2.3.1}

Let $\vec w_0=(u_0,\bar u_0)^T\in (H^s(\BbbR^n))^2$, $s\in \BbbR$. Assume that there
exist $N>1$ and a constant $c_0$ such that if $\lambda(|x|)= \langle \,x \,\rangle^{-N}$ then
$$
\left|\,\mbox{Im}\,b_{llj}(x)\right|\leq c_0\lambda(|x|),\;\;l,m=1,2,\;j=1,..,n,\;x\in\BbbR^n.
$$
Then

(a) If $\vec f\in (L^1([0,T]:(H^s(\BbbR^n)))^2$ then the IVP (\ref{eq2.3.2})
has a unique solution
$\vec w\in C([0,T]:(H^s(\BbbR^n))^2)$ satisfying
\begin{equation}\label{eq2.3.3}
\sup_{0\leq t\leq T}\|\vec w(t)\|_{(H^s)^2}\leq c_1e^{c_2T}\left(\|\vec w_0\|_{(H^s)^2}
+\int_0^T\|\vec f(t)\|_{(H^s)^2}dt\right).
\end{equation}

(b) If $\vec f\in (L^2([0,T]:(H^s(\BbbR^n)))^2$ then the IVP (\ref{eq2.3.2})
has a unique solution $\vec w\in
C([0,T]:(H^s(\BbbR^n))^2)$ satisfying 
\begin{equation}\label{eq2.3.4}
  \begin{array}{l}
\displaystyle\sup_{0\leq t\leq T}\|\vec w(t)\|^2_{(H^s)^2}
+\int_0^T\int_{\BbbR^n}|J^{s+1/2}\vec w(x,t)|^2\lambda(|x|)dxdt\\
\qquad\displaystyle\leq
c_1e^{c_2T}\left(\|\vec w_0\|^2_{(H^s)^2}  +\int_0^T\|\vec f(t)\|^2_{(H^s)^2}dt\right). 
\end{array}
  \end{equation}
  
(c) If $J^{s-1/2} \vec f\in (L^2(\BbbR^n\times
[0,T]:(\lambda(|x|)^{-1}dxdt))^2$ then the IVP (\ref{eq2.3.2}) 
has a unique solution $\vec w\in C([0,T]:(H^s(\BbbR^n))^2)$ satisfying 
\begin{equation}\label{eq2.3.5}
  \begin{array}{l}
\displaystyle\sup_{0\leq
t\leq T}\|\vec w(t)\|^2_{(H^s)^2} +\int_0^T\int_{\BbbR^n}|J^{s+1/2}\vec
w(x,t)|^2\lambda(|x|)dxdt\\ 
\qquad\displaystyle\leq c_1e^{c_2T}\left(\|\vec w_0\|^2_{(H^s)^2} +
\int_0^T\int_{\BbbR^n}|J^{s-1/2}\vec f(x,t)|^2 (\lambda(|x|))^{-1}dx dt\right). 
\end{array}
  \end{equation}
  
  Here $c_1$ depends on $n, s, \nu, (a_{jk})_{j,k=1,..,n}, c_0, (b_{lmj})_{l,m=1,2;j=1,..,n}$, and
 $c_2$ depends in addition on \linebreak $(c_{lm})_{l,m=1,2}$.

\end{theorem}
  
  \begin{corollary}\label{co2.3.2}

Let $s\in \BbbR$ and $u_0\in H^s(\BbbR^n$) and suppose $\lambda$
satisfies  the assumptions of Theorem \ref{th 2.3.1}. Then
 
 (a) If $f\in L^1([0,T]:H^s(\BbbR^n))$, then there exists a unique solution
 $u\in C([0,T]:H^s(\BbbR^n))$ of (\ref{eq2.3.1}) satisfying
 $$
 \sup_{0\leq t\leq T}\|u(\cdot,t)\|_{H^s}\leq c_1e^{c_2T}\left(\|u_0\|_{H^s}+\int_0^T\|
 f(\cdot,t)\|_{H^s}dt\right)
 $$
 
 (b) If $f\in L^2([0,T]:H^s(\BbbR^n))$, then there exists a unique solution
 $u\in C([0,T]:H^s(\BbbR^n))$ of (\ref{eq2.3.1}) satisfying
 $$
 \begin{array}{l}
\displaystyle \sup_{0\leq t\leq T}\|u(\cdot,t)\|^2_{H^s} +\int_0^T\int_{\BbbR^n}|J^{s+1/2}u(x,t)|^2
 \lambda(|x|)dxdt\\
 \qquad\displaystyle\leq c_1e^{c_2T}\left(\|u_0\|^2_{H^s}+\int_0^T\|
 f(\cdot,t)\|^2_{H^s}dt\right);
 \end{array}
 $$
 
 (c) If $J^{s-1/2}f\in L^2(\BbbR^n\times [0,T]:(\lambda(|x|))^{-1}dxdt)$, then there exists a unique solution
 $u\in C([0,T]:H^s(\BbbR^n))$ of (\ref{eq2.3.1}) satisfying
 $$
 \begin{array}{l}
 \displaystyle\sup_{0\leq t\leq T}\|u(\cdot,t)\|^2_{H^s}
 +\int_0^T\int_{\BbbR^n}|J^{s+1/2}u(x,t)|^2
 \lambda(|x|)dxdt\\
 \qquad\displaystyle\leq c_1e^{c_2T}\left(\|u_0\|^2_{H^s}+\int_0^T\int_{\BbbR^n}|J^{s-1/2}
 f(x,t)|^2(\lambda(|x|))^{-1}dxdt\right).
\end{array}
 $$ 
 
 Here $c_1$ depends on $n, s, \nu, (a_{jk})_{j,k=1,..,n}, c_0, 
 (b_{{\mbox{lmj}}})_{l,m=1,2;j=1,..,n}$, and
 $c_2$ depends in addition on \linebreak $(c_{lm})_{l,m=1,2}$.
 \end{corollary}
  
    The following a priori estimate is needed for the proof of Theorem \ref{th 2.3.1}.
  
 \begin {lemma}\label{le2.3.3} 
   
Let $s\in \BbbR$ and suppose $\lambda$ satisfies  the assumptions of
Theorem \ref{th 2.3.1}. Then there exist $c_1$ depending on $n, s$, $\nu,
(a_{jk})_{j,k=1,..,n}$, $c_0, 
$ and finitely many derivatives of $(b_{lmj})_{l,m=1,2;j=1,..,n}$, and $c_2$ depending
in addition on finitely many derivatives of $(c_{lm})_{l,m=1,2}$ such that for all
$\vec w\in (C([0,T]:H^{s+2}(\BbbR^n))\cap C^1([0,T]:H^{s+2}(\BbbR^n)))^2$ the following
four estimates hold : 
 
$$
\begin{array}{l}
   \displaystyle (i)\;\;\;\;\;\sup_{0\leq t\leq T}\|\vec w(t)\|_{(H^s)^2}\leq c_1e^{c_2T}
\Big(\|\vec w_0\|_{(H^s)^2}\\
\;\;\;\;\;\;\;\;\;\;\;\;\;\;\;\;\;\;\;\;
\displaystyle+\int_0^T\|(\partial_t-(iH+B+C))\vec w(\cdot,t)\|_{(H^s)^2}dt\Big),\\
(ii)\displaystyle \;\;\;\;\sup_{0\leq t\leq T}\|\vec w(t)\|_{(H^s)^2}\leq c_1e^{c_2T}
\Big(\|\vec w(\cdot,T)\|_{(H^s)^2}\\
\;\;\;\;\;\;\;\;\;\;\;\;\;\;\;\;\;\;\;\;
\displaystyle+\int_0^T\left\|(\partial_t+(iH+B+C)^*)\vec w(\cdot,t)\right\|_{(H^s)^2}dt\Big),\\
\displaystyle(iii)\;\sup_{0\leq t\leq T}\|\vec w(t)\|^2_{(H^s)^2}
+\int_0^T\int_{\BbbR^n}\left|J^{s+1/2}\vec w(x,t)\right|^2\lambda(|x|)dxdt\\
   \;\;\;\;\;\;\;\;\;\;\;\;\;\;\;\;\;\;\;\;
\displaystyle\leq
c_1e^{c_2T}\left(\|\vec w_0\|^2_{(H^s)^2}  + 
\int_0^T\left\|(\partial_t-(iH+B+C))\vec w(\cdot,t)\right\|^2_{(H^s)^2}dt\right),\\ 
\displaystyle(iv)\;\sup_{0\leq t\leq T}\|\vec w(t)\|^2_{(H^s)^2}
+\int_0^T\int_{\BbbR^n}\left|J^{s+1/2}\vec w(x,t)\right|^2\lambda(|x|)dxdt\\
   \;\;\;\;\;\;\;\;\;\;\;\;\;\;\;\;\;\;\;\;
\displaystyle\leq c_1e^{c_2T}
\|\vec w_0\|^2_{(H^s)^2}\\
    \;\;\;\;\;\;\;\;\;\;\;\;\;\;\;\;\;\;\;\;
\displaystyle+c_1e^{c_2T}
\int_0^T\int_{\BbbR^n}\left|J^{s-1/2}(\partial_t-(iH+B+C))\vec w(\cdot,t)\right|^2
(\lambda(|x|))^{-1}dxdt.
\end{array}
$$

\end{lemma}
    
      We observe that (ii) follows from (i) by applying (i), with $iH+B+C$ replaced by
$(iH+B+C)^*$, to $\vec w(\cdot,T-t)$, so it suffices to prove (i), (iii), and (iv)
of Lemma \ref{le2.3.3}. The idea of the proof is to apply transformations $\Lambda$ and $K$ 
to the system. $\Lambda$ will diagonalize $B$ and essentially transforms the system into two single
 equations where the pseudo-differential calculus applies. The idea of this 
 diagonalization came from the work of H. Chihara \cite{2}. $K$ will eliminate the loss of derivatives of the
 first order terms. This idea is due to S. Doi \cite{7}, \cite{8} and S. Mizohata \cite{26}.

{\sl Proof of Lemma \ref{le2.3.3} \rm}

Let $\vec w\in (C([0,T]:H^{s+2}(\BbbR^n))\cap C^1([0,T]:H^{s}(\BbbR^n)))^2$. 
 Set
$$
h_2(x,\xi)=\sum_{j,k=1}^n\,a_{jk}(x)\xi_j\xi_k,\;\;\;
h_1(x,\xi)=i\,\sum_{j,k=1}^n \p_{x_j}a_{jk}(x)\xi_k,
$$
so that ${\cal L}=\Psi _{h_1+h_2}$. Since $(a_{jk}(x))$ is positive definite there exist $
c=c(a_{jk};\,\p_{x_j}a_{jk})$ and $\,R=R(a_{jk};\,\p_{x_j}a_{jk})$ such that
$$
|h_1(x,\xi)+h_2(x,\xi)|\geq c|\xi|^2,\;\;\;\;\;\;\;\forall \,|\xi|\geq R.
$$
 Choosing $\phi\in C^{\infty}_0(\BbbR^n)$ with $\phi(y)=1$ if $|y|\leq R$ and
$\phi(y)=0$ if $|y|\geq 2R$ we define
$$
\tilde h(x,\xi)=(h_1(x,\xi)+h_2(x,\xi))^{-1}(1-\phi(\xi))\;\;
\;\mbox{and}\;\;\;\tilde{\cal L}=\Psi_{\tilde h}.
$$
So $\tilde h\in S^{-2}$ and 
$$
\tilde {\cal L} {\cal L} = I + \Psi_{r_1},
$$ 
where $r_1=r_1((a_{jk})_{j,k=1,..,n})\in S^{-1}$ by the symbolic calculus in Theorem \ref{th 1.2}
in the sense that for any $l\in\BbbN$, $|r_1|^l_{S^{-1}}$ depends on the ellipticity constant 
$\mu$ and on finitely many derivatives of the $a_{jk}$'s.
      
       We define $B_{diag}$, $B_{anti}$ with symbols in $(S^1)^{2\times2},\,S_{12},\,S_{21}$
with symbols in $S^{-1}$ and $S$ with symbol in $(S^{-1})^{2\times2}$ by
$$ 
\begin{array}{l}
        B_{diag}=
        \left(
        \begin{array}{cc}
        B_{11}&0\\
        0&B_{22}
        \end{array}
        \right),\;\;\;\;B_{anti}=
        \left(
        \begin{array}{cc}
0&B_{12}\\
        B_{21}&0
        \end{array}
        \right),
\end{array}
$$
$$
S_{12}=\frac{1}{2}\,iB_{12}\tilde{\cal L},\;\;\;S_{21}=-\frac{1}{2}\,iB_{21}\tilde{\cal L},
$$
$$
S=\left(
        \begin{array}{cc}
        0&S_{12}
        \\S_{21}&0
        \end{array}
        \right),
$$
and the diagonalizing transform $\Lambda$ with symbol in $(S^0)^{2\times 2}$ by
$$
\Lambda = I-S.
$$
Then $S=S(n;\nu;(a_{jk});(b_{lmj}))$ in the sense above. Letting 
$$
\vec f=
(\p_t-(iH+B+C))\vec w,
$$ 
and applying $\Lambda$, one obtains
\begin{equation}\label{eq2.3.6}
\p_t \Lambda \vec w =i\Lambda H \vec w +\Lambda B \vec w + \Lambda C \vec w+\Lambda\vec f.
 \end{equation}

 We shall show that the system (\ref{eq2.3.6}) is diagonalized modulo operators with symbols in $S^0$. So wewrite
\begin{equation}\label{eq2.3.7}
        \begin{array}{rcl}
i\Lambda H +\Lambda B &=& iH \Lambda  + (i\Lambda H - i H \Lambda) + B - SB\\
        &=&iH \Lambda  + (i\Lambda H - i H \Lambda) + B_{diag} +B_{anti} - SB\\
&=&(iH \Lambda  + B_{diag}\Lambda) + (B_{anti}+ i\Lambda H - i H \Lambda) + (B_{diag}S-S B).
\end{array}
        \end{equation}
 Since the operator in the first parenthesis is diagonalized and the operator in the last
parenthesis has order $0$, it suffices to consider only the operator in the second parenthesis
in (\ref{eq2.3.7}).
 Thus,
$$
\begin{array}{l}B_{anti}+ i\Lambda H - i H \Lambda = B_{anti}+ i H S - i S H=
\\
\\
\qquad=
        \left(
        \begin{array}{cc}
         0&B_{12}\\
         B_{21}&0
         \end{array}
         \right)
          + 
           \left(
        \begin{array}{cc}
{\cal L}&0\\
0&-{\cal L}
 \end{array}
         \right)
\left(
        \begin{array}{cc}
 0&-\frac{1}{2}B_{12}\tilde {\cal L}
 \\ \frac{1}{2}B_{21}\tilde{\cal L}&0
 \end{array}
         \right)
\\
\\
\qquad-
\left(
        \begin{array}{cc}
0&-\frac{1}{2}B_{12}\tilde {\cal L}\\
 \frac{1}{2}B_{21}\tilde {\cal L}&0
 \end{array}
         \right)
\left(
        \begin{array}{cc}
 {\cal L}&0\\
 0&-{\cal L}
 \end{array}
         \right)
\\
\\
\qquad=\left(
        \begin{array}{cc}
0&B_{12}-\frac{1}{2}{\cal L} B_{12}\tilde {\cal L} -\frac{1}{2}B_{12}\tilde{ \cal L}{\cal L}\\
B_{21}-\frac{1}{2}{\cal L} B_{21}\tilde {\cal L} -\frac{1}{2}B_{21}\tilde {\cal L}{\cal L}&0
\end{array}
         \right).
\end{array}
$$
 We observe that ${\cal L} B_{12}\tilde {\cal L}=B_{12}\tilde {\cal L}{\cal L}+\Psi_{r_2}$,
where $r_2\in S^0$ and
$$
B_{12}-\frac{1}{2}{\cal L} B_{12}\tilde {\cal L} -\frac{1}{2}B_{12}\tilde {\cal L}{\cal L}=
B_{12}-B_{12}\tilde{\cal L}{\cal L}-\frac{1}{2}\Psi_{r_2}=-B_{12}\Psi_{r_1}
-\frac{1}{2}\Psi_{r_2}.
$$
A similar calculus argument handles the term involving $B_{21}$. Therefore, we have that 
 $ B_{anti}+ i\Lambda H - i H \Lambda $ has order zero, which allows to conclude that
\begin{equation}\label{eq2.3.8}
\p_t\Lambda \vec w=i H\Lambda \vec w + B_{diag}\Lambda\vec w + \Psi_{r_3}\vec w +\Lambda
\vec f,
\end{equation}
where $r_3\in(S^0)^{2\times 2} $ and 
$$
r_3=r_3(n;\nu;(a_{jk});(b_{lmj});(c_{lm})),
$$ 
in the sense that for any $j_0\in\BbbZ^+$, $|r_3|^{(j_0)}_{(S^0)^{2\times 2}}$ depends
on $n, \nu$ and finitely many of the derivatives of $a_{jk}$, $b_{lmj}$ and $c_{lm}$.

 By Lemma \ref{le 2.2.2} there exists a real-valued $p\in S^0$ and $C>0$, both depending on
$(a_{jk})$ and $c_0$, such that
\begin{equation}\label{eq2.3.9}
\{h_2(x,\xi);p(x,\xi)\}\geq C'c_0\,\lambda(|x|)|\xi|-C,\;\;\;\;\forall x, \xi\in\BbbR^n,
\end{equation}
with $C'=C'(n)$ to be determined. Let
$$
k(x,\xi)=
         \left(
         \begin{array}{cc}
          e^{p(x,\xi)}(1+(h_2(x,\xi))^2)^{s/4}&0\\
0&e^{-p(x,\xi)}(1+(h_2(x,\xi))^2)^{s/4}
          \end{array}
          \right)
$$ 
and $K=\Psi_k$. Note that $e^{\pm p(x,\xi)}(1+(h_2(x,\xi))^2)^{s/4}\in S^s$ and is elliptic,
since $p\in S^0$ is real and
$$
\left(1+\nu^{-2}|\xi|^4\right)^{1/4}\leq\left(1+(h_2(x,\xi))^2\right)^{1/4}\leq \left(1+\nu^{2}|\xi|^4\right)^{1/4}
$$
where $\nu$ is the ellipticity constant of $(a_{jk})$.

The norm $N$ on $(H^s)^2$ is defined by
\begin{equation}\label{eq2.3.10}
(N(v))^2=\|K\Lambda v\|_{L^2}^2+\|v\|^2_{H^{s-1}}.
\end{equation}

 It will be shown that $N$ is equivalent to the standard $H^s$-norm.

Let 
$$ 
\tilde k(x,\xi)=
           \left(
         \begin{array}{cc}
          e^{-p(x,\xi)}(1+(h_2(x,\xi))^2)^{-s/4}&0\\
0&e^{p(x,\xi)}(1+(h_2(x,\xi))^2)^{-s/4}
           \end{array}
          \right)
$$
and $\tilde K=\Psi_{\tilde k}$. Then $\tilde k\in (S^{-s})^{2\times 2}$ and by the symbolic calculus in 
Theorem \ref{th 1.2}
$$
\tilde K K
=I+\Psi_{r_4},
$$
for some $r_4\in (S^{-1})^{2\times 2}$, where $|r_4|^{(j_0)}_{(S^{-1})^{2\times 2}}$ 
depends on $\nu,\,
(a_{jk})$ and $c_0$  for each $j_0\in\BbbN$. Therefore
$$
\tilde K K\Lambda=(I+\Psi_{r_4})(I-S)=I-(S-\Psi_{r_4}+\Psi_{r_4}S),
$$
where $S-\Psi_{r_4}+\Psi_{r_4}S$ has order $-1$. By the Sobolev boundedness (Theorem \ref{th 1.1})
$$
\begin{array}{rcl}
\displaystyle\|v\|^2_{H^s}&\leq& \displaystyle2\big\|\tilde K K \Lambda v\big\|^2_{H^s}+2\left\|(S-\Psi_{r_4}+\Psi_{r_4}S)v\right\|^2_{H^s}
\\
\displaystyle&\leq& \displaystyle c\left(\left\|K \Lambda v\right\|^2_{L^2} + \|v\|_{H^{s-1}}^2\right)\leq c \|v\|_{H^{s}}^2,
\end{array}
$$
for a sufficient large constant $c=c(n,\,s,\,\nu,\,(a_{jk}),\,c_0,\,(b_{lmj}))$ independent of $v$.
This shows the equivalence of the norms.

 Next, we shall estimate the norm $N$ to establish the inequalities (i), (iii) and (iv) in Lemma
\ref{le2.3.3} In the following $c_j$ will denote a constant depending on $n,\,s,\,\nu,\,(a_{jk})$, $c_0$
and finitely
many derivatives of $(b_{lmj})$ and $(c_{lm})$.

 To estimate the second term of $(N(v))^2$ in (\ref{eq2.3.10}) we write
$$
\begin{array}{rcl}
\displaystyle\p_t\|\vec w\|^2_{H^{s-1}}&=&\p_t\langle J^{s-1}\vec w, J^{s-1}\vec w\rangle_{L^2}=
2\,\mbox{Re}\, \langle J^{s-1}\p_t \vec w, J^{s-1}\vec w\rangle_{L^2}
\\
&=&\displaystyle2\,\mbox{Re}\, \langle J^{s-1}(iH\vec w+B\vec w+C\vec w +\vec f), J^{s-1}\vec w\rangle_{L^2}\\
&=&\displaystyle2\,\mbox{Re}\, \langle i H J^{s-1}\vec w, J^{s-1}\vec w\rangle_{L^2}\\
&+&\displaystyle2\,\mbox{Re}\, \langle (i[J^{s-1}H-HJ^{s-1}]+J^{s-1}B+J^{s-1}C)\vec w, J^{s-1}\vec w\rangle_{L^2}\\
&+&\displaystyle2\,\mbox{Re}\, \langle J^{s-1}f, J^{s-1}\vec w\rangle_{L^2}\\
&\leq&\displaystyle c_1\left(N(\vec w)\right)^2+c_2\,\mbox{min}\,\left\{N(\vec f)N(\vec w);\langle (\lambda(|x|))^{-1} J^{s-1/2}\vec f,
J^{s-1/2}\vec f\rangle_{L^2}\right\}.
\end{array}
$$
          
           Above we have used that $H$ is self-adjoint and diagonal, that $N$ and $H^s$-norm are equivalent, 
and that $\lambda$ is bounded above.

 For the first term of $(N(\vec w))^2$ in (\ref{eq2.3.10}) we write 
\begin{equation}\label{eq2.3.11}
           \begin{array}{rcl}
\p_t\|K\Lambda \vec w\|^2_{L^2}&=&\displaystyle\p_t\langle K\Lambda \vec w, K\Lambda \vec w\rangle_{L^2}=
2\,\mbox{Re}\, \langle \p_t K\Lambda \vec w, K\Lambda \vec w\rangle_{L^2}\\
&=&\displaystyle2\,\mbox{Re}\, \langle K \p_t \Lambda \vec w, K\Lambda \vec w\rangle_{L^2}\\
&=&\displaystyle2\,\mbox{Re}\, \langle K(iH+B_{diag})\Lambda \vec w+\Psi_{r_3}\vec w + \Lambda \vec f, K\Lambda \vec
w\rangle_{L^2}\\ 
          & =&\displaystyle2\,\mbox{Re}\, \langle K(iH+B_{diag})\Lambda\vec w, K\Lambda \vec w\rangle_{L^2}+2\,\mbox{Re}\,
\langle K\Psi_{r_3}\vec w, K\Lambda \vec w\rangle_{L^2}\\
&+&\displaystyle2\,\mbox{Re}\, \langle K \Lambda \vec f,K\Lambda \vec w\rangle_{L^2}\\
&\leq&\displaystyle 2\,\mbox{Re}\, \langle K(iH+B_{diag})\Lambda\vec w, K\Lambda \vec w\rangle_{L^2} +c(N(\vec w))^2+
2\,\mbox{Re}\, \langle K \Lambda \vec f,K\Lambda \vec w\rangle_{L^2}\\
&=&\displaystyle I+c(N(\vec w))^2+III,
\end{array}
\end{equation}
since $\left\|K\Psi_{r_3}\vec w\right\|_{L^2}\leq c\|\vec w\|_{H^s}\leq c(N(\vec w))$ by Sobolev boundedness and norm
equivalence.

 We should consider the terms $I$ and $III$ separately.
First we have
$$
\begin{array}{rcl}
I&=&\displaystyle2\,\mbox{Re}\, \langle K(iH+B_{diag})\Lambda\vec w, K\Lambda \vec w\rangle_{L^2}\\
&=&\displaystyle2\,\mbox{Re}\, \langle (iH+B_{diag})K\Lambda\vec w, K\Lambda \vec w\rangle_{L^2}
+2\,\mbox{Re}\, \langle i[KH-HK]\Lambda\vec w, K\Lambda \vec w\rangle_{L^2}\\
&+&\displaystyle2\,\mbox{Re}\, \langle [KB_{diag}-B_{diag}K]\Lambda\vec w, K\Lambda \vec w\rangle_{L^2}
\\
&\leq& \displaystyle2\,\mbox{Re}\, \langle (B_{diag}K+i[KH-HK])\Lambda\vec w, K\Lambda \vec w\rangle_{L^2}
+c(N(\vec w))^2,
\end{array}
$$
since $H$ is self-adjoint and $[KB_{diag}-B_{diag}K]$ is a commutator of diagonal
matrices and therefore has order $s$. Using the commutator formula of the symbolic calculus
on the diagonal matrices $K$ and $H$ it follows that
$$
i[KH-HK]=\Psi_q+\Psi_{r_5},
$$
where $r_5\in S^s$ and $q\in S^{s+1}$ is given by
$$
q=
           \left(
           \begin{array}{cc}
            \{e^p(1+h_2^2)^{s/4};h_2+h_1\}&0\\
0&\{e^{-p}(1+h_2^2)^{s/4};-h_2-h_1\}
\end{array}
            \right).
$$

 Since 
$$
\{e^p(1+h_2^2)^{s/4};h_2+h_1\}=-\{h_2;p\}e^p(1+h_2^2)^{s/4}+\{e^p(1+h_2^2)^{s/4}\,;\,h_1\},
$$
where the last term is in $S^s$, it follows that
$$
q=
           \left(
           \begin{array}{cc}
 -\{h_2;p\}&0\\0&-\{h_2;p\}
\end{array}
            \right)
             k+r_6,
$$
for some $r_6\in S^s$, and thus
$$
i[KH-HK]=-\Psi_{\{h_2;p\}}K+\Psi_{r_7},
$$
with $r_7\in S^s$. Hence,
$$
\begin{array}{rcl}
I&=&\displaystyle2\,\mbox{Re}\, \langle K(iH+B_{diag})\Lambda\vec w, K\Lambda \vec w\rangle_{L^2}
\\
&\leq&\displaystyle 2\,\mbox{Re}\, \langle(B_{diag}-\Psi_{\{h_2;p\}})K\Lambda\vec w, K\Lambda \vec w\rangle_{L^2}+c(N(\vec w))^2.
\end{array}
$$
             
              Next we apply the sharp G\aa rding inequality (Theorem \ref{th 1.3}) to the diagonal 
matrix $B_{diag}-\Psi_{\{h_2;p\}}$. For $l=1,2$,
$$
\begin{array}{l}
\displaystyle\mbox{Re}\, \left(i\sum_{j=1}^n b_{llj}(x)\xi_j-\{h_2;p\}\right)=
              \displaystyle-\sum_{j=1}^n \mbox{Im}\,b_{llj}(x)\xi_j-\{h_2;p\}
\\
\qquad\leq\displaystyle c_0\lambda(|x|)\,\sum_{j=1}^n |\xi_j|-C'c_0\lambda(|x|)|\xi|+c\leq \left(\sqrt{n}-C'\right)
c_0\lambda(|x|)|\xi|+C.
\end{array}
$$

Choosing $C'=1+\sqrt n$ and using that $(1+|\xi|^2)^{1/2}\leq 1+|\xi|$, we obtain
$$
\mbox{Re}\, \left(i\sum_{j=1}^n b_{llj}(x)\xi_j-\{h_2;p\}\right)\leq
-c_0\lambda(|x|)\left(1+|\xi|^2\right)^{1/2}+c_0+C,
$$
and the sharp  G\aa
rding inequality yields
$$
\begin{array}{l}2\,\mbox{Re}\, \left\langle(B_{diag}-\Psi_{\{h_2;p\}})K\Lambda\vec w, K\Lambda \vec w\right\rangle_{L^2} + c(N(\vec w))^2
\\
\qquad\leq - 2\,\mbox{Re} \langle c_0\lambda(|x|)J^1 K \Lambda\vec w, K\Lambda \vec w\rangle_{L^2}+c(N(\vec
w))^2.
\end{array}
$$
Since $(\lambda(|x|))^{1/2}\in C^{\infty}_b$, one has
$$
\lambda\,J^1=\Psi_{\lambda\langle \xi\rangle}=\Psi_{(\lambda)^{1/2}(\langle \xi\rangle)^{1/2}}
\,\Psi_{(\lambda)^{1/2}(\langle \xi\rangle)^{1/2}}
+\Psi_{r_8}, 
$$
with $r_8\in S^0$ and $\langle \xi\rangle= (1+|\xi|^2)^{1/2}$. But
$$
\Psi_{(\lambda)^{1/2}(\langle \xi\rangle)^{1/2}}=(\Psi_{(\lambda)^{1/2}(\langle \xi\rangle)^{1/2}})^*
+\Psi_{r_9},
$$
for some $r_9\in S^{-1/2}$, so
$$
- 2\,\mbox{Re}\, \langle\lambda J^1 K\Lambda\vec w, K\Lambda \vec w\rangle_{L^2}
\leq -2\|(\lambda)^{1/2}J^{1/2}K\Lambda\vec w\|^2_{L^2}+ c(N(\vec w))^2.
$$

 Recalling that $ I=\tilde K K\Lambda+\Psi_{r_{-1}}$, $r_{-1}\in S^{-1}$ and using the symbolic
calculus on diagonal matrices we obtain that
$$
\begin{array}{l}
\displaystyle\left\|\lambda^{1/2}J^{s+1/2}\vec w\right\|_{L^2}\leq\displaystyle \left\|(\lambda^{1/2}J^{1/2})(J^s\tilde K)(K\Lambda)\vec w\right\|_{L^2}
+ cN(\vec w)\\
\qquad\leq\displaystyle \left\|(J^s\tilde K)\left(\lambda^{1/2}J^{1/2}\right)(K\Lambda)\vec w\right\|_{L^2}
+ cN(\vec w)\leq c\left(\left\|(\lambda^{1/2}J^{1/2})(K\Lambda)\vec w\right\|_{L^2}+N(\vec w)\right).
\end{array}
$$

 So the following estimate for the term $I$ in (\ref{eq2.3.11}) is therefore obtained
$$
\begin{array}{rcl}
I&=&2\,\mbox{Re}\, \langle K(iH+B_{diag})\Lambda \vec w,K\Lambda \vec w\rangle_{L^2}\\
&\leq& -c \left\langle \lambda(|x|)J^{s+1/2}\vec w, J^{s+1/2}\vec w\right\rangle_{L^2} +c (N(\vec w))^2.
\end{array}
$$

 The estimate for the term $III$ in (\ref{eq2.3.11}) will depend on which inequality (i), (ii), (iv) in Lemma \ref{le2.3.3} 
is  considered.

 To obtain (i) we write,
$$
III=2\,\mbox{Re}\,\langle K\Lambda \vec f,K\Lambda \vec w\rangle_{L^2}\leq c N(\vec f)\,N(\vec w).
$$
Adding the estimates for the two terms of $\p_t(N(\vec w))^2$ we get
$$
\p_t(N(\vec w))^2\leq c(N(\vec w))^2 + c N(\vec f)\,N(\vec w),
$$
so 
$$
\p_tN(\vec w) \leq c N(\vec w) + c N(\vec f).
$$
Hence,
$$
N(\vec w(t))\leq e^{c\,t}(N(\vec w(0))+c\int_0^t\,N(\vec f(\tau))d \tau),
$$
which proves 
part (i) in 
Lemma \ref{le2.3.3}.

 To obtain (iii), we use again that
$$
III\leq c N(\vec f)\,N(\vec w).
$$

 Therefore, adding estimates one finds that
$$
\p_t(N(\vec w))^2 + c\left \langle \lambda J^{s+1/2}\vec w,J^{s+1/2}\vec w\right\rangle_{L^2} \leq c(N(\vec w))^2 + c
N(\vec f)\,N(\vec w).
$$
 
 Integration from $0$ to $t$ yields (iii).
              
               Finally to obtain (iv) we estimate the term $III$ as follows 
$$ 
\begin{array}{rcl}
III&=&\displaystyle2\,\mbox{Re}\,\langle K\Lambda \vec f,K\Lambda \vec w\rangle_{L^2} 
=2\,\mbox{Re}\,\langle J^{1/2}J^{-1/2} K\Lambda \vec f,K\Lambda \vec w\rangle_{L^2} \\ 
&=&\displaystyle2\,\mbox{Re}\,\langle J^{-1/2} K\Lambda \vec f,J^{1/2} K\Lambda \vec w\rangle_{L^2}
=2\,\mbox{Re}\,\langle \lambda^{-1/2}J^{-1/2} K\Lambda \vec f,\lambda^{1/2} J^{1/2} K\Lambda \vec w\rangle_{L^2}
\\ 
&\leq& \displaystyle2\displaystyle\|\lambda^{-1/2}J^{-1/2} K\Lambda \vec f\displaystyle\|_{L^2}\,\|\lambda^{1/2} J^{1/2} K\Lambda \vec
w\|_{L^2}
\\ 
&=&\displaystyle2\,\left\|\lambda^{-1/2}(J^{-1/2} K\Lambda J^{-s+1/2})(J^{s-1/2} \vec f)\right\|_{L^2} \,
\\
&\cdot&\left\|\lambda^{1/2}\left(J^{1/2} K\Lambda J^{-s-1/2})(J^{s+1/2} \vec w\right)\right\|_{L^2}\\
&\leq& \displaystyle\frac{1}{\epsilon}\,\int_{\BbbR^n}\left|\left(J^{-1/2} K\Lambda J^{-s+1/2}\right)\left(J^{s-1/2} \vec f\right)\right|^2\lambda^{-1}dx
\\
&+&
\displaystyle\epsilon\,\int_{\BbbR^n}\left|\left(J^{1/2} K\Lambda J^{-s-1/2}\right)\left(J^{s+1/2} \vec w\right)\right|^2\lambda dx
\\
&\leq&\displaystyle \frac{c}{\epsilon}\,\int_{\BbbR^n}\left|J^{s-1/2} \vec f\right|^2\lambda^{-1}dx + c\,\epsilon
\,\int_{\BbbR^n}\left|J^{s+1/2} \vec w)\right|^2\lambda dx,
\end{array}
$$
where the last estimate follows from Theorem \ref{th 1.1}  since $J^{-1/2} K\Lambda
J^{-s+1/2}$ and $J^{1/2} K\Lambda J^{-s-1/2}$ have order zero.
$\epsilon$ is a small constant to be chosen. Adding estimates we obtain
$$
\begin{array}{l}
\p_t (N(\vec w))^2 + (c_1-c_2\,\epsilon)\langle \lambda J^{s+1/2}\vec w, J^{s+1/2}\vec w\rangle_{L^2}
\\
\qquad\displaystyle\leq c(N(\vec w))^2 +\left(\frac{c_3}{\epsilon}+c_4\right)\left\langle \lambda^{-1} J^{s-1/2}\vec f, J^{s-1/2}\vec
f\right\rangle_{L^2}.
\end{array}
$$
Then, choosing $\epsilon$ sufficiently small, integrating from $0$ to $t$ we obtain part (iv) of
Lemma \ref{le2.3.3}.
               
               \bigskip
               
              {\sl Proof of Theorem \ref{th 2.3.1}\rm} 

 Uniqueness :
 
 Assume $\vec w\in(C([0,T]:H^s(\BbbR^n)))^2$ is a solution of the system in (\ref{eq2.3.2}) with
  $\vec f=0$ and $\vec w_0=0$. Then $\vec w\in(C([0,T]:H^s(\BbbR^n))\cap C^1([0,T]:H^{s-2}(\BbbR^n)))^2$.
   From part (i) in Lemma \ref{le2.3.3} one concludes that $\vec w=0$. 
   
   Existence :
   
   Case 1 : $\vec f\in ({\cal S}(\BbbR^{n+1}))^2$ and $\vec w_0\in  ({\cal S}(\BbbR^n))^2$:
   
    The conjugate linear functional $l^*$ is defined in the linear
    subspace
    $$
    \left(\partial_t+(iH+B+C)^*\right)\left(C^{\infty}_0(\BbbR^n\times [0,T))\right)^2 \subset \left(L^1[0,T]:H^{-s}(\BbbR^n))\right)^2,
    $$
    by
    $$
    l^*(\vec \eta)=\int_0^T\langle \vec f, \vec \phi\rangle_{L^2\times L^2}dt
    +\langle\vec w_0,\vec \phi(0)\rangle_{L^2\times L^2} 
    $$
    for $\vec \phi \in\left(C^{\infty}_0(\BbbR^n\times [0,T))\right)^2$ and $\vec \eta 
    =- (\partial_t+(iH+B+C)^*)\vec \phi$. This is well defined by the uniqueness part above.
    
     By part (ii) of Lemma \ref{le2.3.3} with $s$ replaced by $-s$ it follows that
     $$
     \begin{array}{rcl}
     |l^*(\vec \eta)|&\leq&\displaystyle \int_0^T \|\vec f\|_{(H^s)^2}\|\vec \phi\|_{(H^{-s})^2}dt 
     +\|\vec w_0\|_{(H^s)^2} \|\vec \phi(0)\|_{(H^{-s})^2}
     \\
     &\leq &\displaystyle c_1e^{c_2T}\left(\|\vec f\|_{(L^1([0,T]:H^s))^2}+\|\vec w_0\|_{(H^s)^2}\right)\,
     \|\vec \eta\|_{(L^1([0,T]:H^{-s}))^2}.
     \end{array}
     $$
               
                Using the Hahn-Banach theorem to extend $l^*$, there exists $\vec w\in (L^{\infty}([0,T]:H^s
     (\BbbR^n)))^2$ such that
    \begin{equation}\label{eq2.3.12}
     \begin{array}{rl}
     -&\displaystyle\int_0^T \langle \vec w, (\partial_t+(iH+B+C)^*)\vec \phi\rangle_{ (H^s)^2\times(H^s)^2}\\
     =&\displaystyle \int_0^T\langle \vec f, \vec \phi\rangle_{L^2\times L^2}dt
    +\langle\vec w_0,\vec \phi(0)\rangle_{L^2\times L^2},\;\;\;\;
    \forall \vec \phi\in\left(C^{\infty}_0(\BbbR^n\times [0,T))\right)^2.
    \end{array}
    \end{equation}
    
     Thus, $(\partial_t-(iH+B+C))\vec w=\vec f$ as distributions for $0<t<T$. 
     From this equation one has that $\partial_t\vec w\in (L^{\infty}([0,T]:H^{s-2}
     (\BbbR^n)))^2$ since $\vec f\in ({\cal S}(\BbbR^{n+1}))^2$, so $\vec w\in (C([0,T]:H^{s-2}
     (\BbbR^n)))^2$. Using the equation once more, $\vec w\in (C^1([0,T]:H^{s-4}
     (\BbbR^n)))^2$, and $\vec w(0)=\vec w_0$ by (\ref{eq2.3.12}).
     
      Since  $\vec w_0\in ({\cal S}(\BbbR^n))^2$, $s$ can be replaced  by $s+4$ in the previous
      argument  and there is a solution $\vec w$ of (\ref{eq2.3.2}) to which Lemma \ref{le2.3.3} parts (i)-(iv) hold.

     Case 2 : $\vec w_0 \in (H^s(\BbbR^n))^2$ :
     
      Choose a sequence $(\vec v_j)$ in ${\cal S}(\BbbR^n))^2$ such that $\vec v_j\to \vec w_0$ in 
       $(H^s(\BbbR^n))^2$.
       
       (A) If $\vec f\in (L^1([0,T]:H^s))^2$, choose a sequence $(\vec f_j)$ in  
       $({\cal S}(\BbbR^{n+1}))^2$ such that $\vec f_j\to \vec f$ in $(L^1([0,T]:H^s(\BbbR^n)))^2$.
       
        By case 1 there is a solution $\vec w_j \in (C([0,T]:H^{s+2}
     (\BbbR^n)))^2$ of (\ref{eq2.3.2}) with $\vec f$ and $\vec w_0$ replaced by $\vec f_j$ and $\vec v_j$ 
     respectively.
     Using Lemma \ref{le2.3.3}, part (i), it follows that 
     $(\vec w_j)$ is a Cauchy sequence in $(C([0,T]:H^{s}(\BbbR^n))^2$  and that the limit $\vec w$
      is a solution of (\ref{eq2.3.2}) satisfying part (i) in Lemma \ref{le2.3.3}.

(B) If $\vec f\in (L^2([0,T]:H^s))^2$, choose a sequence $(\vec f_j)$ in  
       $({\cal S}(\BbbR^{n+1}))^2$ such that $\vec f_j\to \vec f$ in $(L^2([0,T]:H^s(\BbbR^n)))^2$.
       Procceding as in (A), there is a solution $\vec w\in (C([0,T]:H^{s}(\BbbR^n)))^2$ of
       (\ref{eq2.3.2}) satisfying (iii) of Lemma \ref{le2.3.3}.
       
       (C) Let
       $$
       J^{s-1/2}\vec f\in \left(L^2\left(\BbbR^n\times [0,T]:\frac{dxdt}{\lambda(|x|)}\right)\right)^2.
       $$
       
  By Theorem \ref{th 1.1} there exists a sequence $(\vec g_j)$
       in $({\cal S}(\BbbR^{n+1}))^2$, such that $\vec g_j\to J^{s-1/2}\vec f$ in \linebreak
       $\displaystyle\left(L^2\left(\BbbR^n\times [0,T]:\frac{dxdt}{\lambda(|x|)}\right)\right)^2$.
       
        Procceding as in (A) with $\vec f_j$ replaced by $J^{s-1/2}\vec g_j\in \left({\cal S}(\BbbR^{n+1})\right)^2$, there is
         a solution $\vec w\in (C([0,T]:H^{s}(\BbbR^n)))^2$ of
       (\ref{eq2.3.2}) satisfying (iv) of Lemma \ref{le2.3.3}.
       
        This completes the proof of Theorem \ref{th 2.3.1}.
                \bigskip

 \begin{remark} \label{re2.3.4} Suppose that the differential operators $B_{11}$ and $B_{22}$ in the
        entries of $B$ in (\ref{eq2.3.2}) are replaced by pseudo-differential operators $\Psi_{b_{11}}$
         and $\Psi_{b_{22}}$ of order $1$ and suppose 
that $C$ in (\ref{eq2.3.2}) is replaced by a $2\times
2$ matrix of
         pseudo-differential operators of order $0$. Then the conclusion of Theorem \ref{th 2.3.1} still holds
         if 
         $$
         |Re\,b_{ll}(x,\xi)|\leq c_0\lambda(|x|)\langle \xi\rangle,\;\;\;\;l=1,2,\;\;\;\;\forall
\,x,
\xi\in\BbbR^n.
         $$
           The reason is that the application of the sharp G\aa rding inequality in the proof
of Lemma \ref{le2.3.3}  
           goes through in exactly the same way.
           
            Consequently, if $\vec b_1(x)\cdot\nabla_x$ in (\ref{eq2.3.2}) is replaced by $\Psi_b$ with $b\in
S^1$ and
            $c_1, c_2$ in (\ref{eq2.3.2}) are replaced by pseudo-differential operators of order $0$, then the
conclusion
            of Corollary \ref{co2.3.2} holds if
            $$
            |Re\,b(x,\xi)|\leq c_0\lambda(|x|)\langle \xi\rangle,\;\;\;\;\;\forall\, x, \xi\in\BbbR^n.
            $$
             \end{remark}

             This will be useful later.
                
\newpage

\section{A NEW CLASS OF SYMBOLS}\label{sec3}

 As it has been shown in \cite{18, 22}
to obtain  local well posedness for nonlinear Schr\"odinger equations
one relies on certain smoothing
effects for
the associated linear equation with lower order terms (order zero and
one).
In the previous section we have established these smoothing effects in
equations with variable 
second order elliptic
coefficients by using known properties of classical pseudo-differential
operators.
 In an attempt to prove these smoothing effects for the non-elliptic case,
one is led (see \cite{4})
to the study of certain
operators with non-standard symbols. Our goal in this
section is
to study results concerning the $L^2$-boundedness and composition of
operators in this class by using geometric arguments. The elliptic
case of our results were proved by Craig, Kappeler and Strauss in
\cite{4} whose statements we follow. The differences between the elliptic  and non-elliptic settings are highlighted in  Proposition \ref{pro3.1.4} below, -see also  \cite{21}.

\subsection{Symbol Properties}\label{subseq3.1}

To begin with, the symbols of interest will be compared to the classical ones defined in
(\ref{eq2.1.2}) Section \ref{sec2}.

  We recall the
following spaces 
\begin{equation}\label{eq3.1.1}
{\cal S}(\BbbR^n)=\left\{u\in C^{\infty}(\BbbR^n)\,:\,\sup_{x\in\BbbR^n}\langle x\rangle^k|\p_x^{\alpha}u(x)| <\infty,\;k\in\BbbN,\;\alpha\in\BbbN^n\right\}, 
\end{equation} 
with seminorms 
\begin{equation}\label{eq3.1.2}
|u|_{{\cal S},m}=\max_{k+|\alpha|\leq m}\left\|\langle
x\rangle^k\,\p_x^{\alpha}u(x)\right\|_{L^{\infty}(\BbbR^n)},
\end{equation}
and
\begin{equation}\label{eq3.1.3}
C^{\infty}_b(\BbbR^n)=\left\{u\in C^{\infty}(\BbbR^n)\,:\,\sup_{x\in\BbbR^
n}
|\p_x^{\alpha}u(x)| <\infty,\;\alpha\in\BbbN^n\right\},
\end{equation}
with seminorms 
\begin{equation}\label{eq3.1.4}
|u|_{C^{\infty}_b,m}=\max_{|\alpha|\leq m}
\left\|\p_x^{\alpha}u(x)\right\|_{L^{\infty}(\BbbR^n)}.
\end{equation}

 The symbol $a=a(x,\xi)\in C^{\infty}(\BbbR^n\times\BbbR^n)$ will
satisfy certain estimates and the
operator $\Psi_a$ associated with the symbol $a$ will
be defined as
$$
\Psi_au(x)=\int_{\BbbR^n}e^{i\,x\cdot\xi}a(x,\xi)\hat
u(\xi)\frac{d\xi}{(2\pi)^{n/2}},\;\;\;\;\;u\in{\cal S}(\BbbR^n).
$$

\begin {proposition}\label{pro3.1.1}

 (i) Suppose $a$ is a classical symbol of order $m\in \BbbR$, $a\in
S^m_{1,0}$, i.e.
 for $\alpha,\,\beta \in \BbbN^n$
\begin{equation}\label{eq3.1.5}
\left|\p_x^{\alpha}\p_{\xi}^{\beta}a(x,\xi)\right|\leq c_{\alpha,\beta}\langle\xi
\rangle^{m-|\beta|},\;\;\forall
\,x,\,\xi\in\BbbR^n.
\end{equation}
Then $\Psi_a$ is a continuous map from ${\cal S}(\BbbR^n)$ into ${\cal
S}(\BbbR^n)$.

 (ii) Suppose $m\in \BbbR$ and $a$ satisfies that for $\alpha \in \BbbN^n$ 
\begin{equation}\label{eq3.1.6}
\left|\p_x^{\alpha}a(x,\xi)\right|\leq c_{\alpha}\langle\xi \rangle^{m},\;\;\forall
x,\,\xi\in\BbbR^n.
\end{equation}
Then $\Psi_a$ is a continuous map from ${\cal S}(\BbbR^n)$ into
$C^{\infty}_b(\BbbR^n)$.
\end{proposition}

{\sl Proof of Proposition \ref{pro3.1.1}\rm}
See \cite{25}.
\bigskip

 The symbols of interest in this section satisfy estimates of the type
\begin{equation}\label{eq3.1.7}
\left|\p_x^{\alpha}\p_{\xi}^{\beta}a(x,\xi)\right|\leq c_{\alpha\beta}\langle
x\rangle^{|\beta|}
\langle\xi\rangle^{m-|\beta|}. 
\end{equation}

This is better than (\ref{eq3.1.6}) in part (ii) of Proposition \ref{pro3.1.1}, but not as 
good
as (i).
In particular, we will see (Proposition \ref{pro3.1.4}, part (iv)) that there exist
$a$ satisfying (\ref{eq3.1.7}) with $m=0$
and $v\in{\cal S}(\BbbR^n)$ such that $\Psi_a v\notin {\cal S}(\BbbR^n)$.

Let $A(x)$ be a real, symmetric and invertible $n\times n$ 
matrix.  Using a coordinate change (a rotation and dilations) in the
$x$-variable, there are essentially only
the elliptic
and ultrahyperbolic cases
$$
A_e=I_n\;\;\;\mbox{and}\;\;\;A_h=
\left(
\begin{array}{cc}
I_k&0\\
0&-I_{n-k}
\end{array}
\right)
,\;\;k\in\{1,..,n-1\}.
$$
where $I_j$ is the $j\times j$ unit matrix. 

Take $\chi\in C^{\infty}(\BbbR^n)$ with $\chi(t)=0$ for
$|t|\leq 1$ and 
$\chi(t)=1$ for $|t|\geq 2$.

 \begin {definition}\label{def3.1.2}

(i) It will be said that $a\in {\cal S}(\BbbR^n:S^m_{1,0})$ if $a\in
C^{\infty}(\BbbR^n\times
\BbbR^n\times \BbbR^n)$ and $a$ satisfies
\begin{equation}\label{eq3.1.8}
\left|\langle s
\rangle^{\mu}\p_s^{\alpha}\p_x^{\beta}\p_{\xi}^{\gamma}\,a(s;x,\xi)\right|\leq
c_{\mu \alpha \beta
\gamma}\langle
\xi\rangle^{m-|\gamma|},\;\;\forall\,s,x,\xi\in\BbbR^n,\;\forall \mu,
\alpha, \beta, \gamma\in\BbbN^n.
\end{equation}

 (ii) For $a\in {\cal S}(\BbbR^n:S^m_{1,0})$, let
\begin{equation}\label{eq3.1.9}
 \left\{
 \begin{array}{l}
(1)\;b_e(x,\xi)=\chi(|\xi|)a(P(x,A_e\xi);x,\xi),\\
(2)\;b_h(x,\xi)=\chi(|\xi|)a(P(x,A_h\xi);x,\xi),
\end{array}
 \right.
\end{equation}
where $P(y,z) = y-(y\cdot z)z/|z|^2$ for $y,z\in\BbbR^n, \;z\ne 0$, is
the projection
of $y$ onto the hyperplane perpendicular to $z$, (notice that $P(y,z)$ is
homogeneous of degree $0$ in $z$).

\end{definition}
 
 \begin {remark} \label{re3.1.3} (a)  Although we shall work in the class ${\cal
S}(\BbbR^n:S^m_{1,0})$, it will be 
clear that all the results
deduced for this class still hold when just a finite number of 
semi-norms
in ${\cal S}$ and $S^m_{1,0}$ are assumed 
to be finite, i.e. (\ref{eq3.1.8}) with $|\mu| +|\alpha|+|\beta|+|\gamma|\leq N$
for $N$ large enough.

(b)  We observe that if  $a\in {\cal S}(\BbbR^n:S^m_{1,0})$, then
$$
\p^{\alpha}_{\xi}(\xi^{\beta}a(\cdot))\in{\cal S}(\BbbR^n:S^k_{1,0}),\;\;k=m+|\beta|-|\alpha|,
$$
and for $M\in\BbbN$ large enough
 $$
\langle x\rangle^{-M} \chi(|\xi|)a(P(x,A_l\xi);x,\xi)=\langle
x\rangle^{-M} b_l(x,\xi),\;\mbox{with}\;\;\;l=e,\;\mbox{or }\;\;l=h, 
$$
is ``roughly speaking" a symbol in the class $S^m_{1,0}$ (when only finitely
many $\xi$ derivatives are taken into account, which is always the case in the sequel).

(c) Finally notice that if $a\in {\cal S}(\BbbR^n:S^m_{1,0})$ and $b_l$ is
defined as in (\ref{eq3.1.9})
then $b^\beta$ given by
$$b^\beta_l=\langle x\rangle^{-|\beta|}\p^\beta_{\xi}
b_l(x,\xi),\;\mbox{with}\;\;\;l=e,\;\mbox{or }\;\;l=h,$$
is a symbol of the same type and the corresponding bounds in (\ref{eq3.1.8}) are controlled by those of $a$.

(d) The symbols described here will be our basic building blocks in the 
study 
of variable coefficients Schr\"odinger operators, see Section \ref{sec5}.

\end{remark}
 
 The symbols defined in (\ref{eq3.1.9}) satisfy an estimate of the type given
in Proposition \ref{pro3.1.1}.
More precisely, if $u\in{\cal S}(\BbbR^n)$ then
$\Psi_{b_e}u\in{\cal S}(\BbbR^n)$ and 
$\Psi_{b_h}u\in C^{\infty}_b(\BbbR^n)$ and is rapidly decreasing away from
the characteristic
directions, i.e. $Ax\cdot x=0$. In the characteristics directions, $\Psi_
{b_h}u(x)$ decays as
$|x|^{1-n}$ as
$|x|\to\infty$. 

\begin {proposition} \label{pro3.1.4}

 Let $u\in{\cal S}(\BbbR^n)$,  and $a$, $b_e$ and $b_h$ as in definition 
\ref{def3.1.2}

(i) If $|\alpha+\beta|\leq k$, then for all $x$ in $\BbbR^n$
$$
\left|x^{\alpha}\p_x^{\beta}\Psi_{b_e}u(x)\right|\leq 
c_k\,\left(\sup_{|\alpha|\leq k;\,|\beta|\leq k;\,|\gamma|\leq k}\langle s
\rangle^{k}\left|\p_s^{\alpha}\p_x^{\beta}\p_{\xi}^{\gamma}\,a(s;x,\xi)\right|\right)\,|\hat
u|_{{\cal
S},2k+m+n+1}.
$$

(ii) If $c\in[0,1]$ and $|\alpha+\beta|\leq k$, then  if $\displaystyle\frac{|A_hx\cdot
x|}{|A_hx||x|}\geq c$
$$
c^k\left|x^{\alpha}\p_x^{\beta}\Psi_{b_h}u(x)\right|\leq 
c_k\,\left(\sup_{|\alpha|\leq k;\,|\beta|\leq k;\,|\gamma|\leq k}\langle s
\rangle^{k}\left|\p_s^{\alpha}\p_x^{\beta}\p_{\xi}^{\gamma}\,a(s;x,\xi)\right|\right)\,
|\hat u|_{{\cal S},2k+m+n+1}.
$$

(iii) If $a\in C^{\infty}_0(\BbbR^n:S^m_{1,0})$ with $a(s;x,\xi)=0$ if
$|s|>1$ and $|\beta|\leq k$, then
$$
\left|\p_x^{\beta}\Psi_{b_h}u(x)\right|\leq 
c_k\,\left(\sup_{|\alpha|\leq k;\,|\beta|\leq k;\,|\gamma|\leq k}
\left|\p_s^{\alpha}\p_x^{\beta}\p_{\xi}^{\gamma}\,a(s;x,\xi)\right|\right)\,|\hat
u|_{{\cal S},2k+m+n+1}\langle
x\rangle^{1-n},\;\;\forall x\in\BbbR^n.
$$

(iv) There exist $a\in {\cal S}(\BbbR^n:S^m_{1,0})$, in fact, $a\in
C^{\infty}_0(\BbbR^n:S^m_{1,0})$,
$v\in{\cal S}(\BbbR^n)$, and
$c>0$ such that
$$
|\Psi_{b_h}v(x)|\geq c|x|^{1-n},\;\;\mbox{with}\; A_hx\cdot x=0,\;|x|\geq
10
$$
\end{proposition}

In the proof of Proposition \ref{pro3.1.4} we will use the following results.
 
 \begin {proposition} \label{pro3.1.5}

(i) Let $\alpha\in\BbbN^n$. Then
$$
(\xi\cdot\nabla_{\xi})\p_{\xi}^{\alpha}=\p_{\xi}^{\alpha}(\xi\cdot\nabla_
{\xi})-|\alpha|\p_{\xi}^{\alpha}.
$$

(ii) Let $T^t$ denote the transpose of the operator $T$, i.e. $\int Tu\,
v=\int u\,T^tv$. Then
$$
\left[\frac{1}{i\,x\cdot
\xi}(\xi\cdot\nabla_{\xi})\right]^t=\frac{1}{i\,x\cdot \xi}
(-(\xi\cdot\nabla_{\xi})-(n-1)I),
$$
and
$$
(\xi\cdot\nabla_{\xi})\left[\frac{1}{i\,x\cdot
\xi}\right]=-\frac{1}{i\,x\cdot \xi}.
$$

 (iii) Let $\phi$ be a differentiable and homogeneous function of degree
$0$ on $\BbbR^n-\{0\}$. 
Then
$$
(\xi\cdot\nabla_{\xi})\phi(\xi)=0,\;\;\;\xi \ne 0.
$$

(iv) Let $\alpha, \beta\in\BbbN^n$, $k\in\BbbN$, $b=b_e$ or
$b=b_h$. Then
$$
\left|\left(\xi\cdot\nabla_{\xi}\right)^k\p_x^{\alpha}\p_{\xi}^{\beta}b(x,\xi\right)|\leq
c_k\,\left(\sup_{|\gamma|\leq k}\,
\left|\p_s^{\alpha}\p_x^{\beta}\p_{\xi}^{\gamma}\,a(s;x,\xi)\right|\right)\,
\langle x\rangle^{|\beta|} \langle \xi\rangle^{m-|\beta|}.
$$
\end{proposition}
 
 {\sl Proof of Proposition \ref{pro3.1.5}\rm}
 
  (i) and (ii) are easily verified. As for (iii) we have
  $$
  \begin{array}{l}
  \displaystyle\left(\xi\cdot\nabla_{\xi}\right)\left[\phi\left(\frac{\xi}{|\xi|}\right)\right]=
  \sum_{j=1}^n\xi_j \p_{\xi_j}\left[\phi\left(\frac{\xi}{|\xi|}\right)\right]\\ 
 \qquad\displaystyle=
\sum_{j=1}^n\xi_j\sum_{k=1}^n(\p_{\xi_k}\phi)\left(\frac{\xi}{|\xi|}\right)\p_{\xi_
j}
  \left(\frac{\xi_k}{|\xi|}\right)=
\sum_{j,k}(\p_{\xi_k}\phi)\left(\frac{\xi}{|\xi|}\right)\xi_j\left(\frac{\delta_
{jk}}{|\xi|}-
  \frac{\xi_j\xi_k}{|\xi|^3}\right)\\ 
\qquad\displaystyle =
  \sum_{k}(\p_{\xi_k}\phi)\left(\frac{\xi}{|\xi|}\right)\frac{\xi_k}{|\xi|}-
\sum_{k}(\p_{\xi_k}\phi)\left(\frac{\xi}{|\xi|}\right)\xi_k\sum_{j}\frac{\xi^2_j}{|\xi|^3}=0.
  \end{array}
  $$
  
  (iv) follows from (i), (iii) and Definition \ref{def3.1.2}, since $P(x,A\xi)$ is
homogeneous of degree $0$ in $\xi$.
  \bigskip
 
 {\sl Proof of Proposition \ref{pro3.1.4}\rm} 

For simplicity of the exposition we shall 
assume that the constants
$c_{\mu\alpha\beta\gamma}$ which appear in (\ref{eq3.1.8}) are all smaller than 
unity for $|\mu|\leq k$,
$|\alpha|\leq k$, $|\beta|\leq k$ and $|\gamma|\leq k$. 
Also we shall drop the powers of $2\pi$ which appear in the definition of 
$\Psi_b$.

(i) Let $|\alpha+\beta|\leq
k$. We consider three cases:

Case 1 : $|x|\leq 2$. Here
$$
\left|x^{\alpha}\p_x^{\beta}\Psi_{b_e}u(x)\right|\leq c_k\,|\hat
u|_{{\cal S},m+k+n+1}.
$$

Case 2 : $|x|\geq 2\,$ and $\,|P(x,\xi)|\leq |x|/2$. From (\ref{eq3.1.7}), it
suffices to estimate terms of the type
\begin{equation}\label{eq3.1.10}
I=\int e^{i\,x\cdot
\xi}\p_{\xi}^{\alpha_1}\left[{\xi^{\beta_1}}\right]\p_x^{\beta_2}\p_{\xi}^{\alpha_2}
\left[b_e(x,\xi)\right]\p_{\xi}^{\alpha_3}\hat u(\xi)d\xi,
\end{equation}
where $\alpha_1+\alpha_2+\alpha_3=\alpha$, and $\beta_1+\beta_2=\beta$.

 Pythagoras' theorem gives
$$
|x|^2=\frac{|x\cdot \xi|^2}{|\xi|^2} + |P(x,\xi)|^2,
$$
so by hypothesis
$$
\begin{array}{l}
\displaystyle|x\cdot \xi|=|x||\xi|\sqrt{1-\frac{|P(x,\xi)|^2}{|x|^2}}\geq
|x||\xi|\left(1-\frac{|P(x,\xi)|^2}{|x|^2}\right)\\
\qquad\displaystyle\geq |x||\xi|\left(1-\frac{1}{4}\right)\geq \frac{3}{4}|x||\xi|.
\end{array}
$$
Above we have used that $\sqrt{1-t}\geq 1-t$ for $t\in[0,1]$. Now we use
the identity
$$
e^{i\,x\cdot \xi}=\left(\frac{1}{i\,x\cdot
\xi}(\xi\cdot\nabla_{\xi})\right)^k e^{i\,x\cdot \xi}.
$$

 By Proposition \ref{pro3.1.5} (ii), if $k_1+k_2+k_3\leq k$ it suffices to estimate
terms of the type
$$
\begin{array}{l}
\displaystyle II=
\int \frac{e^{i x\cdot \xi}}{(i
x\cdot \xi)^k}\left(\xi\cdot\nabla_{\xi}\right)^{k_1}\left[\p_{\xi}^{\alpha_1}
[\xi^{\beta_1}]\right]\left(\xi\cdot\nabla_{\xi}\right)^{k_2}\p_x^{\beta_2}
\p_{\xi}^{\alpha_2}b_e(x,\xi)\left(\xi\cdot\nabla_{\xi}\right)^{k_3}\p_{\xi}^{\alpha_3}
\hat u(\xi)d\xi.
\end{array}
$$
Using Proposition  \ref{pro3.1.5} (iv) it follows that
$$
\begin{array}{l}
\displaystyle |II|
\leq c_k\int_{|\xi|\geq 1}\frac{1}{|x|^k|\xi|^k}
\langle\xi
\rangle^{|\beta_1-\alpha_1|+m-|\alpha_2|}\,\langle x\,\rangle^{|\alpha_2|}|\hat u|_{{\cal S},2k+m+n+1}
\le c_k|\hat u|_{{\cal S},2k+m+n+1}.
\end{array}
$$

Case 3 : $|x|\geq 2\,$ and $\,|x|<2|P(x,\xi)|$.
 From (\ref{eq3.1.8}) and (\ref{eq3.1.10}) it follows that
 $$
 |I|\leq c_k\,|\hat u|_{{\cal S},m+k+n+1},
 $$
 which proves (i).

(ii) Let $c\in[0,1]$, $|\alpha+\beta|\leq k$. Suppose that
$|x|^{-2}|A_hx\cdot x|\geq c$ with
$c>0$,
(otherwise the statement is trivial).  We consider three cases:

Case 1 : $|x|\leq 4/c$. Hence
$$
c^k|x^{\alpha}\p_x^{\beta}\Psi_{b_h}u(x)|\leq 4^kc_k
\,|\hat u|_{{\cal S},k+m+n+1}.
$$

Case  2 : $|x|\geq 4/c\,$ and $\,|P(x,A_h\xi)|\leq c|x|/2$.

Since 
$$
x=\frac{ x\cdot A_h\xi}{|\xi|^2}\,A_h\xi +P(x,A_h\xi),
$$
by Pythagoras' theorem,
$$
|x|^2 |\xi|^2 = |x\cdot A_h\xi|^2 +|\xi|^2 |P(x,A_h\xi)|^2.
$$
The matrix $A_h$ is symmetric, so
$$
\begin{array}{l}
\displaystyle\left|\xi\cdot A_h x\right|=\left|x\cdot A_h\xi\right|
=|x||\xi|\sqrt{1-\frac{|P(x,A_h\xi)|^2}{|x|^2}}\\
\qquad\displaystyle\geq |x||\xi|\sqrt{\left(1-\frac{c^2}{4}\right)}\geq \frac{\sqrt{3}}{2}|x||\xi|.
\end{array}
$$

 Next we write
$$
\xi =\frac{\xi\cdot A_h x }{|x|^2}\,A_h x + P(\xi,A_h x).
$$

 By Pythagoras' theorem 
$$ 
\left|P(\xi,A_h x)\right|=|\xi|\sqrt{1-\frac{|\xi\cdot A_h x|^2}{|\xi|^2|x|^2}}\leq
|\xi|\sqrt{1-\left(1-\frac{c^2}{4}\right)}=
\frac{c}{2}|\xi|.$$ 

 Using that 
 $$ 
x\cdot \xi=|x|^{-2}\left(\xi\cdot A_h x\right)(A_hx\cdot x)+x\cdot P(\xi,A_h x), 
$$ 
from the above estimates and our hypothesis we
conclude that 
\begin{equation}\label{eq3.1.11}
|x\cdot \xi|\geq c\frac{\sqrt{3}}{2}|x|\cdot|\xi|-|x|\cdot|P(x,A_h\xi)|
\geq c\frac{\sqrt{3}}{2}|x|\cdot|\xi| -\frac{c}{2}|x\cdot|\xi|\geq
\frac{c}{4}|x|\cdot|\xi|
 \end{equation}

Integrating by parts as in case 2 of part (i), (ii) follows.

Case  3 : $|x|\geq 4/c\,$ and $\,|P(x,A_h\xi)|\geq c|x|/2$.

From (\ref{eq3.1.9}) and (\ref{eq3.1.10}) as in case 3 of the elliptic case it follows that
$$
\left|\partial_x^{\beta}\Psi_{b_h}u(x)\right|\leq c^{-k}c_0
\,|\hat u|_{{\cal S},m+k+n+1}.
$$
 
 (iii) It suffices to show the statement for $|x|\geq 10$. Let
$|\beta_1+\beta_2|\leq k$ and consider terms
of the type
$$
I=\int e^{ix\cdot \xi } (i\xi)^{\beta_1}\p_x^{\beta_2}b_h(x,\xi)\hat
u(\xi)d\xi.
$$

 Since  $|P(x,A_h \xi)|\leq 1$ proceeding as before
$$
x=\frac{x\cdot A_h\xi}{|\xi|^2}A_h\xi + P(x,A_h\xi),
$$
and by Pythagoras' theorem
$$
|x|^2|\xi|^2=|x\cdot A_h\xi|^2+|\xi|^2|P(x,A_h\xi)|^2.
$$
The matrix $A_h$ is symmetric, so
$$
|\xi\cdot A_hx|=|x\cdot
A_h\xi|=|x||\xi|\sqrt{1-\frac{|P(x,A_h\xi)|^2}{|x|^2}}\geq
|x||\xi|(1-|x|^{-2}).
$$
 Next we write 
$$
 \xi=\frac{\xi\cdot A_hx}{|x|^2}A_hx+P(\xi,A_hx).
$$
By Pythagoras' theorem
$$
|P(\xi,A_hx)|=|\xi|\sqrt{1-\frac{|\xi\cdot A_hx|^2}{|\xi|^2|x|^2}}\leq
|\xi|\sqrt{1-(1-|x|^{-2})}\leq \sqrt{2}|\xi||x|^{-1},
$$
i.e.
\begin{equation}\label{eq3.1.12}
\frac{|P(\xi,A_hx)|}{|\xi|}\leq \sqrt{2}|x|^{-1}.
\end{equation}
 
 Hence $\xi$ is in a cone $\Gamma_x$ with vertex at the origin, axis given
by $A_hx$ and opening angle
$\theta$ where $sin(\theta)=\sqrt{2}|x|^{-1}$. In particular, $\theta\leq
2|x|^{-1}$, because $|x|\geq
10$.

 Therefore
$$
I=\int_{\Gamma_x}\,e^{i x\xi} (i\xi)^{\beta_1}\p_x^{\beta_2}b_h(x,\xi)\hat u(\xi)d\xi
$$
and
$$
\begin{array}{rcl}
|I|&\leq&\displaystyle c\,|\hat u|_{{\cal
S},k+m+n+1}\,\int_{\Gamma_x}\,\langle \xi\rangle^{-n-1}d\xi\\
&\leq&\displaystyle c\,|\hat u|_{{\cal S},k+m+n+1}
(|x|^{-1})^{n-1}\,\int_{\BbbR^n}\,\langle
\xi\rangle^{-n-1}d\xi\\
&=&\displaystyle c\,|\hat u|_{{\cal S},k+m+n+1} |x|^{-n+1}.
\end{array}
$$

(iv) Choose $\phi\in C^{\infty}_0(\BbbR^n)$ with $\phi\geq 0$,
$\phi(s)=1$ if $|s|\leq 1/4$
 and $\phi(s)=0$ if $|s|\geq 1/2$. Let $a(s;x,\xi)=\phi(s)\langle\xi
\rangle^m$. Choose $v\in
{\cal S}(\BbbR^n)$ such that $\hat v(\xi)=\chi(|\xi|)(1-\chi(|\xi|))^2$
with $\chi(\cdot)$ given
 in Definition \ref{def3.1.2}. Then 
\begin{equation}\label{eq3.1.13}
\Psi_{b_h}v(x)=\int e^{i x\cdot \xi}\phi\left(P(x,A_h\xi)\right)\langle \xi\rangle^m
\left(\chi(|\xi|)\right)^2\left(1-\chi(|\xi|)\right)^2 d\xi.
\end{equation}
 
  Assume that $|x|\geq 10$. An argument similar to that used to deduce
(\ref{eq3.1.12}) shows that 
$$
\left|P(\xi,A_h x)\right|\leq \left(\sqrt {2}\right)^{-1}|\xi||x|^{-1}.
$$
Also one has that $x\cdot \xi = x\cdot P(\xi,A_h x)$ when $ A_hx\cdot
x=0$. Hence 
$|x\cdot \xi |\leq (\sqrt
{2})^{-1}|\xi|\leq \sqrt {2}$ on the $\xi$- support of the integrand in
(\ref{eq3.1.13}). Thus
$$
Re\, \Psi_{b_h}v(x)\geq \cos\sqrt {2}\int \phi\left(P(x,A_h\xi)\right)\langle
\xi\rangle^m
\chi(|\xi|)^2\left(1-\chi(|\xi|)\right)^2 d\xi.
$$
  Let $\Gamma_x$ be the cone defined by
  $$ \Gamma_x=\{\xi\in \BbbR^n\,:\,|\xi|^{-1}|P(\xi,A_h x)|\leq (4\sqrt{2}|x|)^{-1}\}. $$
   Suppose
$\xi\in \Gamma_x$. Since $$ \xi=\frac{\xi\cdot A_hx}{|x|^2}\,A_hx +
P(\xi,A_h x), $$ it follows from Pythagoras' theorem that 
   $$ | x\cdot
A_h\xi|=|\xi\cdot A_hx|=|\xi||x|\sqrt{1-\frac{|P(\xi,A_h x)|^2}
{|x|^2}}\geq|\xi||x|\left(1-\frac{1}{32}|x|^{-2}\right). $$
    Now write 
    $$ x=\frac{
x\cdot A_h\xi}{|\xi|^2}A_h\xi+ P\left(x,A_h \xi\right), $$
     and use Pythagoras'
theorem to get 
     $$ \left|P(x,A_h \xi)\right|^2=|x|^2-\frac{|x\cdot
A_h\xi|^2}{|\xi|^2}\leq
|x|^2-|x|^2\left(1-\frac{1}{32}|x|^{-2}\right)^2\leq\frac{1}{16}. $$
      Therefore 
      $$
\begin{array}{rcl}
       Re\,\Psi_{b_h}v(x&\geq&\displaystyle \cos(\sqrt 2)\int_{\Gamma_x}\langle
\xi\rangle^m (\chi(|\xi|))^2(1-\chi(|\xi|))^2d\xi\\ 
       &=&\displaystyle \cos(\sqrt
2)c_n|x|^{1-n}\,\int_{\BbbR^n} \langle
\xi\rangle^m(\chi(|\xi|))^2(1-\chi(|\xi|))^2d\xi =c\,c_n \cos (\sqrt
2)|x|^{1-n},
       \end{array} 
       $$
       where 
       $$ c=\int_{\BbbR^n}\langle
\xi\rangle^m(\chi(|\xi|))^2(1-\chi(|\xi|))^2d\xi>0. $$ 
This completes the
proof of Proposition \ref{pro3.1.4}. 
  \bigskip
  
  \subsection {$L^2$-boundedness}\label{subseq3.2}

 Now we are ready to establish the $L^2$-boundedness of our ultrahyperbolic operators. 

\begin {theorem} \label{th3.2.1}

Suppose $a\in {\cal S}(\BbbR^n:S^m_{1,0})$ (see Definition \ref{def3.1.2},
(\ref{eq3.1.8})-(\ref{eq3.1.9})). Then there exists $c=c(m)$ and $N=N(n)$ 
such that 
$$
\left\|\Psi_{b_h}u\right\|_{L^2}\leq c \max_{\mu+|\alpha+\beta+\gamma|\leq
N}\left\|\langle s\rangle^{\mu}\langle \xi\rangle^{-m+|\gamma|}
\partial_s^{\alpha}\p_x^{\beta}\p_{\xi}^{\gamma}a\right\|_{L^{\infty}}
\|u\|_{H^m},\;\;\;\;u\in{\cal S}(\BbbR^n). 
$$ 
\end{theorem}
  
  The proof of this theorem in the 2-dimensional case ($n=2$) is more involved than in the higher dimensional case. As it was pointed out in Proposition \ref{pro3.1.4} (iv) one can not expect enough decay in cones around the characteristic directions (i.e. $A_h x\cdot x=0$).  In fact, for $n=2$ the estimate $|x|^{-1}$ is critical and after some decomposition in frequency and space we need to use Cotlar-Stein lemma to glue the pieces together.

{\sl Proof of Theorem \ref{th3.2.1}\rm}

 If $a\in{\cal S}(\BbbR^n:S^m_{1,0})$, then $\langle\xi\rangle^{-m}a\in
{\cal S}(\BbbR^n:S^0_{1,0})$ with 
 $$
\begin{array}{l}
\displaystyle \max_{\mu+|\alpha+\beta+\gamma|\leq N}\left\|\langle
s\rangle^{\mu}\langle \xi\rangle^{|\gamma|}\p_s^{\alpha} \p_x^{\beta}\p_
{\xi}^{\gamma}
\left[\langle\xi
\rangle^{-m}a\right] \right\|_{L^{\infty}}\\
\qquad\qquad\displaystyle \leq c_m \max_{\mu+|\alpha+\beta+\gamma|\leq N}\left\|\langle s\rangle^{\mu}
\langle
\xi\rangle^{-m}\p_s^{\alpha}\p_x^{\beta}
\p_{\xi}^{\gamma}a\right\|_{L^{\infty}}.
\end{array}
$$
  
   Furthermore, $J^m$ is an isometry of $H^m$ onto $L^2$, so it can be
assumed that $m=0$.
   
   Next, we use an argument similar to that in \cite{21} Section 3 (which in the case of the elliptic operators gives a straightforward proof of the $L^2$-continuity). This is based on a simple change of variables -see (\ref{numero})  below. In order to do it we make the following decomposition.   
 Let $\{\phi_j\}_{j\in\BbbN}$ be a smooth partition of unity on $\BbbR$
subordinate to the covering 
   $$ 
   \BbbR=(-1,1)\bigcup \cup_{j=1}^{\infty}\left\{t\in\BbbR :2^{j-2}<|t|<2^j\right\}
   $$
i.e. satisfying $\displaystyle 0\leq \phi_j(x)\leq
1,\;1=\sum_{j=0}^{\infty}\phi_j$, $\mbox{supp}\,\phi_0\subset (-1,1)$ and
$\mbox{supp}\, \phi_j\subset \left\{t\in\BbbR :2^{j-2}<|t|<2^j\right\}$. 

 Since the intervals in the covering have length at least $1$, it can be
assumed that 
   $$
    \left|\frac{d^k}{dt^k}\phi_j(t)\right|\leq
c_k,\;\;\;\;t\in\BbbR,\;k\in\BbbN.  
    $$

 Let $\phi\in C^{\infty}(\BbbR)$ with $\phi(t)=1$ for $t\leq 1$,
$\phi(t)=0$
for $t\geq 2$ and $0\leq \phi\leq 1$.

 Denote by $\hat a$ the Fourier transform of $a$ in the
$s$-variable.  Thus, 
 $$ 
\left \|\langle
\sigma\rangle^{N_1}\p_{\sigma}^{\beta_1}\hat a(\sigma;x,\xi)\right\|_{L^
{\infty}}
 \leq c\max_{N_2+|\beta_2|\leq
2N_1+n+1}\left\|\langle s\rangle^{N_2}\p_s^ {\beta_2}a\right\|_{L^{\infty}}.
$$
    
     Now using that $\sigma\cdot P(x,A_h\xi)=x\cdot P(\sigma,A_h\xi)$ we have
 $$
 \begin{array}{rcl} 
      \Psi_{b_h}u(x)&=&\displaystyle \int e^{i x\cdot
\xi}a\left(P(x,A_h\xi);x,\xi\right)\chi(|\xi|)\hat u(\xi)d\xi\\ 
&=&\displaystyle \int\left(\int e^{i
x\cdot( \xi+P(\sigma,A_h\xi))} \hat a(\sigma;x,\xi)\chi(|\xi|)\hat
u(\xi)d\xi \right)d\sigma\\
&=&\displaystyle \sum_{j=0}^{\infty} \int \left(\int e^{i x\cdot(
\xi+P(\sigma,A_h\xi)}\right)\phi_j(|\sigma|)\hat a(\sigma;x,\xi)\chi(|\xi|)
\left(1-\phi\left(\frac{\delta |\xi|}{2^{j+1}}\right)\right)\hat u(\xi)d\xi)d\sigma\\
&+&\displaystyle \sum_{j=0}^{\infty}\int\left(\int e^{i x\cdot(
\xi+P(\sigma,A_h\xi)}\right)\phi_j(|\sigma|)\hat a(\sigma;x,\xi)\chi(|\xi|)
\phi\left(\frac{\delta |\xi|}{2^{j+1}}\right)\hat u(\xi)d\xi)d\sigma\\ 
      &=&I+II,
\end{array}
       $$
        where $\delta>0$ is a small constant to be chosen. In $I$
we make the change of variables
        \begin{equation}\label{numero}
\eta=\xi+P(\sigma,A_h\xi)=\xi+\left(\sigma-\frac{\sigma\cdot
A_h\xi}{|\xi|^2}A_h\xi\right). 
        \end{equation}
         Then 
         $$
          \left(\frac{\partial \eta}{\partial
\xi}\right)= I - \sum_{k=1}^n\frac{\sigma_k}{|\xi|} M_k(\xi), 
          $$
           where
$M_k$ is a matrix whose entries are homogeneous of degree $0$ in $\xi$.
The determinant function is continuous on $\BbbR^{n^2}$, so 
           $$
\frac{1}{2}\leq \left|\det\left(\frac{\partial \eta}{\partial \xi}\right)\right|\leq
2, 
           $$
            by fixing $\delta>0$ sufficiently small, we recall that
$|\sigma|\sim 2^j,\,\delta|\xi|\geq2^{j+1}$. 
This gives 
$$ \begin{array}{rcl}
\displaystyle I&=&\displaystyle
            \int\,\sum_{j=0}^{\infty}\phi_j(|\sigma|)\,\int e^{i x\cdot \eta}\hat
a(\sigma;x,\xi(\eta)) \chi(|\xi(\eta)|)\\
&\cdot&\displaystyle
            \left(1-\phi(\delta2^{-j-1}|\xi(\eta)|)\right)\hat u(\xi(\eta))
\left|\det\left(\frac{\partial \eta}{\partial \xi}\right)\right|^{-1}\,d\eta\,d\sigma.
\end{array} 
$$             
 Combining Minkowski's integral inequality and the $L^2$-boundedness of
$S^0_{1,0}$
 pseudo-differential operators we get
 $$
 \|I\|_{L^2_x}\leq c \max_{\mu+|\alpha+\beta+\gamma|\leq N}\left\|\langle
s\rangle^{\mu}\langle \xi\rangle^{|\gamma|}
\partial_s^{\alpha}\p_x^{\beta}\p_{\xi}^{\gamma}a\right\|_{L^{\infty}}\,\|u\|_
{L^2}\sum_{j=0}^{\infty}
\int \frac{\phi_j(|\sigma|)}{\langle\sigma \rangle^{n+1}}d\sigma, 
            $$ 
where
$$ \sum_{j=0}^{\infty} \int \frac{\phi_j(|\sigma|)}{\langle\sigma
\rangle^{n+1}}d\sigma <\infty. $$
 As for $II$, let
 $$
 a_{jk}(s;x,\xi)= \phi_k(|s|)\int e^{i \sigma\cdot s}\phi_j(|\sigma|) \hat
a(\sigma;x,\xi) d\sigma,
 \;\;\;\;j, k\in\BbbN.
 $$ Then 
            \begin{equation}\label{eq3.2.1}
            \begin{array}{rcl}
             II&=&\displaystyle\sum_{j,k=0}^{\infty}\int e^{i x\cdot
\xi}a_{jk}\left(P(x,A_h\xi);x,\xi\right)\chi(|\xi|)\phi\left(\delta 2^{-j-1}|\xi|\right) \hat
u(\xi)d \xi\\ 
             &=&\displaystyle\sum_{j,k=0}^{\infty} \Psi_{jk}u(x),
              \end{array} 
              \end{equation}
               where $\Psi_{jk}$ is a pseudo-differential operator with symbol 
               $$
a_{jk}\left(P(x,A_h\xi);x,\xi\right)\chi(|\xi|)\phi\left(\delta 2^{-j-1}|\xi|\right). $$

 First, by Cauchy-Schwarz,
 $$
\left |\Psi_{jk}u(x)\right| \leq c 2^{jn/2}\left\|a_{jk}\right\|_{L^{\infty}}\|u\|_{L^2},
 $$
 so
 \begin{equation}\label{eq3.2.2}
 \left\|\chi_{B_{2^{k+10}}(0)}\Psi_{jk}u\right\|_{L^2} \leq c
2^{(j+k)n/2}\left\|a_{jk}\right\|_{L^{\infty}}\|u\|_{L^2}.
 \end{equation}
 Next, we shall estimate
 $$
\left \|(1-\chi_{B_{2^{k+10}}(0)})\Psi_{jk}u\right\|_{L^2}
 $$
 by using Cotlar-Stein lemma which can be stated as follows.
               \medskip
 
 {\bf Cotlar-Stein lemma.}\rm
 
  \em Let $\{\Psi_l\}_{l=0}^{\infty}$ be a sequence of bounded operators on
$L^2$ and let
 $\{\gamma(l)\}_{l=0}^{\infty}$ be a sequence of positive numbers with
$\displaystyle\sum_{l=0}^{\infty}\gamma(l)<\infty$.
 Suppose 
 $$
\left \|\Psi_{l_1}^* \Psi_{l_2}\right\|,\,\left\|\Psi_{l_1}\Psi_{l_2}^*\right\|\leq
\left(\gamma(l_2-l_1)\right)^2,\;\;l_1, l_2\in\BbbN,\;l_1\leq l_2.
 $$
 Then $\,\displaystyle\sum_{l=0}^{\infty} \Psi_l u$ converges in $L^2$ for any $u\in
L^2$ and
 $$
 \left\|\sum_{l=0}^{\infty} \Psi_l u \right\|_{L^2}\leq \sum_{l=0}^{\infty}
\gamma(l)\,\|u\|_{L^2}.
 $$\rm
 
 For a proof of this lemma we refer to \cite{30}.

 For fixed $j, k\geq 0$, let
 $$
 A_l=\left\{x\in\BbbR^n\,:\,2^{k+l+10}<|x|\leq 2^{k+l+11}\right\},\;\;\;l\in\BbbN,
 $$
 and define $\Psi_l=\Psi_{ljk}$ by
 $$
 \Psi_l u(x)=\chi_{A_l}(x)\Psi_{jk}u(x),\;\;\;u\in{\cal S}(\BbbR^n),\;\;l\in\BbbN.
 $$
 Then $\displaystyle\sum_{l=0}^{\infty}\Psi_l u $ converges to
$\left(1-\chi_{B_{2^{k+10}}(0)}\right)\Psi_{jk}u$ uniformly on compact
 subsets of $\BbbR^n$. In particular in the distribution sense. By
Cauchy-Schwarz
 $$
\left\|\Psi_l u\right\|_{L^2}\leq c
2^{(j+k+l)n/2}\left\|a_{jk}\right\|_{L^{\infty}}\|u\|_{L^2}, \;\;\;u\in{\cal S}(\BbbR^n),
 $$
 so $\Psi_l$ is bounded on $L^2$ for $l\in \BbbN$.
 
  In order to apply Cotlar-Stein lemma, it suffices to consider $
\Psi_{l_1}\Psi_{l_2}^*$ for $l_1\leq l_2$, because
  $(\Psi_{l_2}\Psi_{l_1}^*)^*=\Psi_{l_1}\Psi_{l_2}^*$ and $\Psi_{l_1}^*
\Psi_{l_2}=0$ for
  $l_1\neq l_2$.
  
   Let $u\in{\cal S}(\BbbR^n)$. Since
   $$
   \begin{array}{rcll}
    \Psi^*_lu(x)&=&\displaystyle\int\int e^{i(x-y)\cdot \xi}\bar a_{jk}\left(P(y,A_h\xi);y,\xi\right)\chi(|\xi|
)
   \phi\left(\delta 2^{-j-1}|\xi|\right)\chi_{A_l}(y)u(y)dyd\xi\\
   &=&\displaystyle\int e^{i x\cdot \xi}\left(\int e^{-i y \cdot\xi}\bar
a_{jk}\left(P(y,A_h\xi);y,\xi\right)\chi(|\xi|)
   \phi\left(\delta 2^{-j-1}|\xi|\right)\chi_{A_l}(y)u(y)dy\right)d\xi
   \end{array}
   $$
   it follows that
   $$
   \begin{array}{rcl}
   \langle \Psi^*_{l_1} v,\Psi^*_{l_2} u\rangle&
   =&\displaystyle\int\int\int e^{-i y \cdot\xi}\bar
a_{jk}\left(P(y,A_h\xi);y,\xi\right)\chi(|\xi|)
   \phi\left(\delta 2^{-j-1}|\xi|\right)\chi_{A_{l_1}}(y)v(y)\cdot\\
   &\cdot&\displaystyle e^{i z\cdot\xi} a_{jk}\left(P(z,A_h\xi);z,\xi\right)\chi(|\xi|)
   \phi\left(\delta 2^{-j-1}|\xi|\right)\chi_{A_{l_2}}(z)\bar u(z)d\xi dz dy,
   \end{array}
   $$
   for all $v, u\in{ \cal S}(\BbbR^n)$. So
   $$
   \Psi_{l_1}\Psi_{l_2}^*u(y)=\int\,K_{l_1 l_2}(y,z) u(z)dz,
   $$
   where
   $$
   \begin{array}{rcl}
   K_{l_1 l_2}(y,z) &=&\chi_{A_{l_1}}(y) \chi_{A_{l_2}}(z)\int
e^{i(y-z)\cdot\xi}a_{jk}(P(y,A_h\xi;y,\xi)\cdot\\
   &&\cdot \bar a_{jk}(P(z,A_h\xi);z,\xi)\chi^2(|\xi|)\phi^2(\delta
2^{-j-1}|\xi|) d\xi.
   \end{array}
   $$
   
    It is convenient to fix the following terminology.

    \begin{definition} \label{def3.2.2} 
    
     $\Gamma $ is the cone of angle $\theta\in(0,1/2)$ if there exists
$e_0\in
    \BbbR^n-\{0\}$ such that
    $$
    \Gamma=\left\{x\in\BbbR^n\,:\,|x|^{-1}|P(x,e_0)|\leq \theta\right\}.
    $$
    
    Notice that
    $|\Gamma\cap A_l|\leq c \theta^{n-1} 2^{n(k+l)}$.
    \end{definition}
   
   Choose a collection $\{\Gamma_{m_1}\}$ of cones of angle at most
$2^{-l_1}/16$ such that
    $$
    \BbbR^n=\displaystyle\bigcup_{m_1} \Gamma_{m_1}.
    $$
    Let
    $$
    \Gamma_{m_1}^*=\left\{z\in\BbbR^n\,:\,|z|^{-1}|P(z,x)|\leq
\frac{2^{-l_1}}{16}\;\,\mbox{for some}\;\,x\in\Gamma_{m_1}\right\}.
    $$
     Then $\Gamma_{m_1}^*$ is a cone of angle at most ${2^{-l_1}}/4$, and
it can be
     assumed that the collection $\{\Gamma_{m_1}^*\}$ is locally finite 
on
$\BbbR^n-\{0\}$ by choosing 
     $\{\Gamma_{m_1}\}$ appropriately. Let 
     $$
     S_{m_1}=\Gamma_{m_1}\cap A_{l_1},\;\;\;S_{m_1}^*=\Gamma_{m_1}^*\cap
A_{l_2}.
     $$
   
     \begin{claim} :\label{claim3.2.3}
If $x\in S_{m_1}$, $ z\in A_{l_2}$ and $a_{jk}\left(P(x,A_h\xi);x,\xi\right)
     a_{jk}\left(P(z,A_h\xi);z,\xi\right)\ne 0$, then $z\in S_{m_1}^*$.
     \end{claim}
     
     \sl{Proof of claim \ref{claim3.2.3}}\rm  : The goal is to show that $z\in
S_{m_1}^*$. By assumption,
     $$
    \left |P(x,A_h\xi)\right|,\,\left|P(z,A_h\xi)\right|\leq 2^k.
     $$
      Since
      $$
      x=\frac{x\cdot A_h\xi}{|\xi|^2}A_h\xi + P(x,A_h\xi),
      $$
      by Pythagoras' theorem 
      $$
      |x\cdot A_h\xi|\geq
|x||\xi|\sqrt{1-\frac{\left|P(x,A_h\xi)\right|^2}{|x|^2}}\geq |x||\xi|
\left(1-2^{-2l_1-20}\right).
      $$
      Similarly,
      $$
      \left|z\cdot A_h\xi\right|\geq |z|  |\xi| \left(1-2^{-2l_2-20}\right).
      $$
       Now we use the identity
       $$
       x\cdot z=|\xi|^{-2}\left(x\cdot A_h\xi\right)\left(z\cdot A_h\xi\right)+P(x,A_h\xi)\cdot
P(z,A_h\xi),
       $$
       to obtain
       $$
       |x\cdot z|\geq |x| |z|\left(1-2^{-2l_1-20}\right)\left(1-2^{-2l_2-20}\right)-2^{2k}.
       $$
       But
       $$
       z=\frac{z\cdot x}{|x|^2} x+P(z,x),
       $$
       so
 $$
       \begin{array}{rcl}
       |P(z,x)|^2&=&\displaystyle|z|^2-(z\cdot x)^2|x|^{-2}\leq
|z|^2-|z|^2\left(1-2^{-2l_1-19}\right) + |x|^{-2}2^{2k+1}\\
       &\leq& \displaystyle2^{-2l_1-19}|z|^2+2^{-2l_1-19}\leq 2^{-2l_1-16}|z|^2,
       \end{array}
       $$
       and therefore $z\in \Gamma_{m_1}^*$. This proves claim \ref{claim3.2.3}.
    
 \begin{claim} :
     \label{claim3.2.4}
$$
\sup_{x,z\in\BbbR^n}|K_{l_1 l_2}(x,z)|\leq c\left\|a_{jk}\right\|^2_{L^{\infty}}
2^{nj} 2^{-l_2(n-1)}.$$
\end{claim}

  \sl{Proof of claim \ref{claim3.2.4}}\rm 
       
 We recall that
$$
\begin{array}{rcl}
K_{l_1l_2}(x,z)&=&\displaystyle\chi_{A_{l_1}}(x)\chi_{A_{l_2}}(z)\int e^{i
(x-z)\cdot\xi}a_{jk}\left(P(x,A_h\xi);x,\xi\right)\cdot \\
&\cdot&\displaystyle \bar
a_{jk}\left(P(z,A_h\xi);z,\xi\right)\chi^2(|\xi|)\phi^2\left(\delta2^{-j-1}|\xi|\right)d\xi.
\end{array}
$$
     
      Suppose $z\in A_{l_2}$ (otherwise $K_{l_1l_2}(x,z)=0$) and $\bar
a_{jk}(P(z,A_h\xi);z,\xi)\ne 0$ so that
$|P(z,A_h\xi)|\leq 2^k$. Then
$$
\left|\xi\cdot A_hz\right|=\left|z\cdot A_h\xi\right|\geq |z||\xi|\left(1-2^{-2l_2-20}\right),
$$
as in the proof of claim \ref{claim3.2.3}. Since
$$
\xi=\frac{\xi\cdot A_hz}{|z|^2}A_hz+P\left(\xi,A_hz\right),
$$
one has
$$
\left|P(\xi,A_hz)\right|^2=|\xi|^2-\frac{|\xi\cdot A_hz|^2}{|z|^2}
\leq |\xi|^2-|\xi|^2\left(1-2^{-2l_2-20}\right)^2\leq |\xi|^2 2^{-2l_2-19},
$$
which shows that $\xi$ belongs to a cone $\Gamma_z$ of angle
$2^{-l_2-9}$. Hence
$$
\left|K_{l_1l_2}(x,z)\right|\leq
\left\|a_{jk}\right\|^2_{L^{\infty}}\left|\Gamma_z\cap\left\{\delta|\xi|\leq 2^{j+2}\right\}\right|\leq c
\left\|a_{jk}\right\|^2_{L^{\infty}}2^{nj} 2^{-l_2(n-1)}.
$$
This yields claim \ref{claim3.2.4} .
 
Now to estimate $\displaystyle\left\|\Psi_{l_1}\Psi^*_{l_2} u\right\|_{L^2}$ we consider two
separate cases.

\underline{Case 1} : $l_1\leq l_2\leq l_1+10$.
      
       Then
$$
\begin{array}{rcl}
\|\Psi_{l_1}\Psi^*_{l_2} u\|_{L^2}^2&=&\displaystyle\int\left|\int
K_{l_1l_2}(x,z)u(z)dz\right|^2dx\leq
\sum_{m_1}\int_{S_{m_1}}\left|\int_{S_{m_1}^*}K_{l_1l_2}(x,z)u(z)dz\right|^2dx\\
&\leq&\displaystyle \sum_{m_1}\left|S_{m_1}\right|\left|S_{m_1}^*\right|
\sup_{x,z}\left|K_{l_1l_2}(x,z)\right|^2\left\|\chi_{S_{m_1}^*} u\right\|_{L^2}^2\\
&\leq&\displaystyle
c2^{-l_1(n-1)}2^{n(k+l_1)}2^{-l_1(n-1)}2^{n(k+l_2)}\left\|a_{jk}\right\|^4_{L^{\infty}}2^{2nj-2l_2(n-1)}
\sum_{m_1}\left\|\chi_{S_{m_1}^*} u\right\|_{L^2}^2\\
&\leq&\displaystyle c 2^{2n(k+j)}
2^{(2-n)(l_1+l_2)}\left\|a_{jk}\right\|^4_{L^{\infty}}\|u\|_{L^2}^2.
\end{array}
$$
 
 The first inequality above follows from claim \ref{claim3.2.3}, the second inequality
from Cauchy-Schwarz, the third inequality from
claim \ref{claim3.2.4}, and the fourth inequality from the local finiteness of
$\{S^*_{m_1}\}$. Therefore,
\begin{equation}\label{eq3.2.3}
\left\|\Psi_{l_1}\Psi^*_{l_2}\right\|\leq c 2^{n(k+j)}
2^{(2-n)(l_1+l_2)/2}\left\|a_{jk}\right\|^2_{L^{\infty}},\;\;\;
\mbox{for}\;\;\;l_1\leq l_2\leq l_1+10.
\end{equation}
       
        Before turning to the remaining case $l_2\geq l_1+11$, it is useful to
split $A_{l_2}$ in sectors
$S_{m_2}$ where $m_2$ roughly speaking measures how noncharacteristic the
directions in
$S_{m_2}$ are. More precisely, let
$$
\begin{array}{l}
\displaystyle S_{m_2}=\left\{z\in A_{l_2}\,:\,(m_2-1)2^{-l_2}\leq \frac{|A_hz\cdot
z|}{|z|^2}\leq m_2
2^{-l_2}\right\},\;\;m_2=1,2,..,2^{l_2-1},\\ 
\displaystyle       S_{l_2}^*=\left\{z\in
A_{l_2}\,:\,\frac{1}{2}\leq \frac{|A_hz\cdot
z|}{|z|^2}\leq 1\right\}.
\end{array}
$$
Then $\displaystyle A_{l_2}= S_{l_2}^* \bigcup \bigcup_{m_2=1}^{2^{l_2-1}} S_{m_2}$. The next
result will be used to estimate the volume
of $S_{m_2}$. Recall that
$$
A_h=
\left(
        \begin{array}{cc}
         I_k&0\\
         0&-I_{n-k}
         \end{array}
         \right),\;\;k\in\{1,..,n-1\},
$$
where $I_j$ is the $j\times j$ unit matrix.

\begin {proposition} \label{pro3.2.5}

Suppose that $0<\epsilon \leq a\leq 1/2$, $\epsilon\leq 1/16$. Let
$$
S_{a,\epsilon}= \left\{z\in\BbbR^n\,:\,|z|\leq 1\;\;\mbox{and}\;\;a-\epsilon
\leq \frac{|A_hz\cdot z|}{|z|^2}\leq
a\right\}.
$$
Then
$$
|S_{a,\epsilon}|\leq c_{n,k}\,\epsilon,
$$
where $c_{n,k}$ is independent of $a$.
\end{proposition}
         
        {\sl Proof of Proposition \ref{pro3.2.5}\rm}

 Let $m_j$ denote the Lebesgue measure in $\BbbR^j$. Write 
$$ 
S_{a,\epsilon}=S_{a,\epsilon}^+\cup
S_{a,\epsilon}^-, 
$$ 
where 
$$ 
S_{a,\epsilon}^{\pm}=\left\{z\in\BbbR^n\,:\,|z|\leq
1\;\;\mbox{and}\;\;a-\epsilon\leq \pm \frac{|A_hz\cdot z|}{|z|^2}\leq a\right\}. 
$$
Let $z^+=(z_1,..,z_k)$, $z^-=(z_{k+1},..,z_n)$. By Fubini's theorem
$$
m_n\left(S_{a,\epsilon}^+\right)=\int_{|z^-|\leq
1}m_k\left(S^{+}_{a,\epsilon}(z^-)\right)dm_{n-k}\left(z^-\right),
$$
where
$$
\begin{array}{rcl}
S_{a,\epsilon}^+(z^-)&=&\displaystyle\left\{z^+\in\BbbR^k\,:\,(z^+,z^-)\in
S^+_{a,\epsilon}\right\}\\
&=&\displaystyle \left\{z^+\in\BbbR^k\,:\,|z^+|^2+|z^-|^2\leq 1,\;\;a-\epsilon\leq
\frac{|z^+|^2-|z^-|^2}{|z^+|^2+|z^-|^2}\leq a\right\}.
\end{array}
$$
 By straightforward calculation, 
$$ 
\frac{|z^+|^2-|z^-|^2}{|z^+|^2+|z^-|^2}\in [a-\epsilon,a]\Longleftrightarrow
|z^+|^2\in\left[\frac{1+(a-\epsilon)}{
1-(a-\epsilon)}|z^-|^2,\frac{1+a}{1-a}|z^-|^2\right].  
$$ 
Let 
$$f(t)=\frac{1+t}{1-t},\;\;\;t\in(-\infty,1).  
$$ 
Then 
$$ 
f'(t)=\frac{2}{(1-t)^2}, 
$$ 
so 
$\sup_{[0,1/2]}|f'(t)|=8$. By the fundamental theorem of calculus, 
$$ 
\left|\frac{1+(a-\epsilon)}{1-(a-\epsilon)}-\frac{1+a}{ 1-a}\right|\leq
8\epsilon,
$$
so
$$
\frac{1+(a-\epsilon)}{1-(a-\epsilon)}\geq \frac{1+a}{ 1-a}-8 \epsilon.
$$
Also,
$$
\begin{array}{rcl}
\displaystyle\sqrt{\frac{1+a}{1-a}-8 \epsilon}&=&
\displaystyle\sqrt{\frac{1+a}{1-a}}\cdot \sqrt{1-8\epsilon \frac{1-a}{1+a}}\\
&\geq&\displaystyle
\sqrt{\frac{1+a}{1-a}}\cdot \left(1-8\epsilon\frac{1-a}{1+a}\right)
\geq \sqrt{\frac{1+a}{ 1-a}}\cdot (1-8\epsilon).
\end{array}
$$
 It follows that
$$
S_{a,\epsilon}^+(z^-)\subset\left\{z^+\in\BbbR^k\,:\,
\sqrt{\frac{1+a}{ 1-a}}\cdot (1-8\epsilon)\leq \frac{|z^+|}{|z^-|}\leq
\sqrt{\frac{1+a}{ 1-a}}\,\right\},
$$
and hence
$$
m_k\left(S_{a,\epsilon}^+(z^-)\right)\leq c_k
\left(\frac{1+a}{1-a}\right)^{k/2}|z^-|^k\epsilon
\leq c_k 3^{k/2}|z^-|^k\epsilon.
$$
Therefore,
$$
\begin{array}{rcl}
\displaystyle m_n(S_{a,\epsilon}^+)&=&\displaystyle\int_{|z^-|\leq 1}
m_k\left(S_{a,\epsilon}^+(z^-)\right)dm_{n-k}(z^-)\\
&\leq&
\displaystyle c_k\epsilon \int_{|z^-|\leq 1}|z^-|^kdm_{n-k}(z^-)\leq c_{n,k}\epsilon.
\end{array}
$$
 Replacing $k$ by $n-k$ in the above argument for $S_{a,\epsilon}^+$, one
gets
$$
m_n(S_{a,\epsilon}^-)\leq c_{n,k}\epsilon.
$$
This proves Proposition \ref{pro3.2.5}.
         \bigskip

\begin{claim} :\label{claim3.2.6} 
          Let $m_2\in\{1,2,..,2^{l_2-1}\}$. Then
$$
|S_{m_2}|\leq c 2^{(n-1)l_2+nk}.
$$
\end{claim}
\sl{Proof of claim \ref{claim3.2.6}}  \rm 
This follows by Proposition \ref{pro3.2.5} and homogeneity.
         
         \begin{claim} :\label{claim3.2.7}
          Suppose $l_2\geq l_1+11$, $x\in A_{l_1}$, $z\in
A_{l_2}$, $10\leq m_2\leq 2^{l_2}$,
$$
\frac{|A_hz\cdot z|}{|z|^2} \geq (m_2-1)2^{-l_2},
$$
and
$$
a_{jk}\left(P(x,A_h\xi);x,\xi\right) a_{jk}\left(P(z,A_h\xi);z,\xi\right)\ne 0.
$$
Then $\,|(z-x)\cdot\xi|\geq cm_2 2^k|\xi|$.
\end{claim}
          
\sl{Proof of claim \ref{claim3.2.7}} \rm
 The identity
$$
z=\frac{z\cdot A_h\xi}{|\xi|} A_h\xi + P\left(z,A_h\xi\right)
$$
and Pythagoras' theorem give
$$
\left|\xi\cdot A_hz\right|=\left|z\cdot A_h\xi\right|=|z||\xi|\sqrt{1-\frac{P(z,A_h\xi)}{|z|^
2}}
\geq |z||\xi|\left(1-2^{-2l_2-20}\right).
$$
Also
\begin{equation}\label{eq3.2.4}
\xi=\frac{\xi\cdot A_hz}{|z|^2}A_hz+P(\xi,A_hz),
\end{equation}
so
$$
\left|P(\xi,A_hz)\right|^2=|\xi|^2-\frac{|\xi\cdot A_hz|^2}{|z|^2}
\leq |\xi|^2-|\xi|^2\left(1-2^{-2l_2-20}\right)^2\leq |\xi|^22^{-2l_2-19}.
$$
 Using again (\ref{eq3.2.4}) 
$$
z\cdot \xi=|z|^{-2}\left(\xi\cdot A_hz\right)\left(A_hz\cdot z\right)+z\cdot P\left(\xi,A_hz\right),
$$
and therefore
$$
\begin{array}{rcl}
\displaystyle|z\cdot\xi|&\geq&\displaystyle
|z||\xi|\left(1-2^{-2l_2-20}\right)(m_2-1)2^{-l_2}-|z||\xi|2^{-l_2-9}\\
&\geq&\displaystyle 2^{k+10}|\xi|\left(1-2^{-2l_2-20}\right)(m_2-1)-2^{k+2}|\xi|\geq
2^{k+9}|\xi|m_2.
\end{array}
$$
Another application of (\ref{eq3.2.4}) gives
$$
A_h\xi\cdot \xi=|z|^{-2}\left(\xi\cdot A_hz\right)\left(A_h\xi\cdot A_hz\right) +A_h\xi\cdot
P\left(\xi,A_hz\right).
$$
Here $A_h\xi\cdot A_hz=\xi\cdot z$,
so
$$
\begin{array}{rcl}
\left|A_h\xi\cdot \xi\right|&\geq&\displaystyle
|z|^{-2}|z||\xi|\left(1-2^{-2l_2-20}\right)2^{k+9}|\xi|m_2-|\xi|^2 2^{-l_2-9}\\
&\geq&\displaystyle |\xi|^2\left(1-2^{-2l_2-20}\right)2^{-l_2-2}m_2-|\xi|^2 2^{-l_2-9}\geq
2^{-l_2-3}m_2|\xi|^2.
\end{array}
$$
 Finally,
$$
(z-x)\cdot\xi=\frac{(z-x)\cdot A_h\xi}{|\xi|^2}\left(A_h\xi\cdot \xi\right)+\xi\cdot
P\left(z,A_h\xi\right)-\xi\cdot P\left(x,A_h\xi\right).
$$
 Here
$$
\left|(z-x)\cdot A_h\xi\right|\geq \left|z\cdot
A_h\xi\right|-|x||\xi|\geq|z||\xi|\left(1-2^{-2l_2-20}\right)-|x||\xi|\geq\frac{1}{2}|z||\xi|,
$$
so
$$
|(z-x)\cdot\xi|\geq \frac{1}{2}|z||\xi|m_22^{-l_2-3}-2^{k+1}|\xi|\geq
m_22^{k+6}|\xi|-2^{k+1}|\xi|\geq
m_22^{k+5}|\xi|.
$$
This proves claim \ref{claim3.2.7}.
\bigskip

\begin{claim} : \label{claim3.2.8}
          Suppose $l_2\geq l_1+11$, $10\leq m_2\leq 2^{l_2}$
and
$$
\frac{|A_hz\cdot z|}{|z|^2}\geq (m_2-1)2^{-l_2}.
$$
Then
$$
 |K_{l_1l_2}(x,z)|\leq cm_2^{-2} 2^{(n-2) j-2k-(n-1)l_2} \max_{N\leq
2}\left\|(\xi\cdot\nabla_{\xi})^Na_{jk}\right\|^2_{L^{\infty}},\;\;\forall x\in\BbbR^n.
$$
\end{claim}
\sl{Proof of claim \ref{claim3.2.8}} \rm  For $x\in A_{l_1}$ and $z\in A_{l_2}$
(otherwise $K_{l_1l_2}(x,z)=0$), it follows 
from claim \ref{claim3.2.7} and parts (ii) and (iii) of Proposition \ref{pro3.1.5} that
$$
\begin{array}{rcl}
\left|K_{l_1l_2}(x,z)\right|&=&\displaystyle\left|\int
e^{i(x-z)\cdot\xi}\left(\left(\frac{1}{i(x-z)\cdot\xi}\xi\cdot 
\nabla_{\xi}\right)^2\right)^t\right|\cdot
\\
&\cdot&\displaystyle\left [a_{jk}\left(P(x,A_h\xi);x,\xi\right) \bar
a_{jk}\left(P(z,A_h\xi);z,\xi\right)\chi^2(|\xi|)\phi^2\left(\delta2^{-j-1}|\xi|\right)\right]d\xi\\
&\leq&\displaystyle cm_2^{-2} 2^{(n-2)j-2k-(n-1)l_2}\max_{N\leq 2}\left\|(\xi\cdot
\nabla_{\xi})^Na_{jk}\right\|^2_{L^{\infty}},
\end{array}
$$
since the integrand is $0$ for $\xi$ not in the cone with angle
$2^{-l_2-9}$ (see the proof of claim \ref{claim3.2.3}). 
This proves claim \ref{claim3.2.8}.

 With these results $\left\|\Psi_{l_1}\Psi^*_{l_2}u\right\|_{L^2}$ can be estimated
in the remaining case.

\underline{Case 2} : $l_2\geq l_1+11$.
         
         $$
\left\|\Psi_{l_1}\Psi^*_{l_2}u\right\|^2_{L^2}\leq
\sum_{m_1}\int_{S_{m_1}}\left|\int_{A_{l_2}}K_{l_1l_2}(x,z)u(z)dz\right|^2dx.
$$

Writing
$$
\int_{A_{l_2}}=\int_{S^*_{l_2}}+\sum_{m_2=1}^9\int_{S_{m_2}}+\sum_{m_2=
10}^{2^{l_2-1}}\int_{S_{m_2}},
$$
and then 
$$
\left\|\Psi_{l_1}\Psi^*_{l_2}u\right\|_{L^2}\leq 11(I_1+I_2+I_3),
$$
where
$$
\begin{array}{l}
\displaystyle I_1=\sum_{m_1}\int_{S_{m_1}}\left(\int_{S^*_{l_2}}\left|K_{l_1l_2}(x,z)u(z)\right|dz\right)^
2dx,\\
\displaystyle I_2=\sum_{m_1}\int_{S_{m_1}}\sum_{m_2=1}^9\left(\int_{S_{m_2}}\left|K_{l_1l_
2}(x,z)u(z)\right|dz\right)^2dx,\\
\displaystyle I_3=\sum_{m_1}\int_{S_{m_1}}\left|\sum_{m_2=10}^{2^{l_2-1}}\int_{S_{m_
2}}\left|K_{l_1l_2}(x,z)u(z)\right|dz\right|^2dx.
\end{array}
$$

 Using claims \ref{claim3.2.3} and \ref{claim3.2.8} with $m_2=2^{l_2-1}$ (integration by parts),
$$
\begin{array}{rcl}
|I_1|&\leq &\displaystyle c 2^{nk+l_1} 2^{nk+nl_2-(n-1)l_1}
2^{2(n-2)j-4k-2l_2(n-1)-4l_2}\cdot\\
&\cdot &\displaystyle\max_{N\leq 2}\left\|\left(\xi\cdot \nabla_{\xi}\right)^Na_{jk}\right\|^4_{L^{\infty}}
\sum_{m_1}\left\|\chi_{S_{m_1}^*}u\right\|^2_{L^2}\\
&\leq&\displaystyle
 c 2^{2(n-2)j} 2^{k(2n-4)} 2^{(2-n)l_1} 2^{(-2-n)l_2} 
\max_{N\leq 2}\left\|(\xi\cdot \nabla_{\xi})^Na_{jk}\right\|^4_{L^{\infty}}
\|u\|_{L^2}^2.
\end{array}
$$

 Using claims \ref{claim3.2.4} and \ref{claim3.2.6} (size)
$$
\begin{array}{rcl}
|I_2|&\leq&\displaystyle c 2^{nk+l_1} 2^{nk+(n-1)l_2}
2^{2nj-2l_2(n-1)}\left\|a_{jk}\right\|^4_{L^{\infty}}\sum_{m_1}\left\|\chi_{S_{m_1}^*}u\right\|
^2_{L^2}\\
&\leq&\displaystyle c 2^{2n(k+j)}
2^{l_1-l_2(n-1)}\left\|a_{jk}\right\|^4_{L^{\infty}} \|u\|^2_{L^2}.
\end{array}
$$

 Using claim \ref{claim3.2.6} (size)  and \ref{claim3.2.8} with $m_2\in\{10,11,...2^{l_2-1}\}$
(integration by parts), and the
Cauchy-Schwarz inequality on the sum in $m_2$,
$$
\begin{array}{rcl}
|I_3|&\leq&\displaystyle c
2^{nk+l_1}\left(\sum_{m_2=10}^{2^{l_2-1}}\frac{1}{m_2^2}\right) \left(\sum_{m_2=10}^{2^
{l_2-1}}m_2^2
2^{nk+l_2(n-1)}m_2^{-4}\right)\cdot\\
&\cdot&\displaystyle 2^{2(n-2)j-4k-2l_2(n-1)} \max_{N\leq 2}\left\|(\xi\cdot
\nabla_{\xi})^Na_{jk}\right\|^4_{L^{\infty}}\sum_{m_1}\left\|\chi_{S_{m_1}^*}u\right\|^2_
{L^2}\\
&\leq&\displaystyle c 2^{2(n-2)j} 2^{k(2n-4)} 2^{l_1-l_2(n-1)} \max_{N\leq 2}\left\|(\xi\cdot
\nabla_{\xi})^Na_{jk}\right\|^4_{L^{\infty}}\|u\|^2_{L^2}.
\end{array}
$$

 Combining the above estimates for $I_1, I_2$ and $I_3$ , one obtains for $l_1+11\leq
l_2$
\begin{equation}\label{eq3.2.5}
\left\|\Psi_{l_1}\Psi^*_{l_2}\right\|_{L^2}\leq c 2^{n(k+j)}\max_{N\leq 2}\left\|(\xi\cdot
\nabla_{\xi})^Na_{jk}\right\|^2_{L^{\infty}}2^{(l_1-l_2)/2}2^{-l_2\frac{n-2}{2}}.
\end{equation}
         
          By (\ref{eq3.2.3}), (\ref{eq3.2.5}) is valid for all $l_1\leq l_2$. Notice that from
(\ref{eq3.2.3}) and (\ref{eq3.2.5}) we conclude that the
estimates are much better for $n>2$ than for $n=2$. In fact it is just in
dimension $n=2$ where Cotlar-Stein lemma
is necessary.

 Taking $\gamma(l)=2^{-l/4}$ it follows from Cotlar-Stein lemma that
$\displaystyle\sum_{l=0}^{\infty}\Psi_lu$ converges in $L^2$ for any $u\in{\cal S}(\BbbR^n)$. 
As it was previously noted,
$\displaystyle\sum_{l=0}^{\infty}\Psi_lu$ converges to $(1-\chi_{B_{2^{k+10}}})\Psi_{jk}u$ in the
distribution sense, so 
$$ 
\left\|\left(1-\chi_{B_{2^{k+10}}}\right)\Psi_{jk}u\right\|_{L^2}\leq c 2^{n(k+j)/2}\max_{N\leq
2}\left\|\left(\xi\cdot
\nabla_{\xi}\right)^Na_{jk}\right\|^2_{L^{\infty}} \|u\|_{L^2},\;\;\;\forall u \in{\cal
S}(\BbbR^n).
$$
Combining with (\ref{eq3.2.2})  
$$
\left\|\Psi_{jk}u\right\|\leq c2^{n(k+j)/2}
\max_{N\leq 2}\left\|\left(\xi\cdot
\nabla_{\xi}\right)^Na_{jk}\right\|^2_{L^{\infty}}
\|u\|_{L^2}.
$$

 Finally, one can use the following fact which is easily verified using
standard Fourier
transform arguments.
          
          \begin{lemma} \label{le3.2.9}

 For any $N_1, N_2\in\BbbN$ there exists $N_3\in\BbbN$ such that
$$
2^{N_1(j+k)}\left\|(\xi\cdot
\nabla_{\xi})^{N_2}a_{jk}\right\|_{L^{\infty}}
\leq c_{N_1}\max_{k+|\alpha|\leq N_3}\left\|\langle s \rangle^k(\xi\cdot
\nabla_{\xi})^{N_2}\partial_s^{\alpha}
a\right\|_{L^{\infty}}.
$$
\end{lemma}

 Using Lemma \ref{le3.2.9} with $N_1\geq n/2+1$ one
obtains 
$$
\|II\|_{L^2_x}\leq \max_{k+|\alpha|\leq N_1\,;\, N\le2}\left\|\langle s \rangle^k(\xi\cdot
\nabla_{\xi})^N_1\partial_s^{\alpha}
a\right\|_{L^{\infty}} \|u\|_{L^2_x},
$$
where $II$ was defined in (\ref{eq3.2.1}), and the proof of Theorem \ref{th3.2.1} is
completed.
\bigskip

\subsection { Composition results} \label{subseq3.3}

 First a few facts concerning oscillatory integrals will be listed. For
details and further results, see
\cite{25}, section 1.6.

           \begin {definition} \label{def3.3.1}
 Let $m\in\BbbR$, $\tau\geq 0$. Then ${\cal A}_{\tau}^m$ is the class of
functions
$\phi \in C^{\infty}(\BbbR^n_y\times \BbbR^n_{\xi})$ satisfying
$$
\left|\partial_y^{\alpha}\partial_{\xi}^{\beta}\phi(y,\xi)\right|\leq
c_{\alpha,\beta}\langle \xi\rangle^m
\langle y \rangle^{\tau}.
$$
The class ${\cal A}$ of amplitude functions is defined by
$$
{\cal A}=\bigcup_{m\in \BbbR} \bigcup_{\tau\geq 0}{\cal A}_{\tau}^m.
$$
\end{definition}

\begin {definition} \label{def3.3.2}

Let $a\in{\cal A}$. Then
$$
{\cal O}_s \int \int e^{-iy\cdot \xi}a(y,\xi) dy d\xi=\lim_{\epsilon \to
0}\int \int e^{-iy\cdot \xi}a(y,\xi) 
\tilde\chi(\epsilon y,\epsilon \xi)dy
d\xi,
$$
if $\tilde\chi\in{\cal S}\left(\BbbR^n_y\times \BbbR^n_{\xi}\right)$ and $\tilde\chi(0,0)=1$.
\end{definition}

 The oscillatory integral in Definition \ref{def3.3.2} is well-defined because of the
following lemma
which allows one to integrate by parts and use Lebesgue's dominated
convergence theorem.
           
           \begin {lemma} \label{le3.3.3}

 Let $\tilde\chi\in{\cal S}\left(\BbbR^n_y\right)$ with $\tilde\chi(0)=1$. Then 
           $$ 
           \begin{array}{rl}
(i)&\displaystyle\tilde\chi(\epsilon y)\to 1 \;\;\mbox{in}\;\;\BbbR^n\;\;\mbox{uniformly on compact sets},\\
(ii)&\displaystyle\partial_y^{\alpha}\left[\tilde\chi(\epsilon y)\right]\to 0
\;\;\mbox{in}\;\;\BbbR^n\;\;\mbox{uniformly for}\;\; \alpha\in\BbbN^n-\{0\},\\ 
(iii)&\displaystyle\forall \alpha\in\BbbN^n \;\,\exists
\,c_{\alpha}>0\;\;\mbox{s.t.}\;\\
&\displaystyle\left|\partial_y^{\alpha}\left[\tilde\chi(\epsilon y)\right]\right|\leq c_{\alpha}
\epsilon^{\sigma} \langle y\rangle^{\sigma-|\alpha|},\;\;\forall y \in\BbbR^n,\;0\leq \sigma\leq |\alpha|,\;0< \epsilon <1. 
           \end{array}
            $$
\end{lemma}
            
	            \begin {theorem} \label{th3.3.4}
 Let $a(s;x,\xi)\in {\cal S}(\BbbR^n;S^m_{1,0})$, $\alpha\in\BbbN^n$ and
$\phi\in {\cal S}(\BbbR^n)$. Suppose $N\geq m+|\alpha|$, $N\in \BbbN$. Let
\begin{equation}\label{eq3.3.1}
c_1(x,\xi)=\sum_{|\beta|<N}\frac{i^{-|\beta|}}{\beta!}\phi(x)\partial_{\xi}^{\beta}\left[(i\xi)^{\alpha}\right]\partial_x^{\beta}
b_h(x,\xi),
            \end{equation}
             and 
             \begin{equation}\label{eq3.3.2}
c_2(x,\xi)=\sum_{|\beta|<N}\frac{i^{-|\beta|}}{\beta!} 
\partial_{\xi}^{\beta}b_h(x,\xi)\partial_x^{\beta}\phi(x)(i\xi)^{\alpha}.
\end{equation}
             Let
\begin{equation}\label{eq3.3.3}
E_1=\phi\partial_x^{\alpha}\Psi_{b_h}-\Psi_{c_1},
\end{equation}
and
\begin{equation}\label{eq3.3.4}
E_2=\Psi_{b_h}\phi\partial_x^{\alpha}-\Psi_{c_2}.
\end{equation}
Then there exist $N_1\in \BbbN$ and $c>0$ such that for any $u\in{\cal
S}(\BbbR^n)$ and for $j=1,2$ 
\begin{equation}\label{eq3.3.5}
\begin{array}{rcl}
\|E_ju\|_{L^2}&\leq&\displaystyle c\max_{|\alpha|+|\beta|\leq
N_1}\left\|x^{\alpha}\partial_x^{\beta}\phi\right\|_{L^{\infty}}\cdot \\
&\cdot&\displaystyle\max_{|\gamma_1+\gamma_2+\gamma_3+\gamma_4|\leq
N_1}\left\|s^{\gamma_4}\langle \xi\rangle^{-m+|\gamma_3|}
\partial_s^{\gamma_1}\partial_x^{\gamma_2}\partial_{\xi}^{\gamma_3}a\right\|_{L^{\infty}}
\|u\|_{L^2}.
\end{array}
             \end{equation}
\end{theorem}

{\sl Proof of Theorem \ref{th3.3.4}\rm}

 For simplicity of the exposition we shall drop all the powers
of $2\pi$ which appear in the definition of the operators $\Psi_b$.

Let $u\in{\cal S}(\BbbR^n)$. Then
$$
\begin{array}{rcl}
E_1u(x)&=&\displaystyle\phi(x) \partial_x^{\alpha} \int e^{ix\cdot \xi} b_h(x,\xi)\hat
u(\xi)d\xi\\
&-&\displaystyle\phi(x) \int e^{ix\cdot \xi} \sum_{|\beta|<N}\frac{i^{-|\beta|}}{\beta!}
\partial_{\xi}^{\beta}\left[(i\xi)^{\alpha}\right]\partial_x^{\beta}
b_h(x,\xi)\hat u(\xi)d\xi.
\end{array}
$$
  Notice that 
$$
\partial_x^{\alpha}\left[e^{ix\cdot \xi} b_h(x,\xi)\right]=
\sum_{\beta\leq \alpha}
             \left(
             \begin{array}{c}
             \alpha\\ \beta
             \end{array}
             \right)
i^{|\alpha-\beta|}\xi^{\alpha-\beta}e^{ix\cdot \xi}
\partial_x^{\beta} b_h(x,\xi),
$$
and that if $\beta\leq \alpha$
$$
\frac{i^{-|\alpha|}}{\beta!}\partial_{\xi}^{\beta}\left[(i\xi)^{\alpha}\right]
=\frac{i^{|\alpha-\beta|}}{\beta!}\frac{\alpha!}{(\alpha-\beta)!}\xi^{\alpha-\beta}
=i^{|\alpha-\beta|}
              \left(
             \begin{array}{c}
\alpha\\ \beta
 \end{array}
             \right)
\xi^{\alpha-\beta},
$$
then
$$ 
E_1u(x)=\phi(x)\int e^{ix\cdot \xi}\sum_{N\leq |\beta|,\, \beta\leq \alpha}
 \left(
             \begin{array}{c}
\alpha\\ \beta
              \end{array}
             \right)
i^{|\alpha-\beta|}\xi^{\alpha-\beta}\partial_x^{\beta} b_h(x,\xi)\hat
u(\xi)d\xi.
$$

 If $m>0$, then $E_1=0$.

If $m\leq 0$, $N\leq |\beta|$, and $\beta\leq\alpha$, then
$$
\phi(x) \xi^{\alpha-\beta}\partial_x^{\beta} b_h(x,\xi)\in S^0_{1,0}.
$$
 Hence, there exists $c>0$ and $N_1\in\BbbN$ such that (\ref{eq3.3.5}) holds for
$j=1$.

 Let us go back to $E_2$. We shall assume that $\phi\in C^{\infty}_0(\BbbR^n)$ and vanish outside $B_R(0)$. Then we shall obtain bounds of the form of powers of
$R$. Thus by introducing a partition of unity with respect to dyadic
$x$-annuli and summing the corresponding operators 
one obtains the case $\phi\in{\cal S}(\BbbR^n)$.

Write
$$
E_2u(x)=\left(1-\chi\left((2R)^{-1}|x|\right)\right)E_2u(x) + \chi\left((2R)^{-1}|x|\right)E_2u(x)= I + II.
$$

 To estimate $I$ it suffices to consider classical symbols. Indeed,
integrating by parts with respect to 
 $y$ in the second term
\begin{equation}\label{eq3.3.6}
\begin{array}{rcl}
I&=&\displaystyle\left(1-\chi\left((2R)^{-1}|x|\right)\right) \cdot {\cal O}_s\int\int (e^{i(x-y)\cdot\xi}
b_h(x,\xi)\phi(y)\partial_y^{\alpha}u(y)\\
&-&\displaystyle e^{i(x-y)\cdot\xi}\sum_{|\beta|<N}\frac{i^{-|\beta|}}{\beta!}
\partial_{\xi}^{\beta}
b_h(x,\xi)  \partial_x^{\beta}\phi(x)\partial_y^{\alpha}u(y))dyd\xi.
\end{array}
\end{equation}
 By Taylor's formula
$$
\begin{array}{rcl}
\phi(y)&=&\displaystyle\sum_{|\beta|<N}\frac{1}{\beta!}\partial_x^{\beta}\phi(x)(y-x)^{\beta}\\
&+&\displaystyle N\sum_{|\beta|=N}\frac{1}{\beta!}(y-x)^{\beta}\int_0^1\partial_x^{\beta}
\phi(\theta y+(1-\theta)x)(1-\theta)^{N-1}d \theta.
\end{array}
$$
Using
$$
(y-x)^{\beta}e^{i(x-y)\cdot\xi}=i^{|\beta|}\partial_{\xi}^{\beta}\left[e^{i(x-y)\cdot\xi}\right],
$$
and integrating by parts with respect to $\xi$ in the first term of $I$ in
\ref{eq3.3.6}
$$
\begin{array}{rcl}
I&=&\displaystyle N\left(1-\chi\left((2R)^{-1}|x|\right)\right)\int_0^1{\cal O}_s\int\int e^{i(x-y)\cdot\xi} 
\sum_{|\beta|=N}\frac{i^{-|\beta|}}{\beta!}
\partial_{\xi}^{\beta} b_h(x,\xi)\cdot\\
&\cdot&\displaystyle \partial_x^{\beta}\phi\left(\theta y+(1-\theta)x\right)(1-\theta)^{N-1} 
\partial_y^{\alpha}u(y)dyd\xi d\theta.
\end{array}
$$
             
              For each $\theta\in[0,1]$, the multiple symbol 
 $$
N\left(1-\chi\left((2R)^{-1}|x|\right)\right) \sum_{|\beta|=N}\frac{i^{-|\beta|}}{\beta!}
\partial_{\xi}^{\beta} b_h(x,\xi)\partial_x^{\beta}\phi\left(\theta
y+(1-\theta)x\right)(1-\theta)^{N-1},
$$
is in $S^{m-N}_{1,0}$ uniformly in $\theta$. Since $m-N+|\alpha|\leq 0$,
there exist $c>0$
and $N_1\in \BbbN$ such that (\ref{eq3.3.5}) holds with $E_ju$ replaced by I.

 In  II, note that $\chi((2R)^{-1}|x|) \partial_x^{\beta}\phi(x)=0$, so
that
 $$
 \begin{array}{rcl}
 II&=&\displaystyle(-1)^{|\alpha|}\chi\left((2R)^{-1}|x|\right)\cdot \\
 &\cdot&\displaystyle{\cal O}_s\int \int e^{i(x-y)\cdot\xi}
 \sum_{\beta\leq \alpha}
              \left(
              \begin{array}{c}
               \alpha\\ \beta
               \end{array}
               \right)
               (-i\xi)^{\beta}
 b_h(x,\xi)\partial_y^{\alpha-\beta}\phi(y)u(y)dyd\xi,
 \end{array}
 $$
 by integration by parts with respect to $y$. Therefore, it suffices to
consider the 
 pseudo-differential operator
 $\Psi_{II}$ given by
 $$
 \begin{array}{rcl}
 \displaystyle\Psi_{II}u(x)&=&\displaystyle\chi\left((2R)^{-1}|x|\right)\cdot
 \\
& \cdot&\displaystyle {\cal O}_s\int \int e^{i(x-y)\cdot\xi}\xi^{\beta}
 a\left(P(x,A_h\xi);x,\xi\right)\chi(|\xi|)\phi(y)u(y)dy d\xi,
 \end{array}
 $$
 for $\beta\leq \alpha$. Choose $M\in\BbbN$ such that $2M-m-|\alpha|\geq
n+1$. 
 Integrating by parts with respect to $\xi$ and using $|x-y|\geq |x|/2\geq
R>0$,
 \begin{equation}\label{eq3.3.7}
 \begin{array}{rcl}
\displaystyle \Psi_{II}u(x)&=&\displaystyle
               \chi\left((2R)^{-1}|x|\right)\\
 &\cdot&\displaystyle
               \int\int e^{i(x-y)\cdot\xi}\frac{(-1)^M}{|x-y|^{2M}}
 \Delta^M_{\xi}\left[\xi^{\beta}a\left(P(x,A_h\xi);x,\xi\right)\chi(|\xi|)\right]\phi(y)u(y)dy
d\xi,
\end{array}
\end{equation}
 where the integral converges absolutely. Indeed,
 $$
\begin{array}{l}
 \displaystyle\left|\Delta^M_{\xi}\left[\xi^{\beta}a\left(P(x,A_h\xi);x,\xi)\chi(|\xi|\right)\right]\right|
\leq c\langle x\rangle^{2M}\langle\xi\rangle^{-n-1}
 \max_{|\gamma_1+\gamma_2|\leq 2M}\left\|\langle \xi\rangle^{-m+|\gamma_2|}
  \p^{\gamma_1}_s\p_{\xi}^{\gamma_2}a\right\|_{L^{\infty}}.
\end{array}
  $$
  Therefore, letting
  $$
  |a|_{(2M)}= \max_{|\gamma_1+\gamma_2|\leq 2M}\left\|\langle
\xi\rangle^{-m+|\gamma_2|}
  \p^{\gamma_1}_s\p_{\xi}^{\gamma_2}a\right\|_{L^{\infty}},
  $$
  one has
  $$
  \left|\Psi_{II}u(x)\right|\leq c R^{n/2}\|\phi\|_{L^{\infty}} |a|_{(2M)}
\|u\|_{L^2}.
  $$
 Define $\chi_R(|x|)=\chi\left(|x|/R\right)$. Then  by Cauchy-Schwarz 
  $$
  \left\|(1-\chi_{c_0R})\Psi_{II}u\right\|_{L^2}\leq cR^n\|\phi\|_{L^{\infty}}
|a|_{(2M)} \|u\|_{L^2},
  $$
  where the constant $c_0$ will be fixed below (see (\ref{eq3.3.10})).
   Now we turn to the estimate of
 \begin{equation}\label{eq3.3.8}
   \|\chi_{c_0R}\Psi_{II}u\|_{L^2}.
   \end{equation}
We shall use that  
  \begin{equation}\label{eq3.3.9}
               \begin{array}{rcl}
\displaystyle\frac{1}{|x-y|^{2M}}&=&\displaystyle\frac{1}{\left||x|^2-2(x\cdot y)+|y|^2\right|^M}\\
\\
&=&\displaystyle
\frac{1}{|x|^{2M}}
\frac{1}{\left|1-2\left(\frac{x}{|x|}\cdot\frac{y}{|x|}\right)+\frac{|y|^2}{|x|^2}\right|^M}
\\
\\
&=&\displaystyle\frac{1}{|x|^{2M}}\left(1-\sum_{j=1}^{(2n+2)^M-1}P_j\left(\frac{x}{|x|}\right)Q_j
\left(\frac{y}{|x|}\right)\right)^{-1},
\end{array}
\end{equation}
where $P_j,\,Q_j$ are monomials of degree no bigger than $2M$ with $\deg Q_j\neq
0$ for $j=1,2,...,(2n+2)^M-1$.
Since $|x|>c_0R$ and $|y|<R$ it follows that
$$
\left|P_j\left(\frac{x}{|x|}\right)\right|\leq
a_j,\;\;\;\left|Q_j\left(\frac{y}{|x|}\right)\right|\leq 
b_j\left(\frac{|y|}{|x|}\right)^{deg Q_j}\leq \frac{b_j}{{c_0}^{deg Q_j}},
$$
so
$$
\left|
P_j\left(\frac{x}{|x|}\right) Q_j\left(\frac{y}{|x|}\right)\right|\leq
\frac{a_j b_j}{{c_0}^{deg Q_j}}.
$$
 We take 
 \begin{equation}\label{eq3.3.10}
 c_0\geq 2\,(2n+2)^M \,
 \max\left\{ |a_j b_j| \,:\,j=1,..,(2n+2)^M-1\right\}
 \end{equation}
 and rewrite (\ref{eq3.3.9}) as
 \begin{equation}\label{eq3.3.11}
               \begin{array}{l}
\displaystyle\frac{1}{|x|^{2M}}\frac{1}{1-\displaystyle\sum_j\,P_j\left(\frac{x}{|x|}\right) Q_j\left(\frac{y}{|x|}\right)}
 \\
\\
\qquad\displaystyle=\frac{1}{|x|^{2M}}\left(\sum_{k=0}^{\infty}\left(\sum_{j=1}^{(2n+2)^M-1}P_j\left(\frac{x}{|x|}\right)
 Q_j\left(\frac{y}{|x|}\right)\right)^k\right)\\
\\
 \qquad\displaystyle=\frac{1}{|x|^{2M}}\sum_{k=0}^{\infty}\,\sum_{j_1..j_k=1}^{(2n+2)^M-1}
 \frac{\displaystyle\prod_{i=1}^k
P_{j_i}\left(\frac{x}{|x|}\right)Q_{j_i}\left(\frac{y}{|x|}\right)
 \;|x|^{\deg Q_{j_i}}}{|x|^{\deg Q_{j_i}}}.
 \end{array}
 \end{equation}
 We observe that $\displaystyle|x|^{\deg Q_{j_i}}Q_{j_i}\left(\frac{y}{|x|}\right)$
depends just on $y$ and if $|y|<R$ then
 $$
 \left| |x|^{\deg Q_{j_i}}Q_{j_i}\left(\frac{y}{|x|}\right)\right|
\leq b_{j_i}R^{\deg Q_{j_i}}.
 $$
 So returning to the estimate (\ref{eq3.3.8}) and using the argument in (\ref{eq3.3.7}) it
follows 
 $$
 \begin{array}{l}
\displaystyle \left\|\chi_{c_0R}\Psi_{II}u\right\|_{L^2}=\displaystyle
\left\|\chi_{c_0R}\sum_{k=0}^{\infty}\sum_{j_1..j_k=1}^{(2n+2)^M-1}\prod_{i=1}^k 
 (-1)^M\frac{P_{j_i}\left(\frac{x}{|x|}\right)}{|x|^{\deg
Q_{j_i}}}\,\cdot\right.\\
\\
\qquad\cdot\displaystyle
 \left.\int \int e^{i(x-y)\cdot\xi}\frac{1}{|x|^{2M}}
\Delta_{\xi}^M\left(\xi^{\beta}a(\cdot;\cdot,\cdot)
 \chi(|\xi|)\right)Q_{j_i}\left(\frac{y}{|x|}\right)|x|^{\deg
Q_{j_i}}\phi(y)u(y)dyd\xi\displaystyle\right\|_{L^2}.
\end{array}
$$
               
               Using that 
$$ 
\tilde b(x,\xi)= \chi_{c_0R}(|x|)\frac{1}{|x|^{2M}}
\Delta_{\xi}^M\left(\xi^{\beta}a\left(P(x,A_h\xi);x,\xi\right)\chi(|\xi|)\right)  
$$ 
is a symbol which falls under the scope of Theorem \ref{th3.2.1} -see Remark \ref{re3.1.3}(c), it follows that
$$ 
\begin{array} {l}
\displaystyle\left\|\int e^{ix\cdot \xi}\tilde b(x,\xi\left)(\int e^{-iy\cdot
\xi}Q_{j_i}(y)\phi(y)u(y)dy\right)d\xi\right\|_{L^2_x}\\ 
\qquad\displaystyle\leq c\left\|\int e^{-iy\cdot
\xi}Q_{j_i}(y)\phi(y)u(y)dy\right\|_{L^2_{\xi}}\\ 
\qquad\displaystyle\leq cb_{j_i} R^{\deg
Q_{j_i}}\|u\|_{L^2}.  
\end{array}
$$ 
Also 
$$ 
\begin{array}{l}
\displaystyle\left\|\chi_{c_0R}\sum_{k=0}^{\infty}\sum_{j_1..j_k=1}^{(2n+2)^M-1}\prod_{i=1}^k
 (-1)^M\frac{P_{j_i}\left(\frac{x}{|x|}\right)}{|x|^{\deg
Q_{j_i}}}\right\|_{L^{\infty}}
\leq
\sum_{k=0}^{\infty}\sum_{j_1..j_k=1}^{(2n+2)^M-1}\prod_{i=1}^k\frac{a_{j_i}}{(c_0R)^{\deg
Q_{j_i}}}.
\end{array}
$$
               
                Combining the last two estimates we conclude that
$$
\begin{array}{l}
\displaystyle\left\|\chi_{c_0R}\Psi_{II}u\right\|_{L^2}\\
\qquad\displaystyle\leq c\sum_{k=0}^{\infty}\sum_{j_1..j_k=1}^{(2n+2)^M-1} \prod_{i=1}^k
 \frac{a_{j_i} b_{j_i} R^{\deg Q_{j_i}}}{(c_0R)^{\deg Q_{j_i}}}\|u\|_{L^2}
\leq c\sum_{k=0}^{\infty}\sum_{j_1..j_k=1}^{(2n+2)^M-1}
 \prod_{i=1}^k \frac{a_{j_i}b_{j_i}}{{c_0}^{\deg Q_{j_i}}}\|u\|_{L^2}\\
\\
\qquad\displaystyle\leq
c\sum_{k=0}^{\infty}\left(\sum_{j=1}^{(2n+2)^M-1}\frac{a_{j_i}b_{j_i}}{c_0}\right)^k
 \|u\|_{L^2}
 \leq \sum_{k=0}^{\infty} \frac{1}{2^k} \|u\|_{L^2}\leq c\|u\|_{L^2}
 \end{array}
  $$
and the proof of Theorem \ref{th3.3.4} is completed.  
  \bigskip
                
                \begin {theorem} \label{th3.3.5}
  
   Let $a\in {\cal S}(\BbbR^n:S^m_{1,0})$ and let $b_h$ be defined as
(\ref{eq3.1.9}), 
   $N\in\BbbN$, $N\geq m$. Let
   $$
   c(x,\xi)=\sum_{|\alpha|<N}\frac{i^{-|\alpha|}}{\alpha!}
   \p_x^{\alpha}\p_{\xi}^{\alpha}\bar b_h(x,\xi)
   $$
   and $E=\Psi_{b_h}^*-\Psi_c$. Let $\phi\in {\cal S}(\BbbR^n)$. Then
there
   exists
   $N_1=N_1(n)\in\BbbN$ such that for $u\in{\cal S}(\BbbR^n)$ 
   \begin{equation}\label{eq3.3.12}
                \begin{array}{rcl}
  \displaystyle\left\| \phi Eu\right\|_{L^2}&\leq&\displaystyle c \max_{|\beta|\leq N_1}\left\|\langle
x\rangle^{N_1}
   \p_x^{\alpha}\phi\right\|_{L^{\infty}}\cdot\\
  &\cdot&\displaystyle\max_{\beta_4+|\beta_1+\beta_2+\beta_3|\leq N_1}
  \left \|\langle s\rangle^{\beta_4}\langle \xi\rangle^{-m+|\beta_3|}
   \p_s^{\beta_1}\p_x^{\beta_2}\p_{\xi}^{\beta_3}a\right\|_{L^{\infty}}\|u\|_{L^2}
   \end{array}
   \end{equation}
   \end{theorem}
   
   {\sl Proof of Theorem \ref{th3.3.5}\rm}
   
    We shall prove the Theorem for $ \phi \in C^{\infty}_0(\BbbR^n)$ with
support
    in $B_R(0)$. Introducing a partition of unity and summing the
corresponding 
terms we get the desired result.

    Thus,
$$
\begin{array}{rcl}
\displaystyle\phi(x) Eu(x)&=&\displaystyle{\cal O}_s\int\int e^{i(x-y)\cdot\xi}\phi(x)\{\bar
b_h(y,\xi)-\sum_{|\alpha|<N}
\frac{i^{-|\alpha|}}{\alpha!}\p_x^{\alpha}\p_{\xi}^{\alpha}\bar
b_h(x,\xi)\}u(y)dyd\xi\\
&=&\displaystyle{\cal O}_s\int\int e^{i(x-y)\cdot\xi}\phi(x)\left(1-\chi((2R)^{-1}|y|)\right)\cdot\\
&\cdot&\displaystyle\left\{\bar b_h(y,\xi)-\sum_{|\alpha|<N}
\frac{i^{-|\alpha|}}{\alpha!}\p_x^{\alpha}\p_{\xi}^{\alpha}\bar
b_h(x,\xi)\right\}u(y)dyd\xi \\
&+&\displaystyle{\cal O}_s\int\int e^{i(x-y)\cdot\xi}\phi(x)\chi\left((2R)^{-1}|y|\right)\cdot\\
&\cdot&\displaystyle\left\{\bar b_h(y,\xi)-\sum_{|\alpha|<N}
\frac{i^{-|\alpha|}}{\alpha!}\p_x^{\alpha}\p_{\xi}^{\alpha}\bar
b_h(x,\xi)\right\}u(y)dyd\xi \\
&=&I+II.
\end{array}
$$

First consider $I$. By a Taylor expansion of order $N$,
     $$
     \begin{array}{rcl}
b_h(y,\xi)&=&\displaystyle\sum_{|\alpha|<N}\frac{1}{\alpha!}(y-x)^{\alpha}\p_x^{\alpha}b_h(x,\xi)\\
&+&\displaystyle\sum_{|\alpha|=N}\frac{N}{\alpha!}(y-x)^{\alpha}\int_0^1(1-\theta)^{N-1}\p_x^{\alpha}b_h
     \left(\theta y+(1-\theta)x,\xi\right)d\theta.
     \end{array}
     $$
     Since $(y-x)^\alpha
e^{i(x-y)\cdot\xi}=i^{|\alpha|}\p_{\xi}^{\alpha}e^{i(x-y)\cdot\xi}$,
integration by parts with respect
     to $\xi$ gives
     $$
     \begin{array}{l}
     \displaystyle I={\cal O}_s\int\int
e^{i(x-y)\cdot\xi}\phi(x)\left(1-\chi\left((2R)^{-1}|y|\right)\right) Ni^{-N} \\
     \displaystyle  \sum_{|\alpha|=N}\frac{1}{\alpha!}
\int_0^1(1-\theta)^{N-1}\p_x^{\alpha}\p_{\xi}^{\alpha}\bar b_h  
    \left (\theta y+(1-\theta)x,\xi\right)d\theta u(y)dyd\xi.
     \end{array}
     $$
     
      Because of the compact support  in $x$ and $y$, the multiple symbol
      $$
\displaystyle
                \phi(x)\left(1-\chi\left((2R)^{-1}|y|\right)\right) Ni^{-N}\sum_{|\alpha|=N}\frac{1}{\alpha!} 
      (1-\theta)^{N-1}\p_x^{\alpha}\p_{\xi}^{\alpha}\bar b_h\left(\theta
y+(1-\theta)x,\xi\right) $$
      is in $S^{m-N}_{1,0}\subset S^0_{1,0}$ uniformly in $\theta$, so an
estimate of the type (\ref{eq3.3.12})
      holds with $\phi Eu$ replaced by $I$. A factor  $\langle R\rangle^N$
appears due 
to the differentiation with respect to $\xi$ of the quantity \linebreak
      $P\left(\theta y+(1-\theta)x,A_h\xi\right)$.

  Next we consider $II$. Choose $M\in\BbbN$ such that
$m-2M\leq-n-1$. If 
        $$
        \phi(x)\chi(\left(2R)^{-1}|y|\right)\ne 0,
        $$
        then $|x-y|\geq R$ so by integration by parts with respect to
$\xi$
$$
\begin{array}{rcl}
II&=&\displaystyle\int\int
e^{i(x-y)\cdot\xi}\phi(x)\chi\left((2R)^{-1}|y|\right)\frac{(-1)^M}{|x-y|^{2M}}
\Delta_{\xi}^M \bar b_h(y,\xi)u(y)dyd\xi\\
&-&\displaystyle\int\int
e^{i(x-y)\cdot\xi}\phi(x)\chi\left((2R)^{-1}|y|\right)\frac{(-1)^M}{|x-y|^{2M}}\cdot\\
&&\displaystyle\sum_{|\alpha|<N}\frac{i^{-|\alpha|}}{\alpha!}\p_x^{\alpha}\Delta_{\xi}^M
\p_{\xi}^{\alpha}\bar b_h(x,\xi)u(y)dyd\xi,
\end{array}
$$
        where the integrals converge absolutely by the choice on $M$.
        The multiple symbol in the second term above is in
$S^{m-2M}_{1,0}\subset S^0_{1,0}$ because 
        of the factor $\phi(x)$, so the corresponding pseudo-differential
operator is
$L^2$-bounded as 
        in (\ref{eq3.3.12}).

Therefore, it remains to show the $L^2$-boundedness of the
operator
         $$
         \begin{array}{l}
         \displaystyle\bar \Psi_{II}u(x)=\int\int
e^{i(x-y)\cdot\xi}\phi(x)\chi\left((2R)^{-1}|y|\right)\\
                \qquad \displaystyle\cdot\frac{(-1)^M}{|x-y|^{2M}}
        \Delta_{\xi}^M\left[\bar
a\left(P(y,A_h\xi);y,\xi\right)\chi(|\xi|)\right]u(y)dyd\xi.
         \end{array}
         $$ 
          But $\bar \Psi_{II}$ is the adjoint of the operator $\Psi_{II}$
defined in (\ref{eq3.3.7}) in the proof of Theorem 
          \ref{th3.3.4} (with $\beta=0$ and $\phi$ replaced by $\bar \phi$). Since
$\Psi_{II}$ was there proved to be bounded
          in $L^2$, so is $\bar \Psi_{II}$, and the operator norms are
equal.
          
          This proves Theorem \ref{th3.3.5}.
          \bigskip

  \begin {theorem} \label{th3.3.6}
          
           Let $a_1\in {\cal S}\left(\BbbR^n:S^{m_1}_{1,0}\right)$, $a_2 \in {\cal
S}\left(\BbbR^n:S^{m_2}_{1,0}\right)$, and let $b_1, b_2$
            be the corresponding symbols given in (\ref{eq3.1.9}), Definition \ref{def3.1.2},
with $A=A_h$. Suppose $N\in\BbbN$,
            $N\geq m_1+m_2$ and let
            $$
c(x,\xi)=\sum_{|\alpha|<N}\frac{i^{-|\alpha|}}{\alpha!}\p_{\xi}^{\alpha}\left[b_1(x,\xi)\p_x^{\alpha}\bar
b_2(x,\xi)\right].
            $$
            Let $\phi\in{\cal S}(\BbbR^n)$. Then there exist $c=c(N)$ and
$N_1\in\BbbN$ such that
for any 
            $u\in{\cal S}(\BbbR^n)$ 
            \begin{equation}\label{eq3.3.13}
                \begin{array}{rcl}
             \displaystyle\left\|\phi(\Psi_{b_1}\Psi_{b_2}^*-\Psi_c)u\right\|_{L^2}&\leq&\displaystyle
\max_{|\beta|\leq N_1}\left\|\langle x\rangle^{N_1}
             \p_x^{\beta}\phi\right\|_{L^{\infty}}\cdot\\
            & \cdot&\displaystyle \max_{\beta_4+|\beta_1+\beta_2+\beta_3|\leq N_1}
             \left\| \langle s\rangle^{\beta_4}  \langle
\xi\rangle^{-m_1+|\beta_3|}
\p_s^{\beta_1}\p_x^{\beta_2}\p_{\xi}^{\beta_3}a_1\right\|_{L^{\infty}}\cdot\\
           & \cdot&\displaystyle \max_{\beta_4+|\beta_1+\beta_2+\beta_3|\leq N_1}
             \left\| \langle s\rangle^{\beta_4}  \langle
\xi\rangle^{-m_2+|\beta_3|}
          \p_s^{\beta_1}\p_x^{\beta_2}\p_{\xi}^{\beta_3}a_2\right\|_{L^{\infty}}
\|u\|_{L^2}.
          \end{array}
          \end{equation}
          \end{theorem}
                                        
         {\sl Proof of Theorem \ref{th3.3.6}\rm}

As in the proof of Theorem \ref{th3.3.5}
           we shall assume that $ \phi \in C^{\infty}_0(\BbbR^n)$ with
support
    in $B_R(0)$. Introducing a partition of unity and summing the
corresponding terms 
    we get the desired result.

     We have,
$$
\begin{array}{l}
 \displaystyle\phi(\Psi_{b_1}\Psi_{b_2}^*-\Psi_c)u(x)=\displaystyle 
                {\cal O}_s\int\int
e^{i(x-y)\cdot\xi}\phi(x) \left[b_1(x,\xi)\bar b_2(y,\xi)-c(x,\xi)\right]u(y)dyd\xi\\
\qquad=\displaystyle{\cal O}_s \int\int e^{i(x-y)\cdot\xi}\phi(x)\left(1-\chi((2R)^{-1}|y|)\right)
\left[b_1(x,\xi)\bar b_2(y,\xi)-c(x,\xi)\right]u(y)dyd\xi\\ 
                \qquad+\displaystyle{\cal O}_s \int\int
e^{i(x-y)\cdot\xi}\phi(x)\chi\left((2R)^{-1}|y|\right) \left[b_1(x,\xi)\bar
b_2(y,\xi)-c(x,\xi)\right]u(y)dyd\xi\\ 
               \qquad=I+II. \end{array}
 $$
      First we consider $I$. By a Taylor expansion of order $N$,
      $$
      \begin{array}{rcl}
b_2(y,\xi)&=&\displaystyle
               \sum_{|\alpha|<N}\frac{1}{\alpha!}(y-x)^{\alpha}\p_x^{\alpha}b_2(x,\xi)\\
&+&\displaystyle
               \sum_{|\alpha|=N}\frac{N}{\alpha!}(y-x)^{\alpha}\int_0^1(1-\theta)^{N-1}\p_x^{\beta}
      b_2\left(\theta y+(1-\theta)x,\xi\right)d\theta.
      \end{array}
      $$
      Since
$(y-x)^{\alpha}e^{i(x-y)\cdot\xi}=i^{-|\alpha|}\p_{\xi}^{\alpha}e^{i(x-y)\cdot\xi}$,
      integration by parts with respect to $\xi$ gives 
      $$
     \begin{array}{l}
    \displaystyle I=\int\int e^{i(x-y)\cdot\xi}\phi(x)\left(1-\chi\left((2R)^{-1}|y|\right)\right)\\
\displaystyle
               \sum_{|\alpha|=N}\frac{N}{\alpha!}
i^{-|\alpha|}\int_0^1(1-\theta)^{N-1}\p_{\xi}^{\alpha}\left[b_1(x,\xi)\p_x^{\alpha}\bar
b_2\left(\theta y+(1-\theta)x,\xi\right)\right]d\theta u(y)dyd\xi.
     \end{array}
     $$
     Because of the compact support in $x$ and $y$, the multiple symbol
 $$
\phi(x)(1-\chi((2R)^{-1}|y|))\sum_{|\alpha|=N}\frac{N}{\alpha!}i^{-|\alpha|} 
     \p_{\xi}^{\alpha}[b_1(x,\xi)\p_x^{\alpha}\bar b_2(\theta
y+(1-\theta)x,\xi)]
     $$
     is in $S^{m_1+m_2-N}_{1,0} \subset S^0_{1,0}$ with seminorms
uniformly bounded in $\theta\in[0,1]$, 
     so the estimate (\ref{eq3.3.13}) holds  with the left hand side replaced by
$\|I\|_{L^2_x}$. A factor 
     $ \langle R\rangle^{N_1} $ appears from differentiation with respect
to $\xi$ of the quantity
      $P\left(\theta y+(1-\theta)x,A_h\xi\right)$,  $\theta\in[0,1]$. 
      
       Next consider $II$. Note that if $ \phi(x)\chi\left((2R)^{-1}|y|\right)\ne 0$,
then $|x-y|\geq R$.
       Therefore, choosing $M\in\BbbN$ such that $m_1+m_2-2M\leq -n-1$
and integrating by parts with respect to $\xi$
       \begin{equation}\label{eq3.3.14}
       \begin{array}{rcl}
       II&=&\displaystyle
       \int \int
e^{i(x-y)\cdot\xi}\phi(x)\chi\left((2R)^{-1}|y|\right)\frac{(-1)^M}{|x-y|^{2M}}
       \Delta_{\xi}^M \left[b_1(x,\xi)\bar b_2(y,\xi)\right]u(y)dyd\xi\\
       &-&\displaystyle\int\int
e^{i(x-y)\cdot\xi}\phi(x)\chi\left((2R)^{-1}|y|\right)\frac{(-1)^M}{|x-y|^{2M}}
        \Delta_{\xi}^M c(x,\xi)u(y)dyd\xi,
        \end{array}
        \end{equation}
        where the integrals converge absolutely.
        
                         The multiple symbol in the second term in
(\ref{eq3.3.14}) is in $S^{-n-1}_{1,0}\subset
                          S^0_{1,0}$, so an estimate of the type(\ref{eq3.3.13}) holds.

 The adjoint of the operator corresponding to the first term in (\ref{eq3.3.14}) has
multiple symbol
$$
\bar \phi(y)\chi\left((2R)^{-1}|x|\right) \frac{(-1)^M}{|x-y|^{2M}} \Delta_{\xi}^M
\left[\bar b_1(y,\xi) b_2(x,\xi)\right].
$$
 Replacing $\bar \phi$ by $\phi$ and $\bar b_1(y,\xi)$ by $\xi^\beta$ one obtains a
symbol similar to that of the
operator
$\Psi_{II}$ in (\ref{eq3.3.7}) in the proof of Theorem \ref{th3.3.4}, so we just need to
sketch the proof.

Thus, we define
$$
\begin{array}{rcl}
\Psi u(x)&=&\displaystyle
                \int\int e^{i(x-y)\cdot \xi} \bar
\phi(y)\chi\left((2R)^{-1}|x|\right)\cdot
\\
&\cdot&\displaystyle
                \frac{(-1)^M}{|x-y|^{2M}}\left(\sum_{|\gamma_1|+|\gamma_2|=2M}c_{\gamma_1
\gamma_2}
\p_{\xi}^{\gamma_1} \bar b_1(y,\xi) \p_{\xi}^{\gamma_2}
b_2(x,\xi)\right) u(y)dyd\xi.
\end{array}
$$

Next consider $\tilde\phi\in C^\infty_0$ which vanishes outside the ball of radius $3R/2$ and 
such that
$\tilde\phi\phi=\phi$. Then define
$$
\tilde b_{\gamma_1}(y,\xi)=\tilde\phi(y)\langle \xi\rangle^{|\gamma_1|-m_1}\p_{\xi}^{\gamma_1}
\bar b_1(y,\xi),
$$
and
$$
\tilde b_{\gamma_2}(x,\xi)=\chi\left((2R)^{-1}|x|\right)\,|x|^{-2M}\langle
\xi\rangle^{m_1-|\gamma_1|}\p_{\xi}^{\gamma_2} b_2(x,\xi).
$$ 

Using
the notation introduced in the proof of Theorem \ref{th3.3.4}, see
(\ref{eq3.3.9})-(\ref{eq3.3.11}), we have
$$
\frac{1}{|x-y|^{2M}}=\frac{1}{|x|^{2M}}\sum_{k=0}^{\infty}\sum_{j_1..j_k=1}^{(2n+2)^M-1}\prod_{i=1}^k
P_{j_i}\left(\frac{x}{|x|}\right) Q_{j_i}\left(\frac{y}{|x|}\right). $$

Hence   as in 
(\ref{eq3.3.9})-(\ref{eq3.3.11}) we can reduce ourselves to study $\|\chi((c_0R)^{-1}|x|)\Psi u\|_{L^2}$ for $c_0$ large enough.
Then we have 
$$ 
                \begin{array}{l}
                \displaystyle
                \left \|\chi((c_0R)^{-1}|x|)\Psi u\right\|_{L^2}\\ 
               \qquad\displaystyle \leq
c\sum_{k=0}^{\infty}\sum_{j_1..j_k=1}^{(2n+2)^M-1}
\sum_{|\gamma_1|+|\gamma_2|=2M}\left\|\chi\left((c_0R)^{-1}|x|\right)
\prod_{i=1}^kP_{j_i}\left(\frac{x}{|x|}\right)\frac{1}{|x|^{\deg
Q_{j_i}}}\cdot\right.\\ 
               \qquad\displaystyle\left.
                \cdot \int\int e^{i(x-y)\cdot \xi} \phi(y)
Q_{j_i}(y)
\tilde b_{\gamma_1}(y,\xi) \tilde b_{\gamma_2}(x,\xi) u(y)dyd\xi\,\right\|_{L^2_x}. 
\end{array}
$$ 

Then from Remark \ref{re3.1.3} (b)
$$\begin{array}{l} 
\displaystyle
                \left\|\int\int e^{i(x-y)\cdot \xi} \phi(y)
Q_{j_i}(y)
\tilde b_{\gamma_1}(y,\xi) \tilde b_{\gamma_2}(x,\xi) u(y)dyd\xi\,\right\|_{L^2_x}\\
\qquad\displaystyle
                =\,\left\|\,\int e^{ix\cdot \xi}\,\tilde b_{\gamma_2}(x,\xi)\left(\int F(y,\xi)\,dy\right)\,d\xi\,\right\|_{L^2_x}\\
\qquad\displaystyle
                \leq \,\left\|\int F(y,\xi)\,dy\,\right\|_{L^2_{\xi}},
\end{array}
$$
with
$$F(y,\xi)=e^{-iy\cdot \xi} \phi(y)
Q_{j_i}(y)
\tilde b_{\gamma_1}(y,\xi)   u(y).
$$
Now observe that if
$$a(z,\eta)=\tilde b_{\gamma_1}(\eta,z),
$$
then
$$\left|\p^\alpha_z\p^\beta_\eta a(z,\eta)\right|\leq\frac{c_{\alpha\beta}}{\langle\eta\rangle^{|\beta|}}
R^{|\alpha|+|\beta|}.
$$
Hence  $a\in {\cal S}^ 0_{1,0}$ and therefore
$$\left\|\int F(y,\xi)\,dy\,\right\|_{L^2_{\xi}}\leq cR^{\deg Q_{j_i}+N}\|u\|_{L^2}
$$
for some large $N$ which just depends on the dimension.

Gathering the above information with the argument used in the proof of
Theorem \ref{th3.3.4} one completes the proof.
\bigskip

\newpage

\section {THE BICHARACTERISTIC FLOW}\label{sec4}

  In this section we study the bicharacteristic flows associated to the
second
 order 
 ultrahyperbolic variable coefficient operator ${\cal L}(x)$ and its truncated
version 
 ${\cal L}^R(x)$, see (\ref{eq4.1.12})-(\ref{eq4.1.13}) below. Assuming a non-trapping condition (basic
 assumption) for the flow generated by 
 ${\cal L}(x)$ we prove in subsection \ref{subseq4.1} that the bicharacteristic flow is \lq\lq uniformly
 non-trapping" with respect to the parameter
 $R$. 
The analysis of the flow is moredelicate when the trajectory is not outgoing\footnote{Using the homogeneity property (\ref{eq2.2.8}) in
Section \ref{sec2}, it is enough to consider the flow for $s\ge 0$ and $|\xi|=1$.}, -see Theorem \ref{th4.1.5} below. In that case we  prove
that outside a bounded ball (the same for all
$R$ large enough)  it behaves in the
$x$
 variable as the free flow, but just in dyadic annuli - see Theorem \ref{th4.1.5}( v) for a precise statement.
On the other hand in the outgoing case, the trajectories are in fact perturbations of the free ones, as is proved
by Craig, Kappeler, and Strauss in
\cite{4} whose arguments we follow. The end of subsection \ref{subseq4.1} is devoted to prove that the
non-trapping condition is stable under small perturbation in the coefficients.
  
 In subsection \ref{subseq4.2} we deduce
 several estimates  for the continuous dependence upon 
 the initial data of the flows
 associated to 
 the operators  ${\cal L}^R(x)$
 (with respect the parameter $R$). The arguments in \cite{4} are again very helpful. The
estimates in subsections \ref{subseq4.1} and \ref{subseq4.2} will be
 used in the next section 
 to deduce several properties of the integrating factor $K^R$.

 One of the main differences of the flow in the non-elliptic setting with
respect to the elliptic one
 studied in Section \ref{sec2} is
   that the Hamiltonian  in the elliptic case
$$
h_2(x,\xi)=\displaystyle\sum_{j,k}a_{jk}(x)\xi_k\xi_j 
$$
is
 preserved under the flow. So ellipticity gives
 the
 a priori estimate -see (\ref{eq2.2.7}) in Section \ref{sec2},
 $$
 \nu^{-2}\left|\xi_0\right|^2\leq \left|\Xi(s;x_0,\xi_0)\right|^2\leq \nu^2\left|\xi_0\right|^2.
 $$  
 In particular, one has that the bicharacteristic flow is globally
defined. This is also true in the ultra-hyperbolic case under the "asymptotic flatness" assumption,
but does not follow immediately.

 \subsection{Uniformly non-trapping flows}\label{subseq4.1}

 Let
 $$
 {\cal L}(x)=-\p_{x_j}a_{jk}(x)\p_{x_k},
 $$
 where $A(x)=(a_{jk}(x))$ is a real, symmetric, and non-degenerated
matrix, i.e.
 \begin{equation}\label{eq4.1.1}
 \exists\, \nu>0\;\;\;\;\;\forall \xi\in \BbbR^n\;\;\;\;\;
 \nu^{-1}|\xi|\leq |A(x)\xi|\leq \nu|\xi|.
 \end{equation}
 
 We will assume that there exists a constant coefficient
 operator
 \begin{equation}\label{eq4.1.2}
{\cal L}^0=a^0_{jk}\p^2_{x_jx_k},
 \end{equation}
 with
 \begin{equation}\label{eq4.1.3}
 a_{jk}(x)-a^0_{jk}\in{\cal S}(\BbbR^n),\;\;\;j,k=1,..,n.
 \end{equation}
  
  After a change of variable we can assume, without loss of generality,
that
 \begin{equation}\label{eq4.1.4}
 A_h=(a^0_{jk})= 
 \left(
 \begin{array}{cc}
 I_{k\times k}&0\\
 0&-I_{(n-k)\times (n-k)}
 \end{array}
 \right).
 \end{equation}

  We shall use the notation
 $\tilde \xi=A_h \xi.$
  The bicharacteristic flow is defined by
 \begin{equation}\label{eq4.1.5}
 \left\{
 \begin{array}{l}
 \displaystyle\frac{d}{ds} X_j(s;x_0,\xi_0)=
 2\sum_{k=1}^n a_{jk}\left(X(s;x_0,\xi_0)\right)\,\Xi_k(s;x_0,\xi_0),\\
 \\
 \displaystyle\frac{d}{ds} \Xi_j(s;x_0,\xi_0)=-
 \sum_{k,l=1}^n \p_{x_j}a_{kl}\left(X(s;x_0,\xi_0)\right) 
 \Xi_k(s;x_0,\xi_0) \Xi_l(s;x_0,\xi_0)\\
 \\
 (X(0;x_0,\xi_0),\Xi\left(0;x_0,\xi_0)\right)=(x_0,\xi_0).
 \end{array}
 \right.
 \end{equation}

  From the classical result of ode's we know that the bicharacteristic
flow
 exists in the time interval $s\in(-\delta,\delta)$ with 
 $\delta=\delta(x_0,\xi_0)$, and $\delta(\cdot)$ depending continuously on 
 $(x_0,\xi_0)$.

 By homogeneity of the symbol of ${\cal L}(x)$ (see (\ref{eq2.2.8}) in  Section \ref{sec2}) for any
 $t\in \BbbR$ one has that
 $$
 X(s;,x,t\xi)=X(ts;x,\xi),\;\;\;\;\Xi(s;x,t\xi)=t\Xi(ts:x,\xi).
 $$
  
 \vskip.2in
 \underline{BASIC ASSUMPTION } We shall assume that
${\cal L}(x)$ is
 non-trapping, i.e. for each $(x_0,\xi_0)\in\BbbR^n\times(\BbbR^n-\{0\})$ and for each $\mu>0$
 there exists  $s_0$ with $0<s_0<\delta$ such that\footnote{Using the homogeneity property of the hamiltonian flow
together with the methods given below in Theorem \ref{th4.1.5} and simple compactness and connectivity
arguments, one can show that the seemingly weaker assumption that there exists $s_0\in
(-\delta,\delta)$ such that
$$
 \left|X(s_0;x_0,\xi_0)\right|\geq \mu
 $$
implies our basic assumption.}
 $$ |X(s_0;x_0,\xi_0)|\geq \mu.
 $$

 Our first result is the following.

 \begin {proposition} \label{pro4.1.1}
  \begin{equation}\label{eq4.1.6}
 \frac{d}{ds}\left|X(s;x_0,\xi_0)\right|^2=4\langle
 X(s;x_0,\xi_0);A(X(s;x_0,\xi_0))\Xi(s;x_0,\xi_0)
 \rangle
 \end{equation}
 and
 \begin{equation}\label{eq4.1.7}
 \begin{array}{l}
\displaystyle\frac{d^2}{ds^2}\left|X(s;x_0,\xi_0)\right|^2=\displaystyle
 8\left|A\left(X\left(s;x_0,\xi_0\right)\right)\Xi\left(s;x_0,\xi_0\right)\right|^2\\
 \\
 \qquad+\displaystyle
  8\left\langle X\left(s;x_0,\xi_0\right);\nabla
 A\left(X(s;x_0,\xi_0)\right)A\left(X(s;x_0,\xi_0)\right)\Xi(s;x_0,\xi_0)
 \Xi(s;x_0,\xi_0)\right\rangle\\
 \\
 \qquad-\displaystyle
  4\left\langle X(s;x_0,\xi_0); A\left(X(s;x_0,\xi_0)\right)\nabla
 A\left(X(s;x_0,\xi_0)\right)\Xi(s;x_0,\xi_0)
 \Xi(s;x_0,\xi_0)\right\rangle.
 \end{array}
 \end{equation}
 \end{proposition}

 The proof of Proposition \ref{pro4.1.1} follows directly from (\ref{eq4.1.5}).

 \begin{lemma} \label{le4.1.2}

  There exists $M>0$ which depends only on a finite number of ${\cal
 S}(\BbbR^n)$-seminorms
 of $A(x)-A_h$ and on $\nu$ in (\ref{eq4.1.1}) such that if $|X(s;x_0,\xi_0)|>M$,
then
 \begin{equation}\label{eq4.1.8}
 \displaystyle\frac{d^2}{ds^2}\left|X(s;x_0,\xi_0)\right|^2\geq 4
 \left|A\left(X(s;x_0,\xi_0)\right)\Xi(s;x_0,\xi_0)\right|^2\geq 0.
 \end{equation}
 \end{lemma}
 
  The proof of Lemma \ref{le4.1.2}  follows by combining (\ref{eq4.1.7}) with (\ref{eq4.1.1}) and (\ref{eq4.1.3}).

\begin {lemma} \label{le4.1.3} For any $\tilde M\geq M$, $M$ as in Lemma \ref{le4.1.2}, there exist $s_{\tilde M}$,
$b$, $\tilde b$, and $c_1$ depending only on $\tilde M$ and $A$, such that for $|x_0|\leq \tilde
M$, and $1/2\leq |\xi_0|\leq 7/4$,

(i) 
$s_{\tilde M}\in (0,\delta(x_0,\xi_0))$;

  (ii)
$$
 \displaystyle\frac{d}{ds}\left|X(s;x_0,\xi_0)\right|^2\,{\bigg|_{s=s_{\tilde M}}}\geq b>0\,;
  $$

(iii)
 $$
 \displaystyle\frac{d^2}{ds^2}\left|X(s;x_0,\xi_0)\right|^2\,{\bigg|_{s=s_{\tilde M}}}\geq \tilde b>0\,;
  $$

(iv) for all $s$ with $|s|\leq s_{\tilde M}$,
$$0<c_1^{-1}\leq \left|\Xi(s;x_0,\xi_0)\right|\leq c_1<+\infty\,;
$$

(v) for $s\in(s_{\tilde M},\delta)$
, (ii) holds.
\end{lemma}

 {\sl Proof of Lemma \ref{le4.1.3}\rm}

   For convenience we introduce the notation
 \begin{equation}\label{eq4.1.9}
   N(s;x_0,\xi_0)=\left|X(s;x_0,\xi_0)\right|^2.
   \end{equation}

  Because of the homogeneity property of the flow we can assume $s>0$ and $s_0>0$. We fix
$(x_0,\xi_0)\in\BbbR^n\times(\BbbR^n-\{0\})$, and let 
 $$
 M_1=\max\{M;|x_0|\},
 $$
 where $M$ is as in Lemma \ref{le4.1.2}. By the non-trapping assumption
 $$
 \exists\, \bar s=\bar s(x_0,\xi_0)\;\;\;s.t.\;\;\;N(\bar s;x_0,\xi_0)\geq
 (M_1+1)^2.
 $$
 Also $N(0;x_0,\xi_0)=|x_0|^2\leq M_1^2$. Define
 $$
 s_1=\sup\left\{s\in[0,\bar s]\,:\,N(s;x_0,\xi_0)\leq M_1^2\right\}.
 $$
 Thus, $ s_1\in[0,\bar s),\;N(s_1;x_0,\xi_0)=M_1^2$, and for $
s\in(s_1,\bar s]$ 
 one has $ N(s;x_0,\xi_0)>M_1^2$. By the mean value theorem there exists 
 $ s_2\in (s_1,\bar s)$ such that $ N'(s_2;x_0,\xi_0)>0$. Moreover,   
 $ N(s_2;x_0,\xi_0)>M_1^2\geq M^2$. By continuity there exists a
neighborhood
 $U$ of 
 $ (x_0,\xi_0)$ such that
 \begin{equation}\label{eq4.1.10}
 N'\left(s_2;\bar x,\bar \xi\right)\geq N'\left(s_2;x_0,\xi_0\right)/2>0,\;\;\;
 N\left(s_2;\bar x_0,\bar \xi_0\right)>M^2,\;\;\;
 \forall(\bar x,\bar \xi)\in U.
\end{equation}
  
  \begin{claim} :\label{claim4.1.4}
  
  \begin{equation}\label{eq4.1.11}
N'\left(s;\bar x,\bar \xi\right)\geq
N'(s_2;\bar
 x,\bar \xi),
 \;\;\;\forall s\geq s_2\;
 \;\;\forall(\bar x,\bar \xi)\in U.
   \end{equation}
   \end{claim}
 
  Suppose not, so there exists $s_3>s_2,\;(x',\xi')\in U$, and $\beta>0$
 such that
 $$
 N'(s_3;x',\xi')<\beta<N'(s_2;x',\xi').
 $$
 Let
 $$
 s_4=\inf\left\{s\in[s_2,s_3]\,:\;N'(s;x',\xi')<\beta\right\},
 $$
 then $s_4\in(s_2,s_3]$, and $N'\left(s_4;x',\xi'\right)= \beta$. By the mean value
theorem,
 there exists 
 $s_5\in(s_2,s_4)$ such that 
 $$
 N''\left(s_5;x',\xi'\right)<0.
 $$

  But for $s\in[s_2,s_4),\;N'(s;x',\xi')\geq \beta$, and so
 $N\left(s_5;x',\xi'\right)>N\left(s_2;x',\xi'\right)\geq M^2$, 
 which contradicts Lemma \ref{le4.1.2}, so we have 
  established Claim \ref{claim4.1.4} .
 \vskip.15in 
  Next we cover $K$ by finitely many neighborhoods $U_r,\;r=1,..,N$, and
let
 $$
 a=\max_{r=1,..,N}\,s_2(U_r),\;\;\;\;\;(\mbox{see (\ref{eq4.1.10})-(\ref{eq4.1.11}))},
 $$
 and
 $$
 b=\min_{r=1,..,N}\;\left\{N'(s_2(U_r);\bar x,\bar \xi)\;:\;(\bar x,\bar
\xi)\in
 U_r\right\}.
 $$

 From (\ref{eq4.1.10}) it follows that $b>0$ and from Claim \ref{claim4.1.4} one has for $s\geq
a$
 that
 $$
 N'(s;x,\xi)\geq b,\;\;\;\forall (x,\xi)\in K.
 $$
 Hence only (iv) remains to be proved. The upper bound follows by compactness. For the lower bound
it suffices to see that $\Xi(s;x_0,\xi_0)\neq 0$ for $0<s<s_0.$ But this is a consequence of the
uniqueness of the flow (\ref{eq4.1.5}). 

The proof of Lemma \ref{le4.1.3} is completed. 
 \bigskip

  Let $\theta\in C^{\infty}_0(\BbbR^n),\;0\leq\theta(x)\leq 1$, 
  $supp\,\theta\subset\{x\,:\,|x|\leq 2\}$
  and $\theta(x)=1$ on $\{x\,:\,|x|\leq 1\}$. For each $1<R\leq \infty$ define
 \begin{equation}\label{eq4.1.12}
 \begin{array}{rcl}
 A^R(x)&=&
 \theta(x/R)\,A(x) + \left(1-\theta(x/R)\right)A_h\\
 &=&
 A_h+\theta(x/R)\left(A(x)-A_h\right)=a^R_{jk}(x),\quad \mbox{if}\quad R>1,\quad A^\infty=A;
 \end{array}
 \end{equation}
\begin{equation}\label{eq4.1.13}
{ \cal L}^R(x)=-\p_{x_j} a_{jk}^R(x)\p_{x_k}\quad \mbox{if}\quad R>1,\quad {\cal
L}^\infty={\cal L}.
 \end{equation}

  Notice that there exists $R'>1$ such that for $R>R'$ one has
 \begin{equation}\label{eq4.1.14}
 \nu^{-1}|\xi|/2\leq \left|A^R(x)\xi\right|\leq 2\nu|\xi|\;\;\;\mbox{and}\;\;
 \;A^R(x)-A_h\in {\cal S}(\BbbR^n)
 \end{equation}
 uniformly in R, so that Lemma \ref{le4.1.2} will apply uniformly in $R>R'$.
 
Also
 \begin{equation}\label{eq4.1.15}
 \p_{x_l}a^R_{jk}(x)=\theta(x/R)\,\p_{x_l}a_{jk}(x)+
\displaystyle \frac{\p_{x_l}\theta(x/R)}{R}
\left(
 a_{jk}(x)-a_{jk}^0\right),\;\;j,k=1,..,n,
 \end{equation}
 so for any $\tau>0$ there exists $c_{\tau}$ such that
 \begin{equation}\label{eq4.1.16}
 \left|\p_{x_l}a^R_{jk}(x)\right|\leq
\frac{c_{\tau}}{(1+|x|)^{\tau}},\;\;\;\;\;j,k=1,..,n.
 \end{equation}

 We will always consider $R\geq R_*=\max\{R'; 4M; 4\}$ with $M$ as in Lemma
\ref{le4.1.2}.

  Denote by $\left(X^R(s;x_0,\xi_0),\Xi^R(s;x_0,\xi_0)\right)$ the bicharacteristic 
  flow associated to the operators ${\cal L}^R$.
  
  \begin {theorem} \label{th4.1.5}

  There exist $c_1, c_2, c_3>0$ and $M_1>0$ (sufficiently large) such that
  for any $(x,\xi)\in\BbbR^n\times \BbbS^{n-1}$ and $R>2M_1$ one has 

 (i) the bicharacteristic flow 
 \begin{equation}\label{eq4.1.17}
 \left(X^R(s;x,\xi),\Xi^R(s;x,\xi)\right)\;\;\;\mbox{exists for any}\; s\in \BbbR,
 \end{equation}
 (ii)    
 \begin{equation}\label{eq4.1.18}
  0<c_1^{-1}\leq \left|\Xi^R(s;x,\xi)\right|\leq c_1,\;\;\;\;\forall s\in \BbbR,
 \end{equation}
 (iii) there exists $\tilde s_0\in\BbbR$ such that\footnote{We will say that the trajectory is
outgoing at $\tilde s_0$ when (\ref{eq4.1.19}) holds.}
 \begin{equation}\label{eq4.1.19}
\displaystyle \frac{d\;}{dt}\left|X^R(s;x,\xi)\right|^2\bigg|_{s=\tilde s_0}\geq 0
 \;\;\;\;\mbox{and}\;\;\;\;
 \left|X^R(\tilde s_0;x,\xi)\right|\geq M_1,
 \end{equation}

 (iv) define
 $$
 s_0= \inf \left\{\tilde s_0\geq 0\,:\; (\ref{eq4.1.19}) \;\mbox{holds}\right\},
 $$
 then
 \begin{equation}\label{eq4.1.20}
  \left|X^R(s;x,\xi)\right|^2\geq c_2(s-s_0)^2+M_1^2,\;\;\;\;s\geq s_0,
 \end{equation}
 (v) for  $ 0 <s\leq s_0$ and $k\in\{1,2,..\}$  the set
 \begin{equation}\label{eq4.1.21}
 I_k^R=\left\{s\in [0,s_0]\,:\,2^k\leq\left |X^R(s;x,\xi)\right|\leq 2^{k+1}\right\},
 \end{equation}
 verifies
 \begin{equation}\label{eq4.1.22}
 \left|I^R_k\right|\leq c_32^k.
 \end{equation}
 \end{theorem}
 
{\sl Proof of Theorem \ref{th4.1.5}\rm}

  Recall (\ref{eq4.1.5}), (\ref{eq4.1.14}), and (\ref{eq4.1.16}). Then given any $\tau>1$ we can always
assume, possibly by considering another
$\tilde\tau$ with $1<\tilde\tau<\tau$, that there is
$M_1=2^J$, $J\in \BbbN$, $2^{J(1-\tau)}<1/4$ such that if
 $$
 \left|X^R(\tilde s_0;x,\xi)\right|>M_1=2^J
 $$
 then
 \begin{equation}\label{eq4.1.23}
\displaystyle \frac{d^2\;}{ds^2}\left|X^R(\tilde s_0;x,\xi)\right|^2\geq 4\nu^{-2}\left|\Xi^R(\tilde
 s_0;x,\xi)\right|^2,
 \end{equation}
 and
 \begin{equation}\label{eq4.1.24}
\displaystyle \left|\frac{d\;}{ds}\left|\Xi^R(\tilde s_0;x,\xi)\right|\right|\leq
\frac{\nu^{-1}}{100\left|X^R(\tilde
 s_0;x,\xi)\right|^{\tau}}\left|\Xi^R(\tilde
 s_0;x,\xi)\right|^2.
 \end{equation}
 
  Assume as a first step that $|x|\leq  M_1$ and
$1/2\leq
 |\xi|\leq 7/4$. Since $R\geq M_1$ then $X^R(s;x,\xi)=X(s;x,\xi)$,
 $\Xi^R(s;x,\xi)=\Xi(s;x,\xi)$ as far
 as
 $|X(s;x,\xi)|\leq R$. Then from Lemma \ref{le4.1.3} 
 there exist $c_0>0$ and $s_{M_1}$ such that
 \begin{equation}\label{eq4.1.25}
 \left|X^R(s_{M_1};x,\xi)\right|\geq M_1, 
 \end{equation}
 \begin{equation}\label{eq4.1.26}
 c_0^{-1}\leq |\Xi^R(s;x,\xi)|\leq c_0,\;\;\;\;\;\;\mbox{if}\;\;\;0\leq
s\leq
 s_{M_1},
 \end{equation}
 and
 \begin{equation}\label{eq4.1.27}
\displaystyle\frac{d\;}{ds}\left|X^R(s;x,\xi)\right|^2>b>0,\;\;\; \;\;\;\mbox{if}\;\;\;s>s_{M_1},
 \end{equation}
as far as $|X(s;x,\xi|\leq R.$
 
  We will also use the inequalities -see (\ref{eq4.1.6}),
 \begin{equation}\label{eq4.1.28}
\displaystyle \left|\frac{d\;}{ds}\left|X^R(s;x,\xi)\right|^2\right|\leq
 4\nu\left|X^R(s;x,\xi)\right|\,\left|\Xi^R(s;x,\xi)\right|,
 \end{equation}
 and
 \begin{equation}\label{eq4.1.29}
\displaystyle1-w<\frac{1}{1+w}<\frac{1}{1-w}<1+3w/2,\;\;\;\;\mbox{if}\;\;\;w\in(0,1/4).
 \end{equation}

  Assume now for $j\in\BbbZ$ that $\,2^{j+J}\leq |x|\leq 2^{j+1+J}$ and
$$\displaystyle\frac{d\;}{ds}\left|X^R(s;x,\xi)\right|^2|_{s=0}<0.$$ 

 If $j<0$ we have already proved parts (ii), (iii) and (v) of the theorem.
Thus,
 consider $j=0,1,2,..$
 From (\ref{eq4.1.24}) and as long as
 $ |X^R(t;x,\xi)|>2^{j+J}$ for $t\in[0,s]$ one has
 \begin{equation}\label{eq4.1.30}
 \displaystyle\frac{|\xi|}{1+2^{-(j+J)\tau}|\xi|(s-0)/100\nu}\leq
\left|\Xi^R(s;x,\xi)\right|\leq
 \frac{|\xi|}{1-2^{-(j+J)\tau}|\xi|(s-0)/100\nu}.
 \end{equation}
 Thus, combining (\ref{eq4.1.29}) and (\ref{eq4.1.30}) and if $s<2^{(j+J)\tau}25\nu$,
 \begin{equation}\label{eq4.1.31}
 \displaystyle|\xi|\left(1-2^{-(j+J)\tau}|\xi|\frac{s}{100\nu}\right)\leq|\Xi^R(s;x,\xi)|
 \leq |\xi|\left(1+\frac{3}{2}\,2^{-(j+J)\tau}|\xi|\frac{s}{100\nu}\right).  
 \end{equation}
 So, in particular,
 \begin{equation}\label{eq4.1.32}
\displaystyle \frac{1}{2}|\xi|\leq\left|\Xi^R(s;x,\xi)\right|\leq \frac{7}{4}|\xi|,
 \end{equation}
 as long as  $\left |X^R(t;x,\xi)\right|>2^{j+J}$ for $t\in[0,s]$.

  From (\ref{eq4.1.23}), (\ref{eq4.1.28}), and (\ref{eq4.1.32}) we get (recall $\xi\in\BbbS^{n-1}$) that
 \begin{equation}\label{eq4.1.33}
\displaystyle \frac{d\;}{ds} \left|X^R(s;x,\xi)\right|^2\geq -7\nu \,2^{j+J+1} +2\nu^{-1}s,
 \end{equation}
  as long as $ \left|X^R(t;x,\xi)\right|>2^{j+J}$ for $t\in[0,s]$. 
Then either there is $\tilde s_{j+1}$, $0\leq \tilde s_{j+1}\leq 7\nu^22^{j+J+1}$ such that
$\displaystyle\left|X^R(\tilde s_{j+1};x,\xi\right|\leq2^{j+J}$ or for $0\leq s\leq7\nu^22^{j+J+1}$, $\left|X^R(
s;x,\xi\right|\geq2^{j+J}$. In the second case and from (\ref{eq4.1.33}) there exists $\tilde s_0$ such that (\ref{eq4.1.19})
holds with $0\leq \tilde s_0\leq 7\nu^22^{j+J+1}$. Therefore $s_0$ exists and from (\ref{eq4.1.28}) and (\ref{eq4.1.32})
we have that $\left|X(
s_0;x,\xi)\right|\leq(7\nu^3+1)2^{j+J}$. Therefore (v) holds. Assume now the first case. Define 
$$s_{j+1}=\inf
\left\{\, s \,:0\leq s\leq7\nu^22^{j+J+1}\,\, \mbox{such that}\,\,\, \left|X^R(s;x,\xi)\right|\leq 2^{j+J}\,\right\}.
$$ 
Hence $2^{j+J}\leq\left|X^R(s;x,\xi)\right|\leq(7\nu^3+1)2^{j+J+1}$ for $0\leq s\leq s_{j+1}$.
 Then from 
 (\ref{eq4.1.31}) we get
 $$
 |\xi|\left(1-2^{(j+J)(1-\tau)}\right)\leq \left|\Xi^R(s_{j+1};x,\xi)\right|\leq
 |\xi|\left(1+2^{(j+J)(1-\tau)}\right).
 $$
 
  We repeat the process changing $j$ by $j-1$ and $\xi$ by $
\Xi^R(s_{j+1};x,\xi)=:\xi_j$,
 and taking $s>s_{j+1}$.
 Then as before  either the trajectory becomes outgoing at some $\tilde s_0>s_{j+1}+7\nu
2^{j+J}$ and we obtain (v), or there is
$s_{j}$ such that
 $$
 \left|X^R(s_{j};x,\xi)\right|\leq
2^{j-1+J}\;\;\;\;\mbox{and}\;\;\;s_{j}-s_{j+1}\leq
 7\,\nu^2 2^{j+J},
 $$
and 
 $$
 |\xi_j|\left(1-2^{(j-1+J)(1-\tau)}\right)\leq \left|\Xi^R(s;x,\xi)\right|\leq
 |\xi_j|\left(1+2^{(j-1+J)(1-\tau)}\right),
 $$
for $s_{j+1}<s<s_{j}$.
We keep doing this until either the trajectory becomes outgoing or there are
$l=j,j-1,....,1,\,$
with
$$
\left |X^R(s_{l};x,\xi)\right|\leq
2^{l-1+J}\;\;\;\;\mbox{and}\;\;\;s_{l}-s_{l+1}\leq
 7\,\nu^2 2^{l+J},
 $$
and
 $$
 |\xi_l|\left(1-2^{(l-1+J)(1-\tau)}\right)\leq\left |\Xi^R(s;x,\xi)\right|\leq
 |\xi_l|\left(1+2^{(l-1+J)(1-\tau)}\right).
$$
Here $\xi_l=\Xi^R\left(s_{l+1};x,\xi\right)$ and $s_{l+1}<s<s_{l}$.

Assume this second possibility holds.
  From the condition $2^{J(1-\tau)}<1/4$ it follows that
 $$
 \begin{array}{l}
 \displaystyle\prod_{l=1}^j\left(1\pm2^{(l+J)(1-\tau)}\right)
 =\exp\left(\sum_{l=1}^j \ln \left(1\pm2^{(l+J)(1-\tau)}\right)\right)
 \end{array}
 $$
  Hence, we get 
 $$ 
 1/2\leq \left|\Xi^R(s_j;x,\xi)\right|\leq 7/4. 
 $$
  
Then we can apply the first step and obtain (\ref{eq4.1.25})-(\ref{eq4.1.27}). In particular
 from (\ref{eq4.1.25}) we know that $s_0\leq s_1+s_{M_1}$, and from  (\ref{eq4.1.26})
  $$
 c_0^{-1}\leq\left|\Xi^R(s_0;x,\xi)\right|\leq c_0.
$$
At this point we have proved (iii). By considering the cases $k+1\leq J$
and $k+1>J$ (v) also holds assuming either $\frac{d}{ds}\left|X^R(s,x,\xi)|^2\right|_{s=0}<0$ or $|x|\leq M_1$. Otherwise $s_0=0$ and
(v) is void.
 
It remains to prove (i), (ii), and (iv) for $s\geq s_0$. From (\ref{eq4.1.23}) we know that $|X^R(s,x,\xi)|$ is nondecreasing for
$s>s_0$. Assume that
$$\left|X^R(s_0,x,\xi)\right|\geq2^{L+J},
$$
for some $L=0,1,...$. Define $\xi_0=\Xi^R(s_0;x,\xi).$ We know that
$$\min\left\{c^{-1}_0, 1/2\right\}\leq|\xi_0|\leq\max\left\{c_0,7/4\right\}.
$$
Following the same argument as in (\ref{eq4.1.31}) we know that $\left(X^R(s;x,\xi);\Xi^R(s;x,\xi)\right)$ is defined for $s\leq
s_0+\frac{25\nu}{|\xi_0|}2^{J+L}=s_1$ and from (\ref{eq4.1.23})
$$\left|X^R(s_0;x,\xi)\right|\geq2^{L+J+1};
$$
$$\displaystyle|\xi_0|\left(1-2^{-(L+J)\tau}|\xi_0|\frac{s-s_0}{100\nu}\right)\leq|\Xi^R(s;x,\xi)|\leq
\left|\xi_0\right|\left(1+2^{-(L+J)\tau}|\xi_0|\frac{s-s_0}{100\nu}\right),
$$
and therefore for $s_0\leq s\leq s_1$
$$\displaystyle\frac{|\xi_0|}{2}\leq\left|\Xi^R(s;x,\xi)\right|\leq\frac{7}{4}\left|\xi_0\right|.
$$
Then we repeat the process and construct $\{s_k\}$, $k=1,2,...$ such that
$$s_k=s_{k-1}+\frac{25\nu}{|\xi_0|}2^{J+L+k-1}.
$$
Hence $\left(X^R(s;x,\xi);\Xi^R(s;x,\xi)\right)$ is defined for $0<s<
s_k$ and
$$\left|X^R(s_k;x,\xi)\right|\geq2^{L+J+k};
$$
$$|\xi_k|\left(1-2^{(k+L+J)(1-\tau)}\right)\leq\left|\Xi^R(s;x,\xi)\right|\leq
|\xi_k|\left(1+2^{(k+L+J)(1-\tau)}\right).
$$
Then reasoning as in the  case $s<s_0$ we get that
$$1/2|\xi_0|\leq\left|\Xi^R(s;x,\xi)\right|\leq 7/4|\xi_0|,
$$
for $s_0<s<s_k.$ This proves (i) and (ii). Finally (iv) follows from (i) and (\ref{eq4.1.23}).
The proof of Theorem \ref{th4.1.5} is complete.
\bigskip
  
\begin {lemma} \label{le4.1.6} For any $\psi$ such that 
$$
\displaystyle|\psi(x)|\leq \frac{c_{\tau}}{\langle x\rangle^{\tau}},\;\;\;\tau>1
$$
there exists $\tilde c$ such that
$$
\sup_{s, x, |\xi|=1}\int_0^s \psi\left(X^R(r;x,\xi)\right)dr< \tilde c,
$$
where the constant $\tilde c$ depends just on $\tau$ and $c_{\tau}$.
\end{lemma}

{\sl Proof of Lemma \ref{le4.1.6}\rm}

 It follows from Theorem \ref{th4.1.5}.

\bigskip

\begin {corollary} \label{co4.1.7} Let $A_0(x)$ be as in (\ref{eq4.1.1})-(\ref{eq4.1.4}). Let $B_1(x)$ be an $n\times n$ real matrix
with entries in ${\cal S}(\BbbR^n)$, and define $A_1(x)=A_0(x)+\epsilon B_1(x)$, 
where $\epsilon$ is chosen so that for all $\xi\in \BbbR^n$ $(2\nu)^{-1}\leq|A_1(x)\xi|\leq 2\nu|\xi|$.

Then, there exists $\epsilon_0>0,$ such that $A_1$ 
verifies the basic assumption (i.e. $A_1$ is non-trapping).
Here $\epsilon_0$ depends on $M_1$, $s_{M_1}$, $\nu$, $B_1$ and
 on a finite number of seminorms of the difference (\ref{eq4.1.3}).
\end{corollary}
In order to establish Corollary \ref{co4.1.7}, we need the following elementary o.d.e. lemma.

\begin{lemma} \label{le4.1.8} Let $\vec y_0(s)$ verify
$$\left\{
\begin{array}{rcl}
\displaystyle \frac{d}{ds}\vec y_0(s)&=&\vec f_0(\vec y_0(s)) \qquad \mbox{for} \quad 0<s<T\\
\vec y_0(0)&=&\vec z_0(s).
\end{array}
\right.
$$
Let $\tilde M=\sup_{0<s<T}|\vec y_0(s)|$, and suppose that $\vec f_1$ is given, 
with 
$$K=\sup_{|\vec y|\leq \tilde M}\left|\p \vec y \vec f_1(\vec y)\right|,$$
and let $\epsilon_0>0$ be given. Then there exists 
$ \tilde \epsilon=\tilde \epsilon(\epsilon_0, T, K),$ such
that, if 
\begin{equation}\label{eq4.1.34}
 \sup_{0\leq t\leq T}\left|\vec f_0(\vec y_0(t))-\vec f_1(\vec y_0(t))\right|\leq \tilde \epsilon,\
 \end{equation}
then there exists a unique solution $\vec y_1$ to
$$\left\{
 \begin{array}{rcl}
 \displaystyle \frac{d}{ds}\vec y_1(s)&=&\vec f_1\left(\vec y_1(s)\right) \qquad \mbox{for} \quad 0<s<T\\
\vec y_1(0)&=&\vec z_0(s).
\end{array}
 \right.
$$
Moreover  
$$ \sup_{0\leq s\leq T}\left|\vec y_0(s)-\vec y_1(s)\right|\leq  \epsilon_0,
$$ 
\end{lemma}

The proof is elementary and will be ommitted.

{\sl Proof of Corollary \ref{co4.1.7}\rm}
\bigskip

Assume that $\xi\in \BbbS^{n-1},$ $ |x|\leq 2M_1$. By Theorem \ref{th4.1.5}, 
there exist $s_{2M_1}$ and $c_0,$ such
that \linebreak $\left|X^0(s_{2M_1}; x, \xi\right|\geq 4 M_1,$ and such that 
$c_0^{-1}\leq \left|\Xi^0(s;x,\xi)\right|\leq c_0$
for all $s$. Let now $\tilde M=\sup |(X^0,\Xi ^0)|,$ where the supremum is taken over all 
$0\leq s\leq s_{2M_1},$ $|\xi|=1,$ and $|x|\leq 2M_1.$ Let $K$ now be defined as in Lemma \ref{le4.1.8}
for the system (\ref{eq4.1.5})  associated to $A_1.$ For $\epsilon_0$ to be chosen, 
$T=s_{2M_1},$ let $\tilde \epsilon$ be as in Lemma \ref{le4.1.8}. Choose now $\epsilon$ so small that,
in the terminology of Lemma \ref{le4.1.8}, for $\vec f_0,$ $\vec f_1,$ the system 
(\ref{eq4.1.5}) associated to $A_0,$ A$_1,$ (\ref{eq4.1.34}) holds. Then, if $\epsilon_0$ in Lemma \ref{le4.1.8} is
chosen small enough, we can conclude that there exists $0<\bar s<T$ such that
$\left|X^1\left(s_{2M_1}; x, \xi\right)\right|^2\geq 9 M_1^2;$ 
$(2c_0)^{-1}\leq \left|\,\Xi^1(\bar s;x,\xi)\,\right|\leq 2c_0;$ 
$\displaystyle\frac{d}{ds}\left|X^1(s;x,\xi)\right|^2\bigg|_{s=\bar s}>0;$ and
$\displaystyle\frac{d^2}{ds^2}\left|X^1(s;x,\xi)\right|^2\bigg|_{s=\bar s}\geq2\nu^{-1}\left|A_1(X^1)\cdot \Xi^1\right|,$
as long as $|X^1|\geq M_1.$ Inserting this information in 
the proof of Theorem \ref{th4.1.5} we obtain the conclusions of Thoerem \ref{th4.1.5} for $A_1.$

\subsection{The Continuous Dependence}\label{subseq4.2}

  In this subsection we shall deduce estimates concerning the continuous dependence of the flow associated
to the truncated operator ${\cal L}^R(x)$ with respect to the initial value. These estimates will
be given  in Theorem \ref{th4.2.1} as a function of the parameter $R$.  

 To simplify the notation we shall omit the sub-indexes denoting component or coordinates. 
  So instead of

\begin{equation}\label{eq4.2.1}
\left\{
\begin{array}{l}
\displaystyle\frac{d\,}{ds}X^R_j(s;x,\xi)=2\sum_{k=1}^na^R_{jk}\left(X^R(s;x,\xi)\right)\Xi^R_k(s;x,\xi),\\
\displaystyle\frac{d\,}{ds}\Xi^R_j(s;x,\xi)=-\sum_{k,l=1}^n\p_{x_j}a^R_{kl}\left(X^R(s;x,\xi)\right)\Xi^R_k(s;x,\xi)\Xi^R_l(s;x,\xi),
\end{array}
\right.
\end{equation}
with $a^R_{jk}(\cdot)$ defined in (\ref{eq4.1.12})-(\ref{eq4.1.13}), we shall write
\begin{equation}\label{eq4.2.2}
\left\{
\begin{array}{l}
\displaystyle\frac{d\,}{ds}X^R(s;x,\xi)=2a^R\left(X^R(s;x,\xi)\right)\Xi^R(s;x,\xi),\\
\displaystyle\frac{d\,}{ds}\Xi^R(s;x,\xi)=-\p a^R\left(X^R(s;x,\xi)\right)\Xi^R(s;x,\xi)\Xi^R(s;x,\xi).
\end{array}
\right.
\end{equation}

 Also, we shall omit the variables. 

Hence $\p_{\xi}X^R(\cdot)= \p_{\xi}X^R(s;x,\xi)$, $\p_{\xi}\Xi^R(\cdot)= \p_{\xi}\Xi^R(s;x,\xi)$
 satisfy the system
\begin{equation}\label{eq4.2.3}
\left\{
\begin{array}{rcl}
\displaystyle\frac{d\,}{ds}(\p_{\xi}X^R(\cdot))&=&2\p a^R\left(X^R(\cdot)\right)\left(\p_{\xi}X^R(\cdot)\right)\Xi^R(\cdot)
+2a^R\left(X^R(\cdot)\right)\left(\p_{\xi}\Xi^R(\cdot)\right),\\
\displaystyle\frac{d\,}{ds}\left(\p_{\xi}\Xi^R(\cdot)\right)&=&\displaystyle-\p^2 a^R\left(X^R(\cdot)\right)\left(\p_{\xi}X^R(\cdot)\right)\Xi^R(\cdot)\Xi^R(\cdot)
\\
&-&\displaystyle2\p a^R\left(X^R(\cdot)\right)\Xi^R(\cdot)\left(\p_{\xi}\Xi^R(\cdot)\right).
\end{array}
\right.
\end{equation}

 We observe that  $\p_xX^R(\cdot)= \p_xX^R(s;x,\xi)$, 
 $\p_x\Xi^R(\cdot)= \p_x\Xi^R(s;x,\xi)$  satisfy the system
obtained from (\ref{eq4.2.3}) by substituting $\p_{\xi}$ by $\p_x$ everywhere. 

As
a final simplification $c$ 
will denote in what follows a generic constant independent of $R$ which can change from line to line.  

 By homogeneity of the symbol of ${\cal L}^R(x)$ (see (\ref{eq2.2.8}) in Section \ref{sec2}) for any $R>0$,
\begin{equation}\label{eq4.2.4}
X^R(s;,x,t\xi)=X^R(ts;x,\xi),\;\;\;\;\Xi^R(s;x,t\xi)=t\Xi^R(ts:x,\xi),
\end{equation}
and consequently
\begin{equation}\label{eq4.2.5}
 \begin{array}{l}
\displaystyle
 \left(\p_x^{\alpha}\p_{\xi}^{\beta}X^R\right)(s;x,t\xi)=
 t^{-|\beta|}\left(\p_x^{\alpha}\p_{\xi}^{\beta}X^R\right)(ts;x,\xi),\\
\displaystyle\left(\p_x^{\alpha}\p_{\xi}^{\beta}\Xi^R\right)(s;x,t\xi)=t^{1-|\beta|}\left(\p_x^{\alpha}\p_{\xi}^{\beta}\Xi^R\right)(ts;x,\xi).
\end{array}
\end{equation}
 So we can take $|\xi|=1$ and consider $s\geq 0$. 
 Therefore
combining (\ref{eq4.1.18}) of Theorem \ref{th4.1.5} and (\ref{eq4.2.3}) it can be deduced that
$$\left\{
 \begin{array}{rcl}
\displaystyle\frac{d\,}{ds}\left|\p_{\xi}X^R(\cdot)\right|&\leq & c\left|\p a^R(X^R(\cdot))\right|\,\left|\p_{\xi}X^R(\cdot)\right|
+2\left|a^R\left(X^R(\cdot)\right)\right|\,\left|\p_{\xi}\Xi^R(\cdot)\right|,\\ 
\displaystyle\frac{d\,}{ds}\left|\p_{\xi}\Xi^R(\cdot)\right|&\leq &c\left|\p^2 a^R\left(X^R(\cdot)\right)\right|\,\left|\p_{\xi}X^R(\cdot)\right|
+c\left|\p a^R\left(X^R(\cdot)\right)\right|\,\left|\p_{\xi}\Xi^R(\cdot)\right|.
\end{array}
\right.
$$

  From our hypothesis on the decay of $\p a(\cdot)$ and Lemma \ref{le4.1.6} we can define
\begin{equation}\label{eq4.2.6}
  \left\{
  \begin{array}{rcl}
f(s)&=&\displaystyle
 \left |\p_{\xi}X^R(s;x,\xi)\right|\,\exp\left(c\displaystyle
 \int_0^s|\p a^R(X^R(s;x,\xi))|ds\right)\sim \left|\p_{\xi}X^R(s;x,\xi)\right|,\\
g(s)&=&\displaystyle
 \left |\p_{\xi}\Xi^R(s;x,\xi)\right|\,\exp\left(c\displaystyle
 \int_0^s|\p a^R(X^R(s;x,\xi))|ds\right)
\sim \left|\p_{\xi}\Xi^R(s;x,\xi)\right|.
\end{array}
\right.
  \end{equation}
  From (\ref{eq4.2.5}) it follows that
\begin{equation}\label{eq4.2.7}
  \left\{
  \begin{array}{l}
\displaystyle
  f'(s)\leq c\left|a^R(X^R(s;x,\xi))\right|g(s)\leq  cg(s),\\
\displaystyle
  g'(s)\leq c\left|\p^2a^R(X^R(s;x,\xi))\right|f(s)\leq c\left|\p^2a^R(X^R(s;\cdot))\right|f(s).
\end{array}
\right.
  \end{equation}

We observe that 
if 
 $|x|>R$ then $\left(X^R(s;x,\xi),\Xi^R(s;x,\xi)\right)=\left(x+2s\tilde \xi,\xi\right)$ as long  as
$|x+2s\tilde \xi|>R$, where $A_h\xi=\tilde \xi$.

\vskip.1in

\underline{Case 1}

 Assume $|x|>R$ with $|\xi|=1$.

 We shall assume that there exists a first $s_1>0$ such that $\left|X^R(s_1;x,\xi)\right|=R$, otherwise \linebreak
$\left(X^R(s;x,\xi),\Xi^R(s;x,\xi)\right)=(x+2s\tilde \xi,\xi)$, for all $s>0$. One has that 
$s_1\leq |x|+R$. 
 
  Since for $ s\in (0,s_1),$ 
$$\left(X^R(s;x,\xi), 
\Xi^R(s;x,\xi)\right)=(x+2s\tilde \xi,\xi),$$
then
$\left|\p_{\xi}X^R(s;x,\xi)\right|\leq 2s$ and $\left|\p_{\xi}\Xi^R(s;x,\xi)\right|\leq 1.$ 
 The case $s\geq s_1$ reduces to our next step.
\vskip.1in
\underline{Case 2 ($s<s_0$)}

 If $|x|\leq M_1$ , $M_1$ as in Theorem \ref{th4.1.5}, the analysis is given in case 3. Assume $M_1\leq|x|\leq R,$  $1/2<|\xi|\leq 3/2$
and 
 
\begin{equation}\label{eq4.2.8}
\displaystyle\frac{d\;}{ds}\left|X^R(s;x,\xi)\right|\bigg|_{s=0}\leq 0.
\end{equation}

 We also assume that (\ref{eq4.2.8}) holds for $s\in(0,s_0)$ with $s_0$ 
defined in Theorem \ref{th4.1.5} and therefore comparable in size to $R$
because of (v) of that theorem. Otherwise we would have reached
 the outgoing situation
 (case 4) in an intermediate step.

 Consider the majorized system (\ref{eq4.2.7})
with data $(f(0),g(0))=(a,b),\;a,b>0$.

Integrating we have
$$
g(s)\leq b+ c\int_0^s\left|\p^2a\left(X(\theta;x,\xi)\right)\right|f(\theta)d\theta,
$$so
\begin{equation}\label{eq4.2.9}
  \begin{array}{rcl}
f(s)&\leq&\displaystyle  a+bs+c\int_0^s\left(\int_0^l\left|\p^2a(X(\theta;x,\xi))\right|f(\theta)d\theta\right)dl\\
&\leq &\displaystyle  a+bs+c \int_0^s (s-\theta)\left|\p^2a\left(X(\theta;x,\xi)\right)\right|f(\theta)d\theta.
\end{array}
\end{equation}
  
  Let $s_2=\min\,\left\{s\in(0,c_0R)\,:\,f(s)=2(a+bs)\right\}$. 
\vskip.1in
\underline {Claim} $\left|X^R(s_2;x,\xi)\right|\leq R/2$. 
 Otherwise,  by (\ref{eq4.2.9}) one should have
$$
\displaystyle a+bs_2\leq c\int^{s_2}_0(s_2-\theta)\left|\p^2a(X(\theta;x,\xi))\right|(a+b\theta)d\theta,
$$
so
$$
c\,s_2^2 \frac{2^{\tau+2}}{R^{\tau+2}}\geq \frac{1}{2},
$$
which is a contradiction for $R\geq R_0$ with $R_0$ sufficiently large.

 We repeat the argument assuming that (\ref{eq4.2.8}) holds in the interval $(0,s_j)$ defining
$$
s_{j+1}=\min\,\left\{s\in(s_j,c_0R)\,:\,f(s)=2^j(a+bs)\right\}.
$$

\underline {Claim} $\left |X^R(s_j;x,\xi)\right|\leq R/2^j$.
 Otherwise, we would have
$$
\displaystyle 2c\left(s_{j+1}-s_j\right)^2\frac{2^{(\tau +2)j}}{R^{(\tau+2)}}\geq 1,
$$
and consequently, from (\ref{eq4.1.22}) in Theorem \ref{th4.1.5},  
$$
\displaystyle c\frac{R}{2^j}\geq s_{j+1}-s_j\geq \sqrt{\frac{R^{\tau+2}}{2c2^{(\tau+2) j}}},
$$
 which is a contradiction if $ R\geq 10 c_{\tau} 2^j$.

So we can repeat the argument $k$-times until $R\sim 2^k$, with $f(s_k)=2^k(a+bs_k)\sim R(a+bR)$
since by Theorem \ref{th4.1.5} (v) we have that $s_k\sim R.$

Similarly, one gets that $g(s_k)\leq  R(a+bR)$ and  
$\left|X^R(s_k;x,\xi)\right|\leq \tilde M$, with $\tilde M$ independent of $R$. Restarting
the variable $s$ we are led to the following case.
 
\vskip.1in
 \underline {Case 3}

Assume $ \left|X^R(0;x,\xi)\right|=\left|X(0;x,\xi)\right|\leq \tilde M$. 

 We consider the majorized system (\ref{eq4.2.7}) with data $(a_1,b_1)$. From Theorem \ref{th4.1.5} 
there exists $s^*>0$ (independently of $R$)
such that
$$
\displaystyle \left|X^R(s^*;x,\xi)\right|=\tilde M +1,\;\;\;\;\mbox{with}\;\;\;\frac{d\;}{ds}\left|X^R(s;x,\xi)\right|\bigg|_{s=s^*}>0.
$$
 Integrating the system (\ref{eq4.2.7}) we find that
$$
f(s) \leq a_1e^{cs},\;\;\;g(s)\leq b_1e^{cs}\;\;\;\mbox{for any}\;s\in[0,s^*).
$$
After restarting $s$ we are reduced to the following case.
\vskip.1in
\underline {Case 4} (outgoing, i.e. (\ref{eq4.1.20}) holds and $\displaystyle\frac{d}{ds}|X^R(s;x,\xi)|>0$).

Assume
$$
\displaystyle\left|X^R(0;x,\xi)\right|=|x|\geq \tilde M+1,\;\;\;\;\mbox{with}\;\;\;\frac{d\;}{ds}\left|X^R(s;x,\xi)\right|\bigg|_{s=0}>0.
$$ 
We consider the majorized system (\ref{eq4.2.7}) with data  $(a_2,b_2)$.

 Define
$$
\displaystyle h(s)=\frac{f(s)}{(1+s^2)^2}+\frac{g(s)}{1+s^2}.
$$
Thus,
$$
\begin{array}{rcl}
h'(s)&\leq&\displaystyle  \frac{d_0}{1+s^2}\,\frac{g(s)}{1+s^2}+ \frac{d_1}{1+s^2}(1+s^2)^2\left|\p^2a(X^R(s;x,\xi))\right|\frac{f(s)}{(1+s^2)^2}\\
&\leq &\displaystyle\frac{c}{1+s^2} h(s).
\end{array}
$$
  
  Hence, from Theorem \ref{th4.1.5} one has
$$
\begin{array}{l}
\displaystyle
  \left|\p_{\xi}X^R(s;x,\xi)\right|\leq f(s)\leq c\left(1+s^2\right)^2(a_2+b_2)\leq c\left(1+R^2\right)^2(a_2+b_2),\\
\displaystyle
  \left|\p_{\xi}\Xi^R(s;x,\xi)\right|\leq g(s)\leq c\left(1+s^2\right)^2(a_2+b_2)\leq c\left(1+R^2\right)^2(a_2+b_2),
\end{array}
$$
as far as $\left|X^R(s;x,\xi)\right|\leq R$.

\vskip.1in
\underline {Case 5}

Assume $\left|X^R(0;x,\xi)\right|=|x|\geq R$ with $\displaystyle\frac{d\;}{ds}\left|X^R(s;x,\xi)\right|\bigg|_{s=0}>0$.

Also we have the initial values $\displaystyle
  \left(\left|\p_{\xi}X^R(0;x,\xi)\right|,\left|\p_{\xi}\Xi^R(0;x,\xi)\right|\right)=(a_3,b_3)$.

So we have the solution
$$
\left|\p_{\xi}X^R(s;x,\xi)\right|\leq a_3+cs,\;\;\;\;\left|\p_{\xi}\Xi^R(s;x,\xi)\right|\leq b_3.
$$

 Collecting the information above we get
\begin{equation}\label{eq4.2.10}
\left|\p_{\xi}X^R(s;x,\xi)\right|,\,\left|\p_{\xi}\Xi^R(s;x,\xi)\right|\leq c(1+|s|)R^6,
\end{equation}
 for $R\geq R_0$.

  To estimate $\left(\p_xX^R(s;x,\xi),\p_x\Xi^R(s;x,\xi)\right)$ as above we observe that it satisfies 
  the system obtained from (\ref{eq4.2.3}) by substituting $\p_x$ by $\p_{\xi}$ everywhere. As in the previous case,
  we majorize $\left(\left|\p_{x}X^R(s;x,\xi)\right|,\,\left|\p_{x}\Xi^R(s;x,\xi)\right|\right)$ by the system in (\ref{eq4.2.7}). So the same argument shows that
$$
\left|\p_xX^R(s;x,\xi)\right|,\,\left|\p_x\Xi^R(s;x,\xi)\right|\leq c(1+|s|)R^6,
$$
 for $R\geq R_0$.

    To estimate the higher order derivatives we first observe that if
   $$
   \left(h^{R\alpha\beta}_1(s),h^{R\alpha\beta}_2(s)\right)
   =\left(\p_x^{\alpha}\p_{\xi}^{\beta}X^R(s;x,\xi), \p_x^{\alpha}\p_{\xi}^{\beta}\Xi^R(s;x,\xi)\right),
   $$
   then $h^{R\alpha\beta}_1(s),h^{R\alpha\beta}_2(s))$ satisfies the system   
   $$
\left\{
\begin{array}{rcl}
\displaystyle\frac{d\,}{ds}h^{R\alpha\beta}_1(s)&=&2\p a^R\left(X^R(\cdot)\right)\left(h^{R\alpha\beta}_1(s)\right)\Xi^R(\cdot)\\
&+&\displaystyle
    2a^R\left(X^R(\cdot)\right)\left(h^{R\alpha\beta}_2(s)\right)+Q^{R\alpha\beta}_1,\\
\displaystyle\frac{d\,}{ds}(h^{R\alpha\beta}_2(s))&=&-\p^2 a^R\left(X^R(\cdot)\right)\left(h^{R\alpha\beta}_1(s)\right)\Xi^R(\cdot)\Xi^R(\cdot)
\\
&-&\displaystyle
    2\p a^R\left(X^R(\cdot)\right)\Xi^R(\cdot)\left(h^{R\alpha\beta}_2(s)\right)+Q^{R\alpha\beta}_2,
\end{array}
\right.
$$
    where   
 $$
 \begin{array}{l}
\displaystyle Q^{R\alpha\beta}_j=
  Q^{R\alpha\beta}_j\left((\p^{\gamma}a)_{1\leq |\gamma|\leq |\alpha|+|\beta|+1};
 \left (\p_x^{\nu}\p_{\xi}^{\mu}X^R\right)_{|\nu|+|\mu|<|\alpha|+|\beta|};
  \left(\p_x^{\nu}\p_{\xi}^{\mu}\Xi^R\right)_{|\nu|+|\mu|<|\alpha|+|\beta|}\right),
  \end{array}
  $$
   $j=1,2$ is a polynomial in its variables. In others words, $\left(h^{R\alpha\beta}_1(s),h^{R\alpha\beta}_2(s)\right)$ satisfies a system
similar to that in 
   (\ref{eq4.2.3}) with external forces depending on the previous steps. We observe that each term in $Q^{R\alpha\beta}_j$, 
   $j=1,2$  has a factor of the form $\p^{\gamma}a\left(X^R(s;x,\xi)\right)$ with $|\gamma|\geq 1$. This guarantees that
$\p^{\gamma}a\left(X^R(s;x,\xi)\right)=0$ if $\left|X^R(s;x,\xi)\right|\geq 2R$ and if $|X^R(s;x,\xi)|\leq 2R$
   then $|s|\leq c(|x|+R)$, (see Theorem \ref{th4.1.5}). So when integrating the $Q^{R\alpha\beta}_j$'s one can substitute in 
(\ref{eq4.2.10}) the factor $(1+|s|)$ by $c(|x|+R)$ to get a bound independent of $s$ until the trajecory gets the free regime 
(i.e. $\left|X^R(s;x,\xi\right|>R$) which
provides   a linear in $s$
global bound.  Using a recursive argument we can obtain the following estimates which, although no sharp, suffices
  for our purpose here.
    
  \begin {theorem} \label{th4.2.1}  
For any $(x,\xi)\in \BbbR^n\times \BbbS^{n-1}$ and 
  for any $R\geq R_0$ with $R_0$ sufficiently large, the derivatives of the bicharacteristic flow satisfy : 
  given $\alpha,\beta\in\BbbZ^n$ with $|\alpha|+|\beta|\geq 1$ there exists $
  \mu_{\alpha,\beta}\in \BbbZ^+$ with $ \mu_{\alpha,\beta}\geq |\alpha|+|\beta|$ such that
  $$
\begin{array}{l}
  \displaystyle
   \left |\p^{\alpha}_x\p_{\xi}^{\beta}X^R(s;x,\xi)\right|
   \leq c\left(|s|+(|x|+R)^{\mu_{\alpha,\beta}}\right),\\
\displaystyle
   \left |\p^{\alpha}_x\p_{\xi}^{\beta}\Xi^R(s;x,\xi)\right|
   \leq c\left(|s|+(|x|+R)^{\mu_{\alpha,\beta}}\right),
\end{array}
$$
for any $s\in\BbbR$.
   
   \end{theorem}
 
 Combining the results in Theorem \ref{th4.2.1} with the identities in (\ref{eq4.2.4})-(\ref{eq4.2.5}) we  get the result for
$ (x,\xi)\in \BbbR^n\times \BbbR^n-\{0\}$.

\newpage

   \section{LINEAR ULTRAHYPERBOLIC EQUATIONS}\label{sec5}

In this section we shall deduce the linear estimates to be used in Section \ref{sec6} in the proof of our nonlinear results.
The first result in subsection \ref{subseq5.1} is the ultrahyperbolic version of Doi`s lemma \cite {8},
(see Lemma \ref{le 2.2.2} in Section \ref{sec2}). This will allow us to establish the smoothing effects for linear ultrahyperbolic equations with
variable second order coefficients.

 Subsection \ref{subseq5.2} is concerned with the $L^2$ (and $H^s$) well posedness of the associated linear problem  (\ref{eq1.4}) in the
Introduction. To establish this result we 
follow an indirect approach. As we did in Section \ref{sec4} -see (\ref{eq4.1.12}) and (\ref{eq4.1.13})-, we consider the truncation at infinity $ {\cal L}^R(x)$ of the operator ${\cal L}(x)$,
$$
 {\cal L}^R(x)=\theta(x/R){\cal L}(x)+\left(1-\theta(x/R)\right){\cal L}^0.
$$
For $R$ large enough we consider the bicharacteristic flow $(X^R(s;x,\xi),\Xi^R(s;x,\xi))$ (studied in Section \ref{sec4})  
associated to the operator ${\cal L}^R(x)$ and  the 
corresponding integrating factor $K^R$. To obtain the $L^2$ local
well posedness of the linear problem we combine several estimates for the operator 
${ \cal L}^R(x)$ and its associated \lq\lq errors", as function of $R$, with the
local smoothing effect obtained in subsection \ref{subseq5.1}.

\subsection{Linear Ultrahyperbolic Smoothing}\label{subseq5.1}

 We shall begin this subsection by proving the ultrahyperbolic version of Doi's lemma \cite {8}, (see Lemma \ref{le 2.2.2} in Section \ref{sec2}).

\begin {lemma} \label{le5.1.1}

 Assume that the bicharacteristic flow is non-trapped -see basic assumption in subsection \ref{subseq4.1} of Section
4, and that
$ \nabla a_{jk}(x)=o(|x|^{-1})$ as $|x|\to\infty $ for all
$j,k=1,..,n$. Suppose that
$$
\lambda\in L^1\left([0,\infty))\cap C([0,\infty\right)\;\;\mbox{is positive and
nonincreasing}.
$$
Then there exist $p\in S^0_{1,0}$ and $c>0$ such that
$$
\left(H_{h_2}p\right)(x,\xi)\geq \lambda(|x|)|\xi|-c\;\;\;\forall x,\xi\in\BbbR^n.
$$
\end{lemma}

\sl{Proof of Lemma \ref{le5.1.1}}  \rm

 Let $M>0$ be a constant to be chosen. Let $\psi\in C^{\infty}(\BbbR)$
with $\psi(t)=0$
for $t\leq M^2$, $\psi(t)=1$ for $t\geq (M+1)^2$ and $\psi'(t)\geq 0$ for $t\in
\BbbR$.
Let
$$
p_1(x,\xi)=\langle \xi\rangle^{-1}\psi\left(|x|^2\right)H_{h_2}\left(|x|^2\right)=-4\langle \xi\rangle^{-1}\psi\left(|x|^2\right)\left(x\cdot A_h(x)\xi\right).
$$

 By straightforward calculation -see Proposition \ref{pro4.1.1} of Section \ref{sec4}
$$
\begin{array}{rcl}
\displaystyle
\left(H_{h_2}p_1\right)(x,\xi)&=&\displaystyle
\langle \xi\rangle^{-1}\psi'(|x|^2)\left(H_{h_2}(|x|^2)\right)^2 \\
&+&\displaystyle\psi\left(|x|^2\right)\langle\xi\rangle^{-1}\Bigg[8\left|A(x)\xi\right|^2+8\sum_{j,k,l,m}x_l\p_{x_j}a_{lm}(x)a_{jk}(x)\xi_k\xi_m\\
&-&\displaystyle
4\sum_{j,k,l,m,p}x_ja_{jk}(x)\p_{x_k}a_{lm}(x)\xi_l\xi_m\Bigg].
\end{array}
$$

 With the assumptions on $A(x)$ and $c_1=4\nu^{-2}$ one can fix $M$
sufficiently large such that
$$
\left(Hp_1\right)(x,\xi)\geq c_1\psi\left(|x|^2\right)\langle
\xi\rangle^{-1}|\xi|^2,\;\;\;\;\forall  \,x,\xi\in\BbbR^n.
$$
 Now choose $\phi_1\in C^{\infty}_0(\BbbR^n)$ with $\phi_1(x)=1$ for
$|x|\leq M+1$. For $\xi\ne 0$, let
$$
p_2(x,\xi)=-\int_0^{\infty}\phi_1\left(X(s;x,\xi)\right)\left \langle \Xi(s;x,\xi)\right\rangle
ds.
$$
By theorem \ref{th4.1.5} in Section \ref{sec4}, for each $(x_0,\xi_0)\in\BbbR^n\times \BbbR^n-\{0\}$
there is a
neighborhood ${\cal U}$ of $(x_0,\xi_0)$ such that the integral defining
$p_2(x,\xi)$ is
over a fixed compact interval for all
$(x,\xi)\in {\cal U}$ so $p_2(\cdot,\cdot)$ is smooth. Furthermore, by
homogeneity of the bicharacteristic
flow (see (\ref{eq2.2.8}) in Section \ref{sec2}) and a change of variable
$$
\displaystyle p_2(x,\xi)=-|\xi|^{-1}\int_0^{\infty}\phi_1\left(X\left(s;x,\frac{\xi}{|\xi|}\right)\right)
\left\langle |\xi|\, \Xi\left(s;x,\frac{\xi}{|\xi|}\right)\right\rangle ds.
$$

 Choose $\phi_2\in C^{\infty}(\BbbR^n)$ with $\phi_2(\xi)=0$ for
$|\xi|\leq 1$ and 
$\phi_2(\xi)=1$ for $|\xi|\geq 2$. Let
$$
p_3(x,\xi)=\phi_1(x)\phi_2(\xi)p_2(x,\xi),\;\;\;\;x,\xi\in\BbbR^n.
$$

  Then $p_3\in S^0_{1,0}$ by the support properties of $\phi_1$ and
$\phi_2$, and
$$
\begin{array}{rcl}
\displaystyle (H_{h_2}{p_3})(x,\xi)
&=&\displaystyle
\left[2\sum_{jk}a_{jk}(x)\xi_k\p_{x_j}\phi_1(x)\right]\phi_2(\xi)p_2(x,\xi)\\
&+&\displaystyle
\phi_1(x)H_{h_2}\phi_2(\xi)p_2(x,\xi)+\left(\phi_1(x)\right)^2\phi_2(\xi)\langle
\xi \rangle.
\end{array}
$$

 Now let
$$
p_4(x,\xi)=c_2p_1(x,\xi)+p_3(x,\xi),
$$
with $c_2>0$ sufficiently large, $p_4$ then satisfies
$$
\left|\p_x^{\alpha}\p_{\xi}^{\beta}p_4(x,\xi)\right|
\leq c_{\alpha\beta}\langle
x\rangle\langle
\xi\rangle^{-|\beta|},\;\;\;\alpha,\beta\in\BbbN^n,
$$
and
$$
(H_{h_2}p_4)(x,\xi)\geq c_3|\xi|-c_4,\;\;\;\forall x,\xi\in\BbbR^n,
$$
where $c_{\alpha\beta}, c_3, c_4>0$ are constants.

 To complete the construction of $p(x,\xi)$ we observe that the proof of  Lemma 2.3 in \cite{8}
applies verbatim with
 $q(x,\xi)=p_4(x,\xi)$, since the proof does not depend on the ellipticity
of $A(x)$ assumed
in \cite{8}. With this, the proof of Lemma \ref{le5.1.1} is completed.
\bigskip

 Consider now  systems of the form
\begin{equation}\label{eq5.1.1}
\left\{
\begin{array}{l}
\p_t\vec w=i H\vec w + B\vec w +C\vec w+\vec f,\;\;\;\;(x,t)\in\BbbR^n\times(0,T),\\
\vec w(x,0)=\vec w_0(x).
\end{array}
\right.
\end{equation}

 Here $\vec w$ and $\vec f$ are $\BbbC^2$-value functions on $\BbbR^n\times(0,T)$,
$$
H=
\left(
 \begin{array}{cc}
 {\cal L}& 0\\
 0&-{\cal L}
 \end{array}
 \right),
$$
where $\displaystyle{\cal L}=\sum_{j,k}\p_{x_j}\left(a_{jk}(x)\p_{x_k}\right)$ and
$A(x)=(a_{jk}(x))$ satisfies the
assumptions in subsection \ref{subseq4.1} of Section \ref{sec4},
$$
B=
 \left(
 \begin{array}{cc}
 \Psi_{b_1}&\vec b_2(x)\cdot \nabla\\ 
 \bar{\vec b_2}(x)\cdot
\nabla&-\Psi_{\bar b_1}
  \end{array}
 \right)
$$
where $b_1\in S^1_{1,0}$ is odd in $\xi$ and $\vec b_2\in C^{\infty}_b(\BbbR^n:\BbbR^n)$, and
$$
C=
  \left(
 \begin{array}{cc}
\Psi_{c_{11}}&\Psi_{c_{12}}\\
 \Psi_{c_{21}}&\Psi_{c_{22}}
  \end{array}
 \right)
$$
where $c_{lm}\in S^1_{1,0}$ for $l,m=1,2$.

 Next we consider the ultrahyperbolic linear Schr\"odinger (scalar) equation
\begin{equation}\label{eq5.1.2}
 \left\{
 \begin{array}{l}
 \p_t u=i {\cal L} u + \Psi_{b_1}u+b_2(x)\cdot\nabla \bar u
+\Psi_{c_1}u +\Psi_{c_2}\bar u +f,\\
u(x,0)=u_0(x). 
\end{array}
\right.
 \end{equation}
 Taking $\vec f=(f\,,\bar f)^T$, $\vec w=(u \,,\bar u)^T$ and suitable $c$'s, equation (\ref{eq5.1.2})
is reduced to a system as (\ref{eq5.1.1}).

 \begin {theorem} \label{th5.1.2} Let $s\in\BbbR$. Then
there exists $N=N(n)\in \BbbN$ such that if 
$$
\left|\p_x^{\alpha}\p_{\xi}^{\beta}b_1(x,\xi)\right|
  \leq c_{\alpha \beta}\langle
x\rangle^{-N} \langle \xi \rangle^{1-|\beta|},\;\;\;\;x,\xi\in\BbbR^n,
 $$
and 
$$ 
\left|\p_x^{\alpha}\vec b_2(x)\right|\leq c_{\alpha \beta}\langle
x\rangle^{-N},\;\;\;x\in\BbbR^n, 
$$ 
then there exists $T>0$ so that thefollowing holds: Let $u_0\in H^s$. 

(A) If $ f\in L^1([0,T]:H^s)$, then
there is a unique solution $u \in C\left([0,T]:H^s\right)$ of the IVP (\ref{eq5.1.2})
satisfying 
$$
  \begin{array}{l} 
\displaystyle\sup_{0\leq t\leq T}\|u(t)\|_{H^s}+\left(\int_0^T\int_{\BbbR^n}\left|J^{s+1/2}u(x,t)\right|^2\langle x\rangle^{-N}dxdt\,\right)^{1/2}\\
\qquad\displaystyle\leq c\|u_0\|_{H^s}
+c\int_0^T \|f(t)\|_{H^s}dt. 
  \end{array}
$$ 
(B) If $f\in L^2([0,T]:H^s)$, then
there is a unique solution $u\in C([0,T]:H^s)$ of the IVP (\ref{eq5.1.2}) satisfying
$$ 
   \begin{array}{l} 
\displaystyle\sup_{0\leq t\leq T}\|u(t)\|_{H^s}^2+\int_0^T\int_{\BbbR^n}\left|J^{s+1/2}u(x,t)\right|^2\langle x\rangle^{-N}dxdt\\
\qquad \displaystyle\leq c\|u_0\|_{H^s}^2
+cT\int_0^T \|f(t)\|^2_{H^s}dt.
   \end{array}
$$ 

(C) If 
$$ J^{s-1/2}f\in
L^2\left(\BbbR^n\times[0,T]:\langle x\rangle^{N}dxdt\right)
 $$ 
then there is a
unique solution $u\in C([0,T]:H^s)$ of the IVP (\ref{eq5.1.2}) satisfying
 $$\begin{array}{l}
\displaystyle\sup_{0\leq t\leq T}\|u(t)\|_{H^s}^2+\int_0^T\int_{\BbbR^n}\left|J^{s+1/2}u(x,t)\right|^2\langle x\rangle^{-N}dxdt\\
\qquad\displaystyle\leq c\|u_0\|_{H^s}^2
+cT\int_0^T\int_{\BbbR^n}\left|J^{s-1/2}f(x,t)\right|^2\langle x\rangle^N\,dxdt.
  \end{array}
$$ 
Here $c$ depends on $s,$ $\nu,$ $A(x),$ $b_1,$ $\vec b_2,$ and $c_{lm}.$

\end{theorem}
 The goal is to prove Theorem \ref{th5.1.2}.
 The \it  a priori \rm  estimates needed to prove Theorem \ref{th5.1.2} are
essentially reduced to the case $s=0$ by the following commutator result.

\begin {lemma} \label{le5.1.3}

$$
J^s(iH+B+C)=(iH+\tilde  B +\tilde C)J^s,
$$
where
\begin{equation}\label{orders}
\tilde B=-i s \sum_{j,k,l}\p_{x_j}a_{kl}(x)\p^3_{x_jx_kx_l}J^{-2}
\left(
\begin{array}{cc}
1&0\\
0&-1
\end{array}
\right) +B,
\end{equation}
$$
\tilde C=
\left(
\begin{array}{cc}
 \Psi_{\tilde c_{11}}&\Psi_{\tilde c_{11}}\\
\Psi_{\tilde c_{11}}&\Psi_{\tilde c_{11}}
\end{array}
\right)
$$
and $\tilde c_{lm}\in S^0_{1,0}$ for $l,m=1,2$.

\end{lemma}

A similar result applies to the scalar equation (\ref{eq5.1.2}).

\sl {Proof of Lemma \ref{le5.1.3}}\rm

The lemma follows from the classical $S^0_{1,0}$ pseudo-differential
calculus of subsection \ref{subseq2.1} in Section \ref{sec2}.

\bigskip

 The next step is to obtain an \it a priori \rm estimate for Theorem \ref{th5.1.2}
in the case $s=0$.
The last step will be the reduction to the case $s=0$ (via Lemma
\ref{le5.1.3}). From this point, $c$ will denote a 
constant (not necessarily the 
same at each appearance) depending on $s, \nu, A(x), b_1, \vec b_2,
c_{lm}$ and the hypothesis of Theorem \ref{th5.1.2} is assumed to be
satisfied. In order to simplify notation,
$(H^s)^2=H^s\times H^s$ will be denoted merely by $H^s$.

\begin{lemma} \label{le5.1.4}

 Let $T>0$. For all $\vec w\in C\left([0,T]:H^2)\cap C^1([0,T]:L^2\right)$
$$
\displaystyle\int_0^T\int_{\BbbR^n}\left|J^{1/2}\vec w(x,t)\right|^2\langle x\rangle^{-N}dxdt\leq
c(1+T)\sup_{0\leq t\leq T}\|\vec w(t)\|_{L^2}^2+\int_0^T\|\vec f(t)\|^2_{L^2}dt
$$
and
$$\begin{array}{l}
\displaystyle \int_0^T\int_{\BbbR^n}\left|J^{1/2}\vec w(x,t)\right|^2\langle x\rangle^{-N}dxdt\\
\displaystyle\leq c\left((1+T)\sup_{0\leq t\leq T}\|\vec w(t)\|_{L^2}^2+\int_0^T\int_{\BbbR^n}\left|J^{-1/2}\vec f\right|^2 \langle x\rangle^{N}dxdt\right),
\end{array}
$$
where
$$
\vec f(x,t)= \p_t\vec w-(iH+\tilde B+\tilde C)\vec w
$$
and $\tilde B,\, \tilde C$ were given in Lemma \ref{le5.1.3}.

\end{lemma}

\sl{Proof of Lemma \ref{le5.1.4}}\rm

  By Lemma \ref{le5.1.1} there
exists a real-valued $p\in S^0_{1,0}$ and $c>0$ such that 
$$
(H_{h_2}p)(x,\xi)\geq c' \langle x\rangle^{-N}|\xi|-c,\;\;\;\;x,\xi\in\BbbR^n, 
$$ 
where $c'=c'(s)$ is to be determined. Let 
$$ k(x,\xi)=
\left(
\begin{array}{cc}
\exp(p(x,\xi))&0\\
0&-\exp(p(x,\xi))
\end{array}
\right)
$$ and $K=\Psi_k$. Then $K$
is a diagonal $2\times 2$ matrix of $S^0_{1,0}$ pseudo-differential
operators.
 One can now calculate as follows 
$$ 
\begin{array}{rcl}
\p_t\langle K\vec w,\vec
w\rangle_{L^2\times L^2}&=& \displaystyle
\left\langle K\p_t\vec w,\vec w\right\rangle_{L^2\times
L^2} +\left \langle K\vec w,\p_t\vec w\right\rangle_{L^2\times L^2}\\ 
&=&\displaystyle
\left\langle
\left(i[KH-HK]+K\tilde B +\tilde B^*K\right)\vec w,\vec w\right\rangle_{L^2\times L^2}\\
&+&\displaystyle
\left\langle \left[K\tilde C+\tilde C^*K\right]\vec w,\vec w\right\rangle_{L^2\times L^2} + \left(\left\langle K\vec
f,\vec w\right\rangle_{L^2\times L^2}+ \left\langle K\vec w,\vec f\right\rangle_{L^2\times
L^2}\right)\\
 &=& I + II + III. 
\end{array}
$$ 

 Disregarding symbols of order $0$, the first order symbol of
$i[KH-HK]+K\tilde B +\tilde B^*K$ is
$$
\begin{array}{l}
-e^p 
\left(
\begin{array}{cc}
\nabla_{\xi}h_2\cdot\nabla_xp-\nabla_xh_2\cdot\nabla_{\xi}p&0\\
0&\nabla_{\xi}h_2\cdot\nabla_xp-\nabla_xh_2\cdot\nabla_{\xi}p
\end{array}
\right)
\\
-2se^p
\left(
\begin{array}{cc}
\displaystyle \sum_{jkl}\p_{x_j}a_{kl}(x)\xi_j\xi_k\xi_l\langle
\xi\rangle^{-2}&0\\0&
\displaystyle\sum_{jkl}\p_{x_j}a_{kl}(x)\xi_j\xi_k\xi_l\langle
\xi\rangle^{-2}
\end{array}
\right)
\\
\\
+e^p
\left(
\begin{array}{cc}
 2\,\mbox{Re}\,b_1&2i\,\vec b_2(x)\cdot\xi\\
-2i\bar{\vec b}_2(x)\cdot\xi& 2\,\mbox{Re}\,b_1
\end{array}
\right).
\end{array}
$$

 Choosing $c'(s)$ large enough and using the matrix version of the sharp
G\aa rding inequality in \cite {12}, it follows that
$$
\begin{array}{rcl}
|I|&\leq&\displaystyle -c \left \langle 
 \left(
 \begin{array}{cc}
  \lambda J^1&0
  \\0&\lambda
J^1
  \end{array}
  \right)
  \vec w,\vec w\right \rangle +c\|\vec w\|^2_{L^2}\\
\\
&\leq& \displaystyle-c\int_{\BbbR^n}\left|J^{1/2}\vec w(x,t)\right|^2 \langle x\rangle^{-N}dx +c\|
\vec w\|^2_{L^2}.
\end{array}
$$

 Using that both $K$ and $C$ are of order $0$
$$
|II|\leq c\|\vec w\|^2_{L^2}.
$$
 The  estimate of $III$ is split into two cases according to the desired
norm of $\vec f$. This is identical to the elliptic case given in Section \ref{sec2}.

Combining estimates, integrating in $t$ and using the $L^2$-boundedness of
$K$, Lemma \ref{le5.1.4} follows.

\bigskip

\begin{remark}\label{re5.1.5}  We will use Lemma \ref{le5.1.4} for $\vec
\omega=(J^su,\overline{J^s u})$, where $u$ solves the scalar equation
(\ref{eq5.1.2}).
\end{remark}

\subsection{Linear Ultrahyperbolic $L^2$-Well posedness}\label{subseq5.2}
  
 Our goal in this subsection is to established the following result,

 \begin{lemma} \label{le5.2.1} Let $\,T>0\,$. For all $\,u\in
 C\left([0,T];H^2\right)\cap C^1\left([0,T];H^2\right)\,$ the
 following two estimates hold.

 \begin{itemize}
 \item[(A)]
 $\displaystyle\sup_{0\le t\le T}\|u\|_{L^2}\le
 c\|u(0)\|_{L^2}+c
 T\sup_{0\le t\le T}\|u\|_{L^2}+c\int_0^T\| f\|_{L^2}dt$,
\item[(B)]
 $\displaystyle\sup_{0\le t\le T}\|u\|^2_{L^2}\le
 c\|u(0)\|^2_{L^2}+c
 T\sup_{0\le t\le T}\|u\|^2_{L^2}+c\int_0^T\| f\|^2_{L^2}dt.$

 \end{itemize}
 Here 
 $ f=\p_t u-\left(i{\cal L}u+\Psi_{b_1}u+\vec b_2\cdot\nabla\bar u+\Psi_{c_1}u+\Psi_{c_2}\bar u\right)\,$ analogously to Lemma \ref{le5.1.4}.
\end{lemma}

 The proof of Lemma \ref{le5.2.1} involves several steps. The first
 one is to cut ${\cal L}- {\cal L}_0$ at
 infinity. Therefore for $R$ large enough
 and to be fixed later on we define
$$
\displaystyle a_{jk}^R(x)=\theta \left(\frac{x}{R}\right)a_{jk}(x)+\left(1-\theta\left(\frac{x}{R}\right)\right)a_{jk}^0
=a_{jk}^0+\theta\left(\frac{x}{R}\right)\left(a_{jk}(x)-a_{jk}^0\right),
$$
with $\theta$ a smooth cut off function such that
 $\theta(x)=1$ if $|x|<1$ and $\theta(x)=0$ if $|x|>2$.
Define
$${\cal L}^R(x)=-\sum_{j,k}\frac{\partial}
{\partial x_j}\left(a_{jk}^R(x)\frac{\partial}
{\partial x_k}\right),
$$
$h_2^R$ its corresponding hamiltonian
 and 
$${\cal E}^R(x)={\cal L}(x)-{\cal L}^R(x).
$$
Therefore
$$
\displaystyle{\cal E}^R(x)=-\sum_{j,k}\frac{\partial}
{\partial x_j}\left(e_{jk}(x)\frac{\partial}
{\partial x_k}\right)
$$
with 
$$
\displaystyle e_{jk}(x)=\left(1-\theta\left(\frac{x}{R}\right)\right)\left(a_{jk}(x)-a_{jk}^0\right).
$$
  The main point is to apply a suitable 
 pseudo-differential operator $\,K^R\,$  to
 the corresponding system
 in order to cancel the first order terms.  We begin with the
 definition and some properties
 of $\,K^R\,$ and its symbol. We recall that $\chi\in C^{\infty}(\BbbR)$, $\chi(t)=0$ for $|t|\leq 1$ and
$\chi(t)=1$ for $|t|\geq 2$. 

 \begin{definition} \label{def5.2.2} Let $R=2^{j_0}$ for $j_0\in \BbbN.$ Recall
the definition of $\tilde B$ in Lemma \ref{le5.1.3}. We define
 \begin{itemize}
 \item [(i)]
 $\displaystyle b^R(x,\xi)=s\sum\limits_{j,k,l}\partial_{x_j}a_{jk}^R(x)
 \xi_j\xi_k\xi_l\langle\xi\rangle^{-2}-\mbox{Re}\,b_1(x,\xi),$
 \item [(ii)]
 $ p^R(x,\xi)=\chi\left(\displaystyle\frac 12|\xi|\right)\displaystyle\int_{-\infty}^0
 b^R\left(X^R(\sigma;x,\xi),\Xi(\sigma;x,\xi)\right)d\sigma,$
 \item[(iii)]
 $\displaystyle
  p_e^R(x,\xi)=\frac
 12\left(p^R(x,\xi)+p^R(x,-\xi)\right),$
 \item[(iv)]
 $k^R(x,\xi)=\exp
 \left(p_e^R(x,\xi)\right),$
 \item[(v)]
 $\tilde k^R(x,\xi)=\exp
 \left(-p_e^R(x,\xi)\right),$
 \item[(vi)]
 $B^R=\Psi_{b^R},\, P^R=\Psi_{p^R},\,
 P^R_e=\Psi_{p^R_e},\, K^R=\Psi_{k^R},\,\tilde
 K^R=\Psi_{\tilde k^R},$
 \item[(vii)] $\displaystyle b^R=b_0^R+ \sum_{j=j_0}^\infty\,b_j^R$, with
 $\displaystyle b_j^R(x)=\left(\theta\left(\frac{x}{2^{j+1}}\right)-\theta\left(\frac{x}{2^{j}}\right)\right)b^R(x)$, and
 analogously $p_j^R$, 
 $p_{e\,j}^R$, $k^R_j$, and  $\tilde k^R_j.$
 \end{itemize}
\end{definition}
  
  Some comments about the above definition are in order.
  
  We recall that applying the operator $J^s$ to the equation in (\ref{eq5.1.2}) the symbol of the first order term for $J^s u$ is given by (see (\ref{orders}))
  $$ s\displaystyle \sum_{j,k,l} \partial_{x_j} a_{jk}(x)\xi_j\xi_k\xi_l\langle\xi\rangle^{-2}
  -\mbox{Re}\,b_1(x,\xi).$$
So $b^R (x,\xi)$ in (i) is the approximation (due to our truncation of the second order term) of this symbol.
  
 The symbol of the integrating factor needed to ``cancel'' $\Psi_{b^R}$ is given in (iv). The reason to use the even function $p^R_e(x,\xi)$ in (iii) instead of that in (ii) is to preserve the symmetry needed in the integration by parts use to handle the first order term in $\bar u$. 
 
 Finally, $\tilde K^R(x,\xi)$ is the symbol of a cuasi-inverse of $\Psi_K^R$.
 \begin{lemma} \label{le5.2.3}
  
 \begin{itemize}
 \item[(i)]
 $p_e^R,\,k^R,\,\tilde k^R\,$ are all even in $\,\xi$, and real.
 \item[(ii)]
 $\left(\nabla_{\xi} h_2^R\cdot \nabla_x p^R- \nabla_{x}
 h_2^R\cdot \nabla_{\xi} p^R\right)(x,\xi)=
 b^R(x,\xi)+r_1(x,\xi)\,$
 where
 $\,r_1\in S_{1,0}^{-\infty}.$
 \item[(iii)]
 $\left(\nabla_{\xi} h_2^R\cdot \nabla_x
 p^R_e- \nabla_{x} h_2^R\cdot \nabla_{\xi}p^R_e\right)(x,\xi)=
 b^R(x,\xi)+r_2(x,\xi)\,$
 where
 $\,r_2\in S_{1,0}^{-\infty}.$
 \item[(iv)]
 $\left(\nabla_{\xi}
 h_2^R\cdot \nabla_x k^R- \nabla_{x} h_2^R\cdot \nabla_{\xi} k^R\right)(x,\xi)=
 k^R(x,\xi)
 b^R(x,\xi)+r_3(x,\xi)\,$
 where
 $\,r_3\in S_{1,0}^{-\infty}.$
 \item[(v)]
 $\left(\nabla_{\xi} h_2^R\cdot \nabla_x \tilde k^R- \nabla_{x}
 h_2^R\cdot \nabla_{\xi} \tilde k^R\right)(x,\xi)=
 -\tilde k^R(x,\xi)
 b^R(x,\xi)+r_4(x,\xi)\,$
 where
 $\,r_4\in S_{1,0}^{-\infty}.$
 \item[(vi)]
 Let $\psi\in C^{\infty}(\BbbR)$, $\psi(t)=1$
 if $t\ge \frac 12$ and $\psi(t)=-1$
  if $t\le -\frac 12,\,-1\le
 \psi\le 1$.
 Define
 \begin{equation}\label{eq5.2.1}
a_j^R(z;x,\xi)=
\displaystyle\chi\left(\frac 1{10\,2^j}|x|\right)
 p_{e\,j}^R\left(z+2^j\psi\left(\frac{x\cdot
A_h\xi}{|x||\xi|}\right)
 \frac{A_h\xi}{|\xi|},\xi\right)\,\theta\left(\frac{z}{2^{j+1}}\right),
\end{equation}
 $j=0,1,...$. Then 
 $\,a_j^R\in C^{\infty}_0(B_{2^{j+1}}(0);S^0_{1,0})\,$ and
 \begin{equation}\label{eq5.2.2}
  \begin{array}{l}
 p_{e\,j}^R(x,\xi)=\\
\displaystyle p_{e\,j}^R(x,\xi)\left(
 1-\chi\left(\frac
 1{10\,2^j}|x|\right)\right)
 +a_j^R\left(P(x,A_h\xi);x,\xi\right).
\end{array}
\end{equation}
 \item[(vii)] Set
 $\displaystyle q(x,\xi)=\sum_j\,p_{e\,j}^R\left(1-\chi\left(\frac
 1{10\,2^j}|x|\right)\right)$,
 and 
 $$
\displaystyle a^R(z;x,\xi)=\exp(q)\left(\exp(\sum_j\,a_j^R(z;x,\xi)) -1\right).
 $$ 
 Then
 $q\,\in S^0_{1,0}$, $q(\,\cdot\,,\xi)\in{\cal S}(\BbbR^n)$ uniformly in $\xi$,  $a^R\,\in {\cal S}(\BbbR^n;S^0_{1,0})\,$, and 
 \begin{equation}\label{eq5.2.3}
 k^R(x,\xi)=a^R\left(P(x,A_h\,\xi);x,\xi\right)+\exp(q(x,\xi)).
\end{equation}
\item[(viii)] The seminorms of the remainders $r_k$, $k=1,2,3,4$, and of $a^R_j$ and $q$ grow as $R^{N_0}$ with $N_0\in \mathbb Z^+$ depending just on the dimension.
 \end{itemize}
\end{lemma}
  
  Parts (i)--(iii) are preliminary results needed in the proof of (iv). The crucial point (iv) shows at the level of the symbols that the commutator of $K^R$ and ${\cal L}^R(x)$ ``cancels'' $K^RB^R$. In (vi) and (vii) we prove that the symbol $k^R(x,\xi)$ is in the class introduced in Section \ref{sec3}.

 \sl{Proof of Lemma \ref{le5.2.3}}\rm

 For simplicity of the exposition we shall drop the
 index $R$.

 The proof of (i) is clear.

 As for (ii), it is sufficient to work
 with each $b_j$ and the corresponding $p_j$. We get
 $$
\begin{array}{l}
 \left(\nabla_{\xi}
\displaystyle  h_2^R\cdot \nabla_x p_j- \nabla_{x} h_2^R\cdot \nabla_{\xi }p_j\right)(x,\xi)\\
 \qquad\displaystyle
  =-\int_{-\infty}^0 b_j\left(X(\sigma;x,\xi),\Xi(\sigma;x,\xi)\right)d\sigma\cdot
\left (\nabla_x h_2^R\cdot
 \nabla_{\xi}\left[\chi\left(\frac 12|\xi|\right)\right]\right)\\
\qquad\displaystyle
  +\chi\left(\frac
 12|\xi|\right)
\left(\nabla_{\xi} h_2^R\cdot \nabla_x - \nabla_{x} h_2^R\cdot \nabla_x\right)
 \int_{-\infty}^0 b_j\left(X(\sigma;x,\xi),\Xi(\sigma;x,\xi)\right)d\sigma.
 \end{array}
$$
 For the first factor of the first term
 we use Theorem \ref{th4.1.5} in Section \ref{sec4} to see that is bounded by $c_N2^{-Nj}$, 
and therefore independently of $R$. Then by
Theorem \ref{th4.2.1} in Section \ref{sec4} we
 prove
 that the derivatives $\p_x^{\alpha}\p_{\xi}^{\beta}$  inside the integral
are bounded by 
 $c_{\alpha \beta}\left(|\sigma|+(|x|+R)^{\mu_{\alpha \beta}}\right)$ for $\alpha, \beta\in \BbbN^n$. 
If  $|\sigma|>s_0$ where $s_0$ is given in Theorem \ref{th4.1.5} of Section \ref{sec4},
then $|\sigma|<c|X(\sigma;x,\xi)|<c2^j$. Otherwise
$|\sigma|<c|x|<cR$, 
where the
last inequality follows from the support properties of
$\nabla_xh_2^R$. Therefore
 from the decay of $b_j$ and $\nabla A$ we conclude that
 belongs to $S^{-\infty}_{1,0}$ and the corresponding seminorms
grow like powers of $R$.
  By Lemma \ref{le 2.2.1} in Section \ref{sec2}, the second term equals
 $$
\displaystyle \chi\left(\frac
 12|\xi|\right)b_j(x,\xi)=b_j(x,\xi)+
 \left(\chi\left(\frac 12|\xi|\right)-1\right)
 b_j(x,\xi).
$$
 Notice that the last term above is compactly supported in 
 $\,\xi\,$. 

 (iii)
 $\,h_2^R\,$ is homogeneous of degree 2 in $\,\xi\,$, so
 $$
\begin{array}{l}
 \displaystyle
  \left(\nabla_{\xi} h_2^R\cdot \nabla_x  - \nabla_{x} h_2^R\cdot \nabla_{\xi}
\right)
 [p(x,-\xi)]\\
 \qquad\displaystyle
  =\left(\nabla_{\xi} h_2^R\right)(x,\xi)\cdot \left(\nabla_x p\right)(x,-\xi) + \left(\nabla_{x}
 h_2^R\right)(x,\xi)\cdot (\nabla_{\xi}
 p)(x,-\xi)\\ 
 \qquad\displaystyle
  =\left(-\nabla_{\xi} h_2^R\cdot \nabla_x
 p+\nabla_{x} h_2^R\cdot \nabla_{\xi} p\right)(x,-\xi).
 \end{array}
$$
 By (ii), this
 equals 
 $$
 -b(x,-\xi)-r_1(x,-\xi)
 =b(x,\xi)-r_1(x,-\xi),
 $$
 since $\,b\,$ is odd in $\,\xi\,$. This proves (iii).

 (iv) By the chain rule -see also the proof of (ii),
 $$
\begin{array}{l}
\displaystyle
   \left[\nabla_{\xi} h_2^R\cdot \nabla_x k - \nabla_{x} h_2^R\cdot
 \nabla_{\xi} k\right](x,\xi)\\
 \qquad\displaystyle
   =k(x,\xi)\left(\nabla_{\xi} h_2^R\cdot \nabla_x p_e
 - \nabla_{x} h_2^R\cdot \nabla_{\xi}
 p_e\right) (x,\xi)\\
 \qquad\displaystyle
   =k(x,\xi)\left(b(x,\xi)+r_2(x,\xi)\right),
 \end{array}
$$
 where $\,r_2(x,\xi)\,$ is compactly
 supported in $\,\xi$.

(v) is similar to (iv).

 (vi) From the definition of $a_j$ we have
 $\supp a_j(\cdot;x,\xi)\subset B_{2^{j+1}}(0)$.

By homogeneity,
 $$
\displaystyle
    b_j\left(X(\sigma;x,\xi),\Xi(\sigma;x,\xi)\right)=
 b_j\left(X\left(|\xi|\sigma;x,\frac{\xi}{|\xi|}\right),
 |\xi|\Xi\left(|\xi|\sigma;x,\frac{\xi}{|\xi|}\right)\right),
$$
 so it follows by a
 change of variable that
 \begin{equation}\label{eq5.2.4}
  \begin{array}{rcl}
   p_{e\,j}(x,\xi)&=&
\displaystyle\frac 12 \chi\left(\frac 12|\xi|\right)
\frac 1{|\xi|}\left(\int_{-\infty}^0
 b_j\left(X\left(\sigma;x,\frac{\xi}{|\xi|}\right),
|\xi|\Xi\left(\sigma;x,\frac{\xi}{|\xi|}\right)\right) d\sigma-\right.\\
&-& \displaystyle\left.\int_0^{\infty}
 b_j\left(X\left(\sigma;x,\frac{\xi}{|\xi|}\right),
 |\xi|\Xi\left(\sigma;x,\frac{\xi}{|\xi|}\right)\right) d\sigma\right).
 \end{array}
\end{equation}
Thus combining (\ref{eq5.2.4}) and Theorem \ref{th4.2.1} in Section \ref{sec4} and proceeding as in part (ii)
we obtain that 
$a_j(z;\cdot,\cdot)\in S^0_{1,0}$ uniformly in $z$, and so are in the
derivatives with respect to $z$.

  Let us prove (\ref{eq5.2.2}). We have to see that
 \begin{equation}\label{eq5.2.5}
\begin{array}{c}
 \displaystyle
  \chi\left(\frac 1{10\, 2^j}|x|\right)p_{e\,j}(x,\xi)\\
\displaystyle=\chi\left(\frac 1{10\, 2^j}|x|\right)
 p_{e\,j}\left(P(x,A_h\xi)+
 2^j\psi\left(\frac{x\cdot A_h\xi}{|x||\xi|}\right)
 \frac{A_h\xi}{|\xi|},\xi
 \right)\theta\left(\frac{P(x,A_h\xi)}{2^{j+1}}\right).
 \end{array}
\end{equation}
 Consider different cases.

 $1^{\underline {\mbox o}}\,$ $\,|x|\le 10\,2^j\,$: Here
 $\displaystyle\,\chi\left(|\frac {|x|}{10\,2^j}|\right)=0\,$ so
 RHS=LHS in (\ref{eq5.2.5}).

 $2^{\underline {\mbox o}}\,$ $\,|\xi|\le 2\,$: $RHS=LHS=0.$

 $3^{\underline
 {\mbox o}}\,$ $\,|P(x, A_h\xi)|\ge 2^j\,$ and $\,|x|\geq
10\,2^j\,$: First note
that
 $$
   \displaystyle
P\left(P(x,A_h\xi)+
 2^j\psi\left(\frac{x\cdot A_h\xi}{|x||\xi|}\right)
 \frac{A_h\xi}{|\xi|},A_h\xi
 \right)=
 P(x, A_h\xi).
$$
 Next, if $\, |P(y, A_h\xi)|\ge 2^jR\,$ for some $\,y\in
 \BbbR^n\,$, then 
 $\,b_j(X(\sigma;y,\xi),\Xi (\sigma;y,\xi))=0\,$ for all
 $\,\sigma\in \BbbR\,$ because $\,b_j\,$ has
 $\,x$-support in $\,B_{2^j}(0)\,,$ and the bicharacteristics are lines. It follows that
 RHS=LHS=0 in (\ref{eq5.2.5}).

 $4^{\underline {\mbox o}}\,$ $\,|x|\ge 10\,2^j,\,|\xi|\ge 2\,,\,|P(x,
 A_h\xi)|\le 2^j\,$:
 Here $\displaystyle\,\chi(|\xi|/2)=1\,$ and $\theta(\frac{P(x,
 A_h\xi)}{2^{j+1}})=1.$

 $$
   \displaystyle
|x\cdot A_h\xi|=|x||\xi|
 \sqrt{1-\frac{|P(x,
 A_h\xi)|^2}{|x|^2}}\ge
 \frac{99}{100}|x||\xi|. 
$$
 Now split in two subcases according to the sign of
 $\,x\cdot A_h\xi$.

 $4^{\underline {\mbox o}}\,(a)\,$ 
 $\displaystyle x\cdot A_h\xi\geq
 \frac{99}{100}|x||\xi|$: 
 Here $\,\psi\left(\frac{x\cdot
 A_h\xi}{|x||\xi|}\right)=1$.
 Suppose $\displaystyle\,\sigma\le\frac
 12\left(x\cdot\frac{A_h\xi}{|\xi|}-2^{j+1}\right)$. 
Then
 $$
   \displaystyle
\left(x-2\sigma
 \frac{A_h\xi}{|\xi|}\right)\cdot
 \frac{A_h\xi}{|\xi|}=x\cdot
 \frac{A_h\xi}{|\xi|}-2\sigma\ge 2^{j+1}.
$$
 Hence
 $\displaystyle\,X\left(\sigma;x,\frac{\xi}{|\xi|}\right)=x-2\sigma
\frac{A_h\xi}{|\xi|}\,$
 and
 $$ 
b_j\left(X\left(\sigma;x,\frac{\xi}{|\xi|}\right),
 |\xi|\Xi\left(\sigma;x,\frac{\xi}{|\xi|}\right)\right)=0. 
$$
 
 By (\ref{eq5.2.4}) and since
$\displaystyle\,\frac 12\left(\frac{x\cdot A_h\xi}{|\xi|}-2^{j+1}\right)>0\,$,
 $$
\begin{array}{l}
 p_{e\,j}(x,\xi)=\\
\displaystyle-\frac 12\chi
 \left(\frac 12|\xi|\right)\frac
 1{|\xi|}\int_{1/2(x\cdot A_h\xi/|\xi|-2^{j+1})}^{\infty}
 b_j\left(X\left(\sigma;x,\frac{\xi}{|\xi|}\right),
 |\xi|\Xi\left(\sigma;x,\frac{\xi}{|\xi|}\right)\right)d\sigma.
 \end{array}
$$
 Doing the
 change of variable 
 $\displaystyle\,\tau=\sigma-\frac 12\left(x\cdot
 \frac{A_h\xi}{|\xi|}-2^{j+1}\right)$,
 $$
\begin{array}{l}
 p_{e\,j}(x,\xi)=\\
\displaystyle-\frac 12\chi
 \left(\frac 12|\xi|\right)\frac 1{|\xi|}\int_{0}^{\infty}
 b_j\left(X\left(\tau+\frac 12\left(x\cdot
 \frac{A_h\xi}{|\xi|}-2^{j+1}\right);x,\frac{\xi}{|\xi|}\right),
 \right.\\
 \displaystyle,\left. |\xi|\Xi\left(\tau+\frac 12\left(x\cdot
\frac{A_h\xi}{|\xi|}-2^{j+1}\right);x,\frac{\xi}{|\xi|}\right)\right)d\tau.
 \end{array}
$$
 Now,
 $$
\begin{array}{c}
\displaystyle X\left(\tau+\frac 12\left(x\cdot
 \frac{A_h\xi}{|\xi|}-2^{j+1}\right);x,\frac{\xi}{|\xi|}\right)\\ 
\displaystyle=X\left(\tau;
 X\left(\frac 12\left(x\cdot
 \frac{A_h\xi}{|\xi|}-2^{j+1}\right);x,\frac{\xi}{|\xi|}\right),
 \Xi\left(\frac
 12\left(x\cdot
 \frac{A_h\xi}{|\xi|}-2^{j+1}\right);x,\frac{\xi}{|\xi|}\right)\right)\\
\displaystyle=X\left(\tau;P(x,A_h\xi)+2^{j+1}\frac{A_h\xi}{|\xi|},\frac{\xi}{|\xi|}\right),
 \end{array}
$$
 and similarly,
 $$\displaystyle \Xi\left(\tau+\frac 12\left(x\cdot
 \frac{A_h\xi}{|\xi|}-2^{j+1}\right);x,\frac{\xi}{|\xi|}\right)=
 \Xi
\left(\tau;P(x,A_h\xi)+2^{j+1}\frac{A_h\xi}{|\xi|},\frac{\xi}{|\xi|}\right) $$
 If $\,\tau\leq 0\,$ then
 $$
   \displaystyle
b_j\left(X\left(\tau;P(x,A_h\xi)+2^{j+1}\frac{A_h\xi}{|\xi|},\frac{\xi}{|\xi|}\right)
,|\xi|\Xi\left(\tau;P(x,A_h\xi)+2^{j+1}\frac{A_h\xi}{|\xi|},\frac{\xi}{|\xi|}\right)\right)=0.
 $$

 Hence by (\ref{eq5.2.4})
 $$
   \displaystyle
p_{e\,j}(x,\xi)=
 p_{e\,j}\left(P(x,A_h\xi)+
 2^{j+1}\psi\left(\frac{x\cdot
 A_h\xi}{|x||\xi|}\right)\frac{A_h\xi}{|\xi|},\xi\right).
$$
 Therefore, (\ref{eq5.2.5})
 holds.
This finishes case $4^{\underline {\mbox{o}}}\,(a)$.

 $4^{\underline {\mbox o}}\,(b)\,$ 
$x\cdot
 A_h\xi\le
 -\frac{99}{100}|x||\xi|$: 
 Here $\displaystyle\psi\left(\frac{x\cdot
 A_h\xi}{|x||\xi|}\right)=-1$.

 Suppose $\displaystyle\,\sigma\ge\frac
 12\left(x\cdot\frac{A_h\xi}{|\xi|}+2^{j+1}\right)\,$. Then
 $$\displaystyle
   \left(x-2\sigma
 \frac{A_h\xi}{|\xi|}\right)\cdot
 \frac{A_h\xi}{|\xi|}=x\cdot
 \frac{A_h\xi}{|\xi|}-2\sigma\le -2^{j+1}.$$
 Hence
 $\,X\left(\sigma;x,\frac{\xi}{|\xi|}\right)=x-2\sigma
\displaystyle\frac{A_h\xi}{|\xi|}\,$
 and
 $\displaystyle b_j\left(X\left(\sigma;x,\frac{\xi}{|\xi|}\right),
 |\xi|\Xi\left(\sigma;x,\frac{\xi}{|\xi|}\right)\right)=0. $
 Therefore since 
 $\displaystyle\,\frac 12\left(2^{j+1}+\frac{x\cdot A_h\xi}{|\xi|}\right)<0\,$,
$$
\begin{array}{l}
 p_{e\,j}(x,\xi)\\
\displaystyle=-\frac 12\chi
 \left(\frac 12|\xi|\right)\frac
 1{|\xi|}\int_{-\infty}^{1/2(2^{j+1}+x\cdot A_h\xi/|\xi|)}
 b_j\left(X\left(\sigma;x,\frac{\xi}{|\xi|}\right),
 |\xi|\Xi\left(\sigma;x,\frac{\xi}{|\xi|}\right)\right)d\sigma.
 \end{array}
$$
 Let
 $\displaystyle\,\tau=\sigma-\frac 12\left(2^{j+1}+x\cdot
\frac{A_h\xi}{|\xi|}\right)\,$. Then
 $$
\begin{array}{l}
 p_{e\,j}(x,\xi)=\\
\displaystyle-\frac 12\chi
 \left(\frac 12|\xi|\right)\frac
 1{|\xi|}\int_{-\infty}^{0}
 b_j\left(X\left(\tau+\frac 12\left(2^{j+1}+x\cdot
 \frac{A_h\xi}{|\xi|}\right);x,\frac{\xi}{|\xi|}\right),
 \right.\\
\displaystyle,\left.
 |\xi|\Xi\left(\tau+\frac 12\left(2^{j+1}+x\cdot
 \frac{A_h\xi}{|\xi|}\right);x,\frac{\xi}{|\xi|}\right)\right)d\tau.
 \end{array}
$$
 But
 $$
   \displaystyle
 X\left(\tau+\frac 12\left(2^{j+1}+x\cdot
 \frac{A_h\xi}{|\xi|}\right);x,\frac{\xi}{|\xi|}\right)
=X\left(\tau;P(x,A_h\xi)-2^{j+1}\frac{A_h\xi}{|\xi|},\frac{\xi}{|\xi|}\right)
 $$
 and
 $$\displaystyle \Xi\left(\tau+\frac 12\left(2^{j+1}+x\cdot
 \frac{A_h\xi}{|\xi|}\right);x,\frac{\xi}{|\xi|}\right)=
 \Xi
\left(\tau;P(x,A_h\xi)-2^{j+1}\frac{A_h\xi}{|\xi|},\frac{\xi}{|\xi|}\right). 
$$
 Since
 $$ 
   \displaystyle
\Bigg|X\left(
\tau;P(x,A_h\xi)-2^{j+1}\frac{A_h\xi}{|\xi|},\frac{\xi}{|\xi|}\right)\Bigg|
 \ge
 2^{j+1}\,\mbox{ for }\,\tau\geq 0, 
$$
 it follows from (\ref{eq5.2.4}) that
 $$
\displaystyle
   p_{e\,j}(x,\xi)=
 p_{e\,j}\left(P(x,A_h\xi)+
 2^{j+1}\psi\left(\frac{x\cdot
 A_h\xi}{|x||\xi|}\right)\frac{A_h\xi}{|\xi|},\xi\right)
$$
 and (\ref{eq5.2.5}) holds.

 For
 the proof of (vii) notice first that $p_{e\,j}\left(1-\chi\left(\frac
 1{10\,2^j}|x|\right)\right)$ has compact support in $x$ and therefore
belongs to
 $S^0_{1,0}$. Thanks to the decay properties of $b_j$ we can summ in $j$
to
obtain
 the same
 property for $q(x,\xi)$. The other conclusions follow by inspection
 from (vi) and the
 definition of $k$.

 Finally notice that (viii) follows from the previous steps.
This completes the  proof of
Lemma \ref{le5.2.3}.

 \bigskip

 \begin{definition} \label{def5.2.4} The equivalence relation $\,\cong\,$ is
 given by

 $$
A_1\cong A_2\quad\Leftrightarrow\quad(A_1-A_2)\mbox{is }\,L^2\mbox{-bounded
 with norm }\le C(R),
$$
where $C(R)$ grows at most polynomially in $R.$

 \end{definition}

 \begin{lemma} \label{le5.2.5} 

 Let $q\in
 S^1_{1,0}$ with $q(\cdot,\xi)\in {\cal S}(\BbbR^n)$
 uniformly in $\xi$. Let 
 $a\in S(\BbbR^n;S^0_{1,0})$ and

 $$ d(x,\xi)=a(P(x,
 A_h\xi);x,\xi)\chi(|\xi|).$$ Then

 \begin{itemize}
 \item[(i)]
 $\,\Psi_{qd}\cong \Psi_q \Psi_d\cong \Psi_d\Psi_q,$

 \item[(ii)]
 $\,\Psi_{qk^R}\cong \Psi_q K^R\cong K^R\Psi_q\cong \Psi_q
(K^R)^{\ast}\cong
 (K^R)^{\ast}\Psi_q,$

 \item[(iii)] $\,\Psi_{q\tilde k^R}\cong \Psi_q
 \tilde K^R\cong \tilde K^R\Psi_q
 \cong \Psi_q  (\tilde K^R)^{\ast}\cong
 (\tilde K^R)^{\ast}\Psi_q,$

 \item[(iv)] $\,i[(K^R)^{\ast} {\cal L}^R-{\cal L}^R
 (K^R)^{\ast}]\cong \Psi_{k^Rb^R},$

 \item[(v)] $\,i[ (\tilde K^R)^{\ast} {\cal
 L}^R-{\cal L}^R  (\tilde K^R)^{\ast}]\cong \Psi_{\tilde k^Rb^R}.$
\end{itemize}

The constants for the above inequalities are bounded
by $R^N$ for some fixed power $N$, and by some fixed number of
seminorms of $q$ and $a.$

 \end{lemma}

 \sl{Proof of Lemma \ref{le5.2.5}} \rm
The proof is based on the calculus developed in Section \ref{sec3}.

(i) Let
$\phi(x)=\left(1+|x|^2\right)^{-N}=\langle x\rangle^{-2N}$, and
$\displaystyle\frac{1}{\phi} q=\tilde q$. Then for $q=\phi \tilde q$ one has
 $$
\begin{array}{c}
 \Psi_q \Psi_d=\Psi_{\tilde q\phi}\Psi_d \cong \Psi_{\tilde q}\phi\Psi_d
  \cong
 \Psi_{\tilde
qd}=\Psi_{qd}.
 \end{array}
$$
 Here it was used that $\,\phi d$ behaves as a classical symbol because for our purposes just finitely many derivatives
of the symbol
in $\xi$ are needed. This number of derivatives determines the choice of $N$.
Similarly,
$$
\begin{array}{c}
 \Psi_d \Psi_q\cong
 \Psi_d\phi(I-\Delta)J^{-2}\Psi_{\tilde q}\\
 \cong
 \Psi_{d\phi\langle\xi\rangle^{2}} \Psi_{\langle
 \xi\rangle)^{-2}\tilde q}\\
 \cong
 \Psi_{d\phi\langle \xi\rangle^{2}\langle
 \xi\rangle^{-2}\tilde q}=\Psi_{\tilde qd}.
 \end{array}
$$
In this case we used that $\,\phi(I-\Delta)\,$ is a partial
 differential operator with decay in the coefficients and
 that $\,\phi(x) \langle \xi\,\rangle^2$ behaves as a 
symbol in $\, S^2_{1,0}\,$ for all $N$.

 (ii) $\,\Psi_{qk^R}\cong\Psi_q K^R\cong
 K^R\Psi_q\,$ follows from the decomposition in Lemma \ref{le5.2.3} (vii) of
 $\,k^R\,$ into a
 sum of an $\,S^0_{1,0}\,$ symbol and a symbol of the type $\,d\,$ in
(i). For
 the
 remainder of (ii),

 $$
   \begin{array}{c}
 \Psi_q (K^R)^{\ast}=(K^R\Psi^{\ast}_q)^{\ast}\cong
 (K^R\Psi_{\overline q})^{\ast}\\
 \cong (\Psi_{k^R\overline q})^{\ast}\cong
 \Psi_{\overline k^R q}=\Psi_{k^R q}
 \end{array}
   $$
 because 
 $\,k^R\overline q\in
 S^1_{1,0}\,$ and $\,k^R\,$ is real--valued. Similarly,

 $$
   \begin{array}{c}
  (K^R)^{\ast}
 \Psi_q=(\Psi^{\ast}_q K^R)^{\ast}\cong (\Psi_{\overline q} K^R)^{\ast}\\
 (\Psi_{\overline q k^R})^{\ast}\cong
 \Psi_{ q k^R}.
 \end{array}
   $$

 (iii) is similar
 to (ii).

 (iv) $\,k^Rb^R\in S^1_{1,0}\,$ and is real--valued so, taking adjoints,
 it suffices to show that
$$i[K^R{\cal L}^R-{\cal L}^R K^R]\cong \Psi_{k^Rb^R}.$$
 Now write
 $${\cal L}^R=\left[{\cal L}^R-(A_h\nabla)\cdot\nabla\right]+
 (A_h\nabla)\cdot\nabla$$
 to
 see that $\,{\cal L}^R\,$ is a compactly supported perturbation of a
constant
 coefficient operator. By
 Theorem \ref{th3.3.4} in Section \ref{sec3} and the decomposition in Lemma \ref{le5.2.3} (vii),
 it follows that
 $$i[K^R{\cal L}^R-{\cal L}^R K^R]\cong 
 \Psi_{\nabla_{\xi}h_2^R\cdot
 \nabla_{x}k^R-
 \nabla_{x}h_2^R\cdot \nabla_{\xi} k^R}\cong\Psi_{b^Rk^R}$$
 using Lemma \ref{le5.2.3} (iv) in the last equivalence.

 (v) is similar to (iv) using Lemma \ref{le5.2.3} (v)
 instead of Lemma \ref{le5.2.3} (iv).

 This completes the proof of Lemma \ref{le5.2.5}.
\bigskip

In order to prove Lemma \ref{le5.2.1} we still need some technical results.
 Recall that $\,{\cal L}={\cal L}^R+{\cal E}^R\,$
and that  $\,K^R\,$ and $b^R$ were given in Definition \ref{def5.2.2}.

\begin{lemma} \label{le5.2.6} There exists $N_0$ large enough such that

\begin{itemize}
\item[(i)]
$\|(K^R)^* u\|_{L^2}=O\left(R^{N_0}\|u\|_{L^2}\right),$

\item[(ii)]
$\|\,i\left[{\cal L}^R,(K^R)^*\right]u+
(K^R)^*\Psi_{b^R} u\,\|_{L^2}\leq c
R^{N_0}\|u\|_{L^2},$

\item[(iii)]
$\|(K^R)^* b_2(x)\nabla \overline u-
b_2(x)\nabla \overline{(K^R)^{\ast} u}\,\|_{L^2}
\leq cR^{N_0}\|u\|_{L^2}.
$

\item[(iv)]
$\|(K^R)^*\Psi_{Im  b_1} u-
\Psi_{Im  b_1} {(K^R)^{\ast} u}\,\|_{L^2}
\leq cR^{N_0}\|u\|_{L^2}.
$

\end{itemize}
\end{lemma}

   \sl{Proof of Lemma \ref{le5.2.6}}\rm

Part (i) follows from the descomposition of Lemma \ref{le5.2.3} (vii) and Theorem \ref{th3.2.1} in Section \ref{sec3}. 
Part (ii) follows from Lemma \ref{le5.2.3} (iv), Theorems \ref{th3.3.4}--\ref{th3.3.6} in Section \ref{sec3} and that
$k^R$
is real.
As for (iii) follows from Theorems \ref{th3.3.4}-\ref{th3.3.6} in Section \ref{sec3} and that $k^R$ is even in
$\xi$. Finally (iv) follows from Lemma \ref{le5.2.5} (ii).
\bigskip

\begin{lemma} \label{le5.2.7} 

Given $\,M$ there are $R$ large enough and $N(M)$
such
that
$$
|\langle i\left[{\cal E}^R,(K^R)^{\ast}\right] u,(K^R)^{\ast}
u\rangle|\le
R^{-M}\|J^{1/2}(u)\langle x\rangle^{-N}\|^2_{L^2}+
O\left(R^{N(M)=N_0}\|u\|^2_{L^2}\right),
$$
   with $N_0$ as in Lemma \ref{le5.2.6}.
\end{lemma}
   
   \begin{remark} \label{renuevo}
   Here $M$ can be taken arbitrary large since the coefficients are in $\cal S$. However it suffices to assume (\ref{eq4.1.16}) in Section \ref{sec4} for some $\tau$ sufficiently large.
   \end{remark}

\sl{Proof of Lemma \ref{le5.2.7}}\rm

We have
$$i[{\cal E}^R, (K^R)^{\ast}]=
i\left({\cal E}^R (K^R)^{\ast}-(K^R)^{\ast} {\cal E}^R\right).
$$
Take adjoints to get
\begin{equation}\label{eq5.2.6}
   -i\left(K^R ({\cal E}^R)^{\ast}-({\cal E}^R)^{\ast} K^R\right).
\end{equation}
Recall that
$$
\displaystyle{\cal E}^R=-({\cal L}-{\cal L}^R)=\sum_{j,k}\frac{\partial}
{\partial x_j}\left(e_{jk}\frac{\partial}
{\partial x_k}\right),
$$
with 
$$
   \displaystyle
e_{jk}=\left(1-\theta(\frac{x}{R})\right)\left(a_{jk}(x)-a_{jk}^0\right).
$$
Thus we can see ${\cal E} ^R$ as a second order differential operator with
coefficients
of the form $\displaystyle e_{jk}(x)=\,\frac 1{R^M}\tilde e_{jk}^R(x)\,$,
and $\tilde e_{jk}^R$ with decay uniform in $R$.
Then
$$
   \displaystyle
({\cal E}^R)^{\ast}=-e_{jk}^R(x)\frac{\partial^2}{\partial x_j
\partial x_k}+\vec {b_3}\cdot \nabla, 
$$
for some 
$\vec{b_3}$ with the right decay. Therefore this
term gives
bounds of the type
$$ 
O\left(R^{N_0}\|u\|^2_{L^2}\right)
$$
in (\ref{eq5.2.6}), just using Lemma \ref{le5.2.3} (ii).

Now we use Theorem \ref{th3.3.4} in Section \ref{sec3} and 
Lemma \ref{le5.2.3} (iv) to get
$$
K^R({\cal E}^R)^{\ast}-({\cal E}^R)^
{\ast} K^R
=\Psi_{\alpha^R}+\mbox{zero order terms},$$
with 
$$\alpha^R(x,\xi)=\sum_{j,k}(\nabla_x
e_{jk}\xi_j\xi_k\cdot
\nabla_{\xi} k^R
-e^R_{jk}\nabla_{\xi}(\xi_j\xi_k)\cdot \nabla_x k^R),
$$
and the zero order terms have bounds in $L^2$ which grow as $R^{N_0}.$

Notice that  $\alpha^R(x,\xi)\,\in S^1_{1,0}$,
with uniform $O(\,R^{-M}\,)$ decay in $\,x\,$. 
We need to study 
$\,\langle \Psi_{\alpha^R}^{\ast} u,(K^R)^{\ast} u\rangle.$ 

 We observe that
$\,\Psi_{\alpha^R}^{\ast}\,=\Psi_{\overline{\alpha^R}}+ \mbox{zero
order}$. Therefore we will work with $\Psi_{{\alpha^R}}$, being analogous
the
calculations for $\Psi_{\overline{\alpha^R}}$. Define $\phi(x)=\langle
\,x\,\rangle^{-N}$. Then
\begin{equation}\label{eq5.2.7}
\begin{array}{rcl}
R^M\langle \Psi_{{\alpha^R}}\,u\, ,\, (K^R)^\ast
\,u\,\rangle&=&R^M\langle \phi \phi^{-2}
\Psi_{{\alpha^R}}J^{1/2}J^{-1}J^{1/2}\,u\, ,\,\phi (K^R)^\ast
\,u\,\rangle\\
&\cong&\langle \phi J^{1/2}\phi^{-2}
\Psi_{{\alpha^R}}J^{-1}J^{1/2}\,u\, ,\,\phi (K^R)^\ast
\,u\,\rangle\\
&\cong&\langle \phi J^{1/2}\phi\phi^{-3}
\Psi_{{\alpha^R}}J^{-1}J^{1/2}\,u\, ,\,\phi (K^R)^\ast
\,u\,\rangle\\
&\cong&\langle \phi J^{1/2} \phi^{-3} \,
\Psi_{{\alpha^R}}J^{-1}\phi\, J^{1/2}\,u\, ,\,\phi (K^R)^\ast
\,u\,\rangle\\
&\cong&\langle  J^{1/2} \phi\,\phi^{-3} \,
\Psi_{{\alpha^R}}J^{-1}\phi\, J^{1/2}\,u\, ,\,\phi (K^R)^\ast
\,u\,\rangle\\
&\cong&\langle   \,\phi^{-3} \, \Psi_{{\alpha^R}}J^{-1}\phi\,
J^{1/2}\,u\, ,\,\phi
\,J^{1/2}\phi(K^R)^\ast
\,u\,\rangle\\
&\leq& \| \,\,\phi^{-3} \, \Psi_{{\alpha^R}}J^{-1}\phi\,
J^{1/2}\,u\,\|_{L^2}\|\,\phi
\,J^{1/2}\phi(K^R)^\ast
\,u\,\|_{L^2}\\
&\leq &C \|\phi\, J^{1/2}
\,u\,\|_{L^2}\,\|\phi\, J^{1/2}\phi(K^R)^\ast
\,u\,\|_{L^2}\,
\end{array}
\end{equation}
as desired. Lemma \ref{le5.2.7} is proved.
\bigskip

Next, with $\,\tilde K^R\,$ as in Definition \ref{def5.2.2}, we define

$$E^R=I-\tilde K^R (K^R)^{\ast}.$$

\begin{lemma} \label{le5.2.8} 

 There exists $N_0$ such that

\begin{itemize}
\item[(i)]
$\|E^Ru\|_{L^2}\le C
R^{N_0}\|u\|_{L^2},$

\item[(ii)] Let $q\in \,S^1_{1,0}$ with $q(\,.\,,\xi)\,\in{\cal S}(\BbbR^n)$ uniformly in $\xi$. Then
$$\|\Psi_q E^Ru\|_{L^2} +\|E^R\Psi_qu\|_{L^2}\le C R^{N_0}\|u\|_{L^2}.$$
As a consequence
$$
\|E^R \Psi_{b_1} u\|_{L^2}+\|E^R \vec b_2.\nabla \overline
u\|_{L^2}\le C R^{N_0}\|u\|_{L^2},
$$

\item[(iii)]
$|\,\left\langle [{\cal L}^R,E^R]u, E^R u\right\rangle|\le C
R^{N_0}\|u\|^2_{L^2},$

\item[(iv)]
$
|\left\langle i\left[{\cal L},E^R\right]u, E^R u\right\rangle|\le
C R^{N_0}\|u\|^2_{L^2}.$
\end{itemize}
\end{lemma}

\sl{Proof of Lemma \ref{le5.2.8}}\rm
   
   Part (i) follows from the descomposition of Lemma \ref{le5.2.3} (vii) and Theorem \ref{th3.2.1} in Section \ref{sec3}.

As for (ii)

$$
\Psi_q E^R=\Psi_q -\Psi_q \tilde K^R(K^R)^\ast \,\cong\,\Psi_q -\Psi_{q
\tilde k}
(K^R)^\ast\,\cong\,\Psi_q -\Psi_{q \tilde k k}=0,
$$
by using Lemma \ref{le5.2.3} (iii) and (ii). Similarly
$$
E^R\Psi_q =\Psi_q -\tilde K^R(K^R)^\ast \Psi_q \,\cong\,\Psi_q -\tilde
K^R\Psi_{q 
k}
\cong\,\Psi_q -\Psi_{q \tilde k k}=0.
$$

For (iii) it is enough to prove that $[{\cal
L}^R,E^R]u$ is $L^2$-bounded. But
$$\begin{array}{c}
i[{\cal
L}^R,E^R]=i\left({\cal
L}^R\tilde K^R(K^R)^\ast-\tilde K^R(K^R)^\ast{\cal
L}^R\right)\\
=i\left(({\cal
L}^R\tilde K^R-\tilde K^R{\cal
L}^R)(K^R)^\ast+i\tilde K^R(\,{\cal L}^R(K^R)^\ast-(K^R)^\ast{\cal
L}^R\,)\right)\\
\,\cong\,\Psi_{(\tilde k^Rb^R)}\,(K^R)^\ast-\tilde K^R\Psi_{(k^Rb^R)}\\
\,\cong\,\Psi_{ (\tilde k^Rb^Rk^R)}-\Psi_{(k^Rb^R\tilde k^R)}=0.\\
\end{array}
$$

Finally let us prove (iv).
We have
${\cal L}={\cal L}^R+{\cal E}^R,$ and $E^R=I-\tilde K^R (K^R)^{\ast}.$
Part (iii) gives that
$\,\left\langle [{\cal
L}^R,E^R]u, E^R u\right\rangle\,$ has the right bound. Hence we need to
understand 
$[{\cal E}^R,\tilde K^R (K^R)^{\ast}].$

We have
$$
\begin{array}{c}
{\cal E}^R \tilde K^R (K^R)^{\ast}-\tilde K^R (K^R)^{\ast}{\cal E}^R=\\
\left({\cal E}^R \tilde K^R -\tilde K^R {\cal E}^R\right) (K^R)^{\ast}+
\tilde K^R \left({\cal E}^R  (K^R)^{\ast} -(K^R)^{\ast} {\cal E}^R\right). 
\end{array}
$$
For both 
$\left({\cal E}^R \tilde K^R -\tilde K^R {\cal E}^R\right)$  and $
\left({\cal E}^R  (K^R)^{\ast} -(K^R)^{\ast} {\cal E}^R\right)$ we can use
Theorems \ref{th3.3.4}--\ref{th3.3.6} in Section \ref{sec3} 
so that they can be written as $\Psi_{\beta_j} + \mbox{zero order
terms}\,$ with $\Psi_{\beta_j},\,\, j=1,2\,$ classical first order pseudo-differential operators with the right decay in $x$ as  we did in the
proof of
Lemma \ref{le5.2.7}. But from Lemma \ref{le5.2.5}
$$ \Psi_{\beta_1} (K^R)^{\ast}\,\cong\,\Psi_{\beta_1k^R},\qquad\qquad 
\tilde K^R\Psi_{\beta_2} \,\cong\,\Psi_{\beta_2\tilde k^R},$$
and we can apply part (ii).

\bigskip

\sl{Proof of Lemma \ref{le5.2.1}}\rm
   
Notice first that part A follows making $f=0$ in part B and using 
Duhamel's principle.
We study the problem
$$
\left\{
   \begin{array}{l}
\partial_t u=i{\cal L} u+\Psi_{b_1} u+\vec {b_2}(x)\cdot \nabla
\overline u+\Psi_{c_1} u+\Psi_{c_2}\overline u+f\\
 u(x,0)=u_0.
\end{array}
\right.
   $$
We also have  ${\cal L}={\cal L}^R+{\cal E}^R\,.$ Applying the operator
$(K^R)^{\ast}$ and observing that if $s=0$ in Definition \ref{def5.2.2}  then $\psi_{b^R}=\psi_{-\mbox{Re}\, b_1})$, we
get from Lemma \ref{le5.2.6}
$$
\begin{array}{rcl}
\partial_t (K^R)^{\ast} u&=&\displaystyle
i\left[{\cal L}^R,(K^R)^{\ast}\right] u+(K^R)^{\ast}\Psi_{b_1} u+
i {\cal L}^R (K^R)^{\ast} u\\
&+&\displaystyle i{\cal E}^R (K^R)^{\ast} u+
i\left[{\cal E}^R,(K^R)^{\ast}\right] u
+(K^R)^{\ast} \vec b_2(x)\cdot\nabla\overline u\\
&+&  \displaystyle(K^R)^{\ast}f+
\mbox
{zero order terms}\\
&=&\displaystyle i{\cal L} (K^R)^{\ast} u+\vec b_2(x)\cdot\nabla\overline
{(K^R)^{\ast} u}+
i\left[{\cal} E^R,(K^R)^{\ast}\right] u
\\
 &+& \displaystyle i\Psi_{Im b_1}(K^R)^{\ast}u +(K^R)^{\ast}f+
\mbox
{zero order terms}.
\end{array}
$$
Define $v^R=(K^R)^{\ast} u\,$.
Then we have 
$$
\begin{array}{c}
\partial_t\langle v^R,v^R\rangle=
i\langle {\cal L} v^R,v^R\rangle+
i\left\langle \left[{\cal E}^R,(K^R)^{\ast}\right] u,v^R\right\rangle+
\left\langle \vec b_2(x)\cdot\nabla\overline{v^R},v^R\right\rangle +
\left\langle i\Psi_{Im b_1}v^R,v^R\right\rangle\\
+O\left(R^{N_0}\|u\|_{L^2}\|v^R\|_{L^2}\right)+
O\left(R^{N_0}\|u\|^2_{L^2}\right)+ O(\|f\|^2_{L^2})\\
=i\langle {\cal L} v^R,v^R\rangle+
i\left\langle \left[{\cal E}^R,(K^R)^{\ast}\right]
u,v^R\right\rangle+
\left\langle i\Psi_{Im b_1}v^R,v^R\right\rangle\\
+O\left(R^{N_0}\|u\|_{L^2}\|v^R\|_{L^2}\right)
+O\left(R^{N_0}\|u\|^2_{L^2}\right)+ O(\|f\|^2_{L^2}),
\end{array}
$$
where the last step follows by integration by parts.
Taking the real part of both sides and using Lemma \ref{le5.2.7} and Garding's inequality
for $\mbox{Re}\, \left\langle i\Psi_{Im b_1}v^R,v^R\right\rangle$ we get after integration in the temporal variable
$$
\begin{array}{c}
\sup_{0\le t\le T}\|v^R(t)\|^2_{L^2}\le
CR^{N_0}\|u(0)\|^2_{L^2}\\ \\
+R^{-M}\int_0^T
\|J^{1/2}(u)\langle x\rangle^{-N}\|^2_{L^2}dt+
CR^{N(M)+N_0}T\sup_{0<t<T}\|u(t)\|^2_{L^2}+C\int_0^T\|f\|^2_{L^2}
\end{array}
$$
for $\,T\,$ small.

Similarly
$$
\partial_t E^R u=i{\cal L} E^R u+i[{\cal L},E^R]u
+E^R\Psi_{ b_1} u+E^R \vec b_2\cdot\nabla \overline u +E^Rf.
 $$
Then from Lemma \ref{le5.2.8} we have
$$
\sup_{0<t<T}\|E^R u\|^2_{L^2}\le C R^{N_0}\|u(0)\|^2_{L^2}+
C R^{N_0} T\sup_{0<t<T}\|u(t)\|^2_{L^2}+C\int_0^T\|f\|^2_{L^2}.
$$
 Now using $I=E^R+\tilde K^R(K^R)^\ast$  and the previous estimates we get
$$
   \begin{array}{c}
\displaystyle\sup_{0<t<T}\|u(t)\|^2_{L^2}\le 
C R^{N_0}\sup \|v\|^2_{L^2}+
C \sup \|E^R u\|^2_{L^2} \\
\displaystyle\le C R^{N_0}\left\{C R^{N_0}\|u(0)\|^2_{L^2}+
R^{-M}\int_0^T \|J^{1/2}(u)\langle x\rangle^{-N}\|^2_{L^2}
dt+\right.\\\displaystyle+\left.R^{N(M)}_T\sup_{0<t<T}\|u(t)\|^2_{L^2}\right\}+
CR^{N_0} T\sup_{0<t<T}\|u(t)\|^2_{L^2}+C\int_0^T\|f\|^2_{L^2}.
\end{array}
   $$
Since $\,N_0\,$ is fixed and $\,M\,$ arbitrary, this combined  with Lemma \ref{le5.2.8} finishes the proof.

\bigskip

\sl{Proof of Theorem \ref{th5.1.2}}\rm

Parts A and B follow by the a priori estimates in Lemma \ref{le5.1.4} and Lemma \ref{le5.2.1} just as we
did in Theorem 4.1 of Section \ref{sec2} for the elliptic case.

So just part C remains to be proved. It is enough to prove
the case s=0 by Lemma \ref{le5.1.3}. It can be assumed that $u_0=0$ since (B) solves the case $f=0$ and that
$f\in{\cal S}(\BbbR^{n+1})$. Let $u$ be the solution of (\ref{eq5.1.2}) given by (B) and let $\varphi
\in{\cal S}(\BbbR^{n})$. The family of problems (\ref{eq5.1.2}) is invariant under
adjoints and time reversal, so (B) yields a solution $v$ of
$$ 
\left\{ 
\begin{array}{l}
\p_t v=-\big (i {\cal L}  + \Psi_{b_1}+\vec b_2(x)\cdot\nabla \bar{ (\cdot)}
+\Psi_{c_1} +\Psi_{c_2}\bar{ (\cdot)} \big)^\ast v,\\
v(x,0)=v_0(x). 
\end{array} 
\right.
 $$
Now
$$
\begin{array}{c}
\displaystyle\partial_t\langle u(T),v(T)\rangle=\langle u(0),v(0)\rangle+\int_0^T\big(\langle \p_tu(t),v(t)\rangle+
\langle u(t),\p_tv(t)\rangle\big)\,dt\\
\displaystyle
i{\cal L}u  + \Psi_{b_1}u+\vec {b_2}(x)\cdot\nabla \bar u
+\Psi_{c_1}u +\Psi_{c_2}\bar{ u}+f,v\rangle\\
+\displaystyle\int_0^T\langle u,-\big (i {\cal L}  + \Psi_{b_1}+\vec b_2(x)\cdot\nabla \bar{ (\cdot)}
+\Psi_{c_1} +\Psi_{c_2}\bar{ (\cdot)} \big)^\ast v\rangle\,dt\\
\displaystyle=\int_0^T\langle f,v\rangle\,dt=\int_0^T\big\langle\langle x\rangle^{N/2}\, J^{-1/2}f, 
\langle x\rangle^{-N/2}\,J^{1/2}v\big\rangle\,dt\\
\displaystyle\le \big(\int_0^T\int_{\BbbR^n}|J^{-1/2}f|^2\langle x\rangle^{N}\,dxdt\big)^{1/2}
\big(\int_0^T\int_{\BbbR^n}|J^{1/2}v|^2\langle x\rangle^{-N}\,dxdt\big)^{1/2}\\
\displaystyle\le C\big(\int_0^T\int_{\BbbR^n}|J^{-1/2}f|^2\langle x\rangle^{N}\,dxdt\big)^{1/2}\,\|\varphi\|_{L^2},
\end{array}
$$
where the last inequality follows from (B). Hence
$$\displaystyle\sup_{[0,T]}\|u\|^2_{L^2}\le C\int_0^T\int_{\BbbR^n}|J^{-1/2}f|^2\langle x\rangle^{N}\,dxdt.
$$
Finally use Lemma \ref{le5.1.4}, and the proof of Theorem \ref{th5.1.2} is complete.

\newpage

\section{NONLINEAR EQUATIONS}\label{sec6}

The local smoothing results for linear Schr\"odinger equations will now 
be applied to obtain local well-posedness in weighted Sobolev spaces with 
high Sobolev index for a quite general class of nonlinear Schr\"odinger equations 
with initial data 
in the Schwartz class ${\cal S}(\BbbR^n)$. This follows the contraction
 mapping scheme as in \cite{18} and \cite{22}.
 
\subsection{Linear solutions and weights}\label{subseq6.1}

Suppose $\,u\,$ is a solution of the linear Schr\"odinger equation
\begin{equation}\label{eq6.1.1}
\left\{
\begin{array}{l}
\partial_t u=i{\cal L} u+\vec b_1\cdot \nabla u+
\vec {b_2}\cdot \nabla 
\overline u+c_1 u+c_2\overline u+f,\\
u(x,0)=u_0(x).
\end{array}
\right.
\end{equation}

By Lemma  \ref{le5.1.3} and lemma \ref{le5.1.4} in Section \ref{sec5}, $\,w=J^s u\,$ is a solution of the equation
\begin{equation}\label{eq6.1.2}
\left\{
\begin{array}{l}
\partial_t w=i{\cal L} w+\Psi_{b_3}w+\vec {b_2}\cdot \nabla
\overline w+\Psi_{r_1}w+\Psi_{r_2}\overline w+J^s f\\
 w|_{t=0}=J^s u_0.
\end{array}
\right.
\end{equation}
where 
$$
\displaystyle\Psi_{b_3}=\vec {b_1}\cdot \nabla-i s\sum\limits_{j,k,l}\partial_{x_j}a_{kl}
\partial^3_{x_j x_k x_l}J^{-2},\,\,\,r_1,r_2\in S^0_{1,0}.
$$
Let $W_j\,$ denote the solution operator of (6.$j$), $j=1,2$, with $f=0\,$. 
Then one has $J^sW_1 u_0=W_2 J^s u_0\,$. A similar result is needed when $J^s\,$ is replaced by the weight
$(1+|x|^2)^N\,$. It is useful to obtain a few results concerning commutators of 
weights and classical pseudo-differential operators.

\begin{lemma} \label{le6.1.1} 
 
 Let $\,p\in S^m_{1,0}\,,\, \alpha\in\BbbN^n_0\,$. Then
$$
\displaystyle
x^{\alpha}\Psi_p f-\Psi_p\left[x^{\alpha} f\right]=
\sum_{0<\beta\le \alpha}
\left(
\begin{array}{c}
 {\alpha}\\
 {\beta}
 \end{array}
 \right)
 \Psi_{i^{\beta}\partial^{\beta}_{\xi}p}
\left[x^{\alpha-\beta} f\right]\,,\,f\in {\cal S}.
$$
\end{lemma}
\sl{Proof of Lemma \ref{le6.1.1}}\rm
 
Using integration by parts and Leibniz\rq{} rule,
$$
\begin{array}{l}
\displaystyle x^{\alpha}\Psi_p f(x)=
\int x^{\alpha} e^{i x\cdot \xi} p(x,\xi)\hat f(\xi)d\xi\\
\displaystyle=\int (-i \partial_{\xi})^{\alpha}\left[e^{i x\cdot \xi}\right] p(x,\xi)\hat f(\xi)d\xi=
\int e^{i x\cdot \xi} (i \partial_{\xi})^{\alpha}
\left[p(x,\xi)\hat f(\xi)\right]d\xi\\
\displaystyle=\int e^{i x\cdot \xi} \sum_{\beta\le \alpha}
 \left(
 \begin{array}{c}
 {\alpha}\\{\beta}
 \end{array}
 \right)
(i \partial_{\xi})^{\beta}
 p(x,\xi)
 (i\partial_{\xi})^{\alpha-\beta}\hat f(\xi)d\xi\\
\displaystyle=\int e^{i x\cdot \xi} 
\sum_{\beta\le \alpha}
  \left(
 \begin{array}{c}
 {\alpha}\\{\beta}
  \end{array}
 \right)
(i \partial_{\xi})^{\beta} p(x,\xi) 
(x^{\alpha-\beta} f)^{\hat{}}(\xi)d\xi\\
\displaystyle=\Psi_p\left[x^{\alpha} f\right](x)+\sum_{0<\beta\le \alpha}
 \left(
 \begin{array}{c}
  {\alpha}\\{\beta}
  \end{array}
 \right)
\Psi_{i^{\beta} \partial_{{\xi}}^{\beta}p} 
\left[x^{\alpha-\beta} f\right](x).
\end{array}
$$
This proves Lemma \ref{le6.1.1}.
\bigskip

 \begin{lemma} \label{le6.1.2} 

Let $p\in S^m_{1,0}\,,\, N\in\BbbN\,$. Then
$$
\begin{array}{c}
\displaystyle(1+|x|^2)^N\Psi_p f=\Psi_p\left[(1+|x|^2)^N f\right]+2N\sum_j
\Psi_{i\partial_{\xi_j}p}\left[x_j(1+|x|^2)^{N-1}f\right]\\
\displaystyle+\sum_{|\alpha+\beta|\le 2N,|\alpha|\ge 2,|\beta|\le 2N-2}c_{\alpha\beta}
\Psi_{\partial^{\alpha}_{\xi}p}\left[x^{\beta}f\right]\,,\;\;\;\,f\in {\cal S}(\BbbR^n).
\end{array}
$$
\end{lemma}
\sl{Proof of Lemma \ref{le6.1.2}}\rm
 
Lemma \ref{le6.1.1} yields the identity
$$\displaystyle(1+|x|^2)\Psi_p -\Psi_p(1+|x|^2)=2 \sum_j\Psi_{i\partial_{\xi_j^p}}x_j-\Psi_{\Delta_{\xi}p}.
$$
Lemma \ref{le6.1.2} follows by induction of $\,N\,$ and further applications of Lemma \ref{le6.1.1}.
\bigskip
 
Now we study weighted Sobolev norms of solutions of (\ref{eq6.1.1}) and (\ref{eq6.1.2}). Let $\,u(t)=W_1(t)u_0\,$ be the
solution of the linear equation
$$
\left\{
 \begin{array}{l}
 \partial_t u=i{\cal L} u+\Psi_{b_1}u+\vec b_2
\cdot \nabla\overline u+\Psi_{c_1}u+\Psi_{c_2}\overline u\\
 u(x,0)=u_0(x).
\end{array}
\right.
$$
 
\begin{lemma} \label{le6.1.3} 

 Let $\,N\in \BbbN\,,\,s\in \BbbR\,$. Suppose 
$\,\langle x\rangle^{2N}u_0\in H^{s+2N}\,.$ Then
\begin{equation}
\displaystyle\sup_{0\le t\le T}\|\langle x\rangle^{2N}W_1(t)u_0\|^2_{H^s}
\le\sum_{j=0}^{2N}c_jT^j\|\langle x\rangle^{2N-j}u_0\|^2_{H^{s+j}},
\label{a}\end{equation}
and
\begin{equation}\sup_{0\le t\le T}\|\langle x\rangle^{2N}W_1(t)u_0\|^2_{H^s}\le
c(1+T^{2N})\|\langle x\rangle^{2N}u_0\|^2_{H^{s+2N}}\,.\label{b}
\end{equation}

\end{lemma}

\sl{Proof of Lemma \ref{le6.1.3}}\rm 

 Let $\,u=W_1u_0\,$. Then
$$
\begin{array}{l}
\partial_t\left[\langle x\rangle^{2N}u\right]=\langle x\rangle^{2N}\partial_t u\\
=\langle x\rangle^{2N}\left\{i{\cal L} u+\Psi_{b_1}u+\vec b_2\cdot \nabla
\overline u+\Psi_{c_1}u+\Psi_{c_2}\overline u\right\}.
\end{array}
$$
Using Lemma \ref{le6.1.2}, this equals to
$$
\begin{array}{l}
i{\cal L}\langle x\rangle^{2N}u+\Psi_{b_1}\langle x\rangle^{2N}u+
\vec b_2\cdot \nabla \langle x\rangle^{2N}\overline u+\Psi_{\tilde c_1}\langle x\rangle^{2N}u+
\\
\displaystyle+\Psi_{\tilde c_2}\langle x\rangle^{2N}\overline u+
+2N\sum_j\Psi_{\partial_{\xi_j}h_2}x_j\langle x\rangle^{2N-2}u.
\end{array}
$$
Hence $\,\langle x\rangle^{2N} W_1 u_0\,$ satisfies a linear equation with initial data $\,\langle 
x\rangle^{2N} u_0\,$ and forcing term

$$
\displaystyle f=2N\sum_j\Psi_{\partial_{\xi_j}h_2}x_j\langle x\rangle^{2N-2}W_1 u_0.
$$
Hence,
$$
\begin{array}{c}
\displaystyle\sup_{0\le t\le T}\|\langle x\rangle^{2N} W_1(t) u_0\|^2_{H^s}\le
c\|\langle x\rangle^{2N} u_0\|^2_{H^s}\\
\displaystyle+c\int_0^T\sum_j\|x_j\langle x\rangle^{2N-2} W_1(t) u_0\|^2_{H^{s+1}}dt\\
\le c\Big\{\|\langle x\rangle^{2N} u_0\|^2_{H^s}+
T\sum_j \sup_{0\le t\le T}\|x_j\langle x\rangle^{2N-2} W_1(t) u_0\|^2_{H^{s+1}}\Big\}.
\end{array}
$$
Using Lemma \ref{le6.1.1} instead of Lemma \ref{le6.1.2} and arguing as above, it follows that 
$\,x_j\langle x\rangle^{2N-2} W_1 u_0\,$ satisfies a linear equation with initial data 
$\,x_j\langle x\rangle^{2N-2} u_0\,$ and forcing term
$$
f=\Psi_q \langle x\rangle^{2N-2}W_1u_0\,,
$$
where $\,q\in S^1_{1,0}\,$. Therefore,
$$
\begin{array}{l}
\displaystyle\sup_{0\le t\le T}\|x_j \langle x\rangle^{2N-2} W_1(t) u_0\|^2_{H^{s+1}}\\
\displaystyle\le c \|x_j \langle x\rangle^{2N-1} u_0\|^2_{H^{s+1}}+c\int_0^T \|\Psi_q 
\langle x\rangle^{2N-2} W_1(t) u_0\|^2_{H^{s+1}}dt\\
\displaystyle\le c \|\langle x\rangle^{2N-1} u_0\|^2_{H^{s+1}}+cT\sup_{0\le t\le T} \|\langle x\rangle^{2N-2} W_1(t) u_0\|^2_{H^{s+2}},
\end{array}
$$
and
$$
\begin{array}{c}\displaystyle\sup_{0\le t\le T} \|\langle x\rangle^{2N} W_1(t) u_0\|^2_{H^{s}}\le
c\|\langle x\rangle^{2N} u_0\|^2_{H^{s}}
\\+c_1T\|\langle x\rangle^{2N-1} u_0\|^2_{H^{s+1}}+
c_2T^2\|\langle x\rangle^{2N-2} W_1(t) u_0\|^2_{H^{s+2}}.
\end{array}
$$
Now apply this result $\,N-1\,$ times with $\,N\,$ replaced by $\,N-1\,,\,N-2\,,\dots,1\,,$ and the
proof of (\ref{a}) is complete.

For (\ref{b}), it suffices to note that
$$
\begin{array}{l}
\|\langle x\rangle^{2N-j} u_0\|_{H^{s+j}}=
\|\langle x\rangle^{-j}\langle x\rangle^{2N} u_0\|_{H^{s+j}}\\
\le c \|\langle x\rangle^{2N} u_0\|_{H^{s+j}}
\le c \|\langle x\rangle^{2N} u_0\|_{H^{s+2N}}
\end{array}
$$
if $\,j\in\{0,1,\dots,2N\}\,$ since $\,\langle x\rangle^{-j} \in S^0_{1,0}$. This proves Lemma \ref{le6.1.3}.
\bigskip

\subsection{The nonlinear Cauchy problem}\label{subseq6.2}

Consider the initial value problem 
\begin{equation}\label{eq6.2.1}
\left\{
\begin{array}{l}
\partial_t u=i{\cal L} u+\vec b_1\cdot \nabla u+
\vec b_2\cdot \nabla\overline u+c_1 u+c_2\overline u+
P(u,\nabla u,\overline u,\nabla \overline u),\\
u(x,0)=u_0(x),
\end{array}
\right.
\end{equation}
where
$$
\displaystyle{\cal L} u(x)=\sum\limits_{j,k}\partial_{x_j}(a_{jk}(x)\partial_{x_k}u(x)),\;\; \;\;\;A(x)=(a_{jk}(x))_{j,k=1,\dots,n}
$$ 
is a real, symmetric $n\times n\,$ matrix, and 
$P\,$ is any polynomial with no linear or constant terms. 
Concerning the variable coefficients, it will be assumed that
$$
a_{jk}\,,\vec b_1\,,\,\vec b_2\,,\, c_1\,,\,c_2\in C^{\infty}_b.
$$
Assume further that the matrix $\,A(x)=(a_{jk}(x)))_{j,k=1,\dots,n}\,$ is positive definite or
invertible. There will be additional hypotheses on $\,a_{jk}\,,\,\vec b_1\,$ and
$\,\vec b_2\,$ in each of the two cases. More precisely,

\bf Elliptic case \rm

Suppose
$$
\nu^{-1}|\xi|^2\le |A(x)\xi\cdot \xi|\le \nu |\xi|^2,\,\;\;\;\;\,x,\xi\in\BbbR^n.
$$
Then assume in addition the following:
\begin{itemize}
\item[(a)] $A(x)\,$ generates a bicharacteristic flow with non-trapped bicharacteristics.
\item[(b)] There exist $N>1$ and a constant $C$ such that if $\lambda(|x|)=\langle x \rangle^{-N}$ 
$$
|\nabla a_{jk}(x)|\,,\,|\mbox{Im}\,\vec b_1(x)|\le C\lambda(|x|)\,,\;\;\;\;\;\,x\in\BbbR^n.
$$
\end{itemize}

\bf{Ultrahyperbolic case}\rm 

Suppose

$$
\nu^{-1}|\xi|\le |A(x)\xi|\le \nu |\xi|\,,\,x,\xi\in\BbbR^n.
$$
Then assume in addition the following:
\begin{itemize}
\item[(a)] $a_{jk}(x)-a^0_{jk}\in{\cal S}(\BbbR^n)$ for $\,j,k=1,..,n$
where $A^0=(a^0_{jk})$ is a real symmetric  $n\times n$ constant matrix.
\item[(b)] The bicharacteristics are non-trapped.
\item[(c)]   $\vec b_1, \,\vec b_2\in {\cal S}(\BbbR^n:\BbbC^n).$
\end{itemize}

\begin{remark} \label{re6.2.1} It will be clear from our previous results and the proof below that assumptions in (a) and (c)
in the hyperbolic case can be relaxed so that one only needs a finite number of seminorms in (\ref{eq3.1.2}) in Section \ref{sec3}.
\end{remark}

Under the above assumptions, the following result holds.

\begin{theorem} \label{th6.2.2}

 Let $u_0\in {\cal S}(\BbbR^n)$ and $s\ge n+4N+13$. Then there exists \newline
$T=T(\|u_0\|_{H^s},\|\langle x\rangle ^{2N}J^{s-3/2}u_0\|_{L^2})$ such that (\ref{eq6.2.1}) 
has a unique solution $\,u\,$ defined in the time 
interval $\,\left[0,T\right]\,$ satisfying
$$
u\in C^{\infty}\left([0,T]:H^s(\BbbR^n)\cap L^2\left(\langle x\rangle^Ndx\right)\right).
$$
Let
$$
X^s_T=\left\{u\in C([0,T];H^s(\BbbR^n))\,:\,\max_{j=1,2,3}\lambda_j(w)<\infty\right\}
$$
where
$$
\begin{array}{l}
\displaystyle\lambda_1(w)=\sup_{0\le t\le T}\|w\|_{H^s}\\
\displaystyle\lambda_2=\int_0^T\int_{\BbbR^n}|J^{s+1/2} w(x,t)|^2 \langle x\rangle^{-2N}dxdt\\
\displaystyle\lambda_3(w)=\sup_{0\le t\le T}\|\langle x\rangle^{2N} \partial_t w(t)\|_{H^{s/2+n/2+3}}.
\end{array}
$$
Then for every $\,u_0\in {\cal S}(\BbbR^n)\,$ there exists a neighborhood $\,{\cal U}\,$ of $\,u_0\,$ in $\,S\,$ and a 
$\,T^{\prime}>0\,$ such that the map data $\to$ solution of (\ref{eq6.2.1}) is continuous from $\,{\cal U}\,$ into
$\,X^s_{T^{\prime}}$. 
\end{theorem}

\begin{remark}\label{re6.2.3}  The classical pseudo-differential theory in subsection 2.1 of Section \ref{sec2} and the new
operator calculus in Section \ref{sec3} both basically rely on Taylor expansions of finite order and
finitely many integrations by parts. Consequently, by Sobolev\rq{}s theorem, the assumption $\,u_0\in {\cal S}\,$ 
in Theorem \ref{th6.2.2} can be relaxed to $\,\langle x\rangle^{2N}J^{s_1}u_0\in L^2\,$ for some large
$\,s_1=s_1(n,s)\,$. The solution $\,u\,$ of (\ref{eq6.2.1}) is then in $\,C\left([0,T]:H^s(\BbbR^n))\cap L^2\left(\langle x\rangle^Ndx\right)\right)\,$. 
It is an interesting problem to determine the optimal regularity and decay of $\,u_0\,$ needed in specific examples of
equation (\ref{eq6.2.1}).
\end{remark}

\sl{Proof of Theorem \ref{th6.2.2}}\rm 

Let $\,s_0\in 2\BbbN\,$ with $\,s_0+\frac 12\ge n+4N+13\,$. Let
$\,v=J^{s_0}u\,$. Then $\,u\,$ solves (\ref{eq6.2.1}) if and only if $\,v\,$ solves
$$
\left\{
\begin{array}{l}
\partial_t v=i{\cal L} v+\Psi_{\tilde b_1}v+\vec b_2\cdot \nabla
\overline v+\Psi_{\tilde c_1}v+\Psi_{\tilde c_2}\overline v\\
\quad +J^{s_0} \left[P(J^{-s_0}v,\nabla J^{-s_0}v,J^{-s_0}\overline v,\nabla J^{-s_0}\overline v)\right]\\
 v(x,0)=J^{s_0} u_0(x).
\end{array}
\right.
$$
Using Leibniz\rq{} rule,
$$
\begin{array}{l}
J^{s_0} [P(J^{-s_0}v,\nabla J^{-s_0}v,J^{-s_0}\overline v,\nabla J^{-s_0}\overline v)]\\
\displaystyle=\sum_j Q_{1,j}(J^{-s_0}v,\nabla J^{-s_0}v,J^{-s_0}\overline v,\nabla J^{-s_0}\overline v)\partial_{x_j} v\\
\displaystyle+\sum_j Q_{2,j}(J^{-s_0}v,\nabla J^{-s_0}v,J^{-s_0}\overline v,\nabla J^{-s_0}\overline v)\partial_{x_j} \overline v\\
\displaystyle+R_1(\{\partial^{\alpha}_x J^{-s_0}v,\partial^{\alpha}_x J^{-s_0}\overline v\}_{|\alpha|\le s_0/2+1})\Psi_{p_1} v\\
\displaystyle+R_2(\{\partial^{\alpha}_xJ^{-s_0}v,\partial^{\alpha}_x J^{-s_0}\overline
v\}_{|\alpha|\le s_0/2+1})\Psi_{p_2}\overline v, 
\end{array}
$$
where $\,Q_{1,j}\,,\,Q_{2,j}\,,\,j=1,\dots,n\,$, and $\,R_1\,,\,R_2\,$ are polynomials with no constant 
terms and $\,p_1\,,\,p_2\in S^0_{1,0}\,$. The right hand linear factors in the $\,2n+2\,$ terms above 
arise from the highest order derivative in each term of the Leibniz sum. Now let
$$
\begin{array}{l}
\tilde Q_1(v)=\\
\displaystyle\sum_j[Q_{1,j}(J^{-s_0}v,\nabla J^{-s_0}v,J^{-s_0}\overline v,\nabla J^{-s_0}\overline v)-Q_{1,j}(u_0,\nabla
u_0,\overline{u_0},\nabla \overline{u_0})]\partial_{x_j} v,
\\
\tilde Q_2(v)=\\
\displaystyle\sum_j\left[Q_{2,j}(J^{-s_0}v,\nabla J^{-s_0}v,J^{-s_0}\overline v,\nabla J^{-s_0}\overline v)-Q_{2,j}(u_0,\nabla
u_0,\overline{u_0},\nabla \overline{u_0})\right]\partial_{x_j} v,
\\
\tilde R_1(v)=\\
\left[R_1(\{\partial_x^{\alpha} J^{-s_0}v,\partial_x^{\alpha}
 J^{-s_0}\overline v\}_{|\alpha|\le s_0/2+1})-R_1(\{\partial_x^{\alpha} u_0,\partial_x^{\alpha}\overline {u_0}\}_{|\alpha|\le s_0/2+1})\right]\Psi_{p_1}v,
\\
\tilde R_2(v)=\\
\left[R_2(\{\partial_x^{\alpha}J^{-s_0}v,\partial_x^{\alpha} J^{-s_0}\overline v\}_{|\alpha|\le
s_0/2+1})-R_2(\{\partial_x^{\alpha} u_0,\partial_x^{\alpha}\overline {u_0}\}_{|\alpha|\le s_0/2+1})\right]\Psi_{p_2}\overline v,
\end{array}
$$
and
$$
\begin{array}{l}
\tilde b_3(x,\xi)=\tilde b_1(x,\xi)+\sum_j Q_{1,j}(u_0,\nabla u_0,\overline{u_0},\nabla \overline{u_0})(x) i \xi_j,\\
\vec {\tilde b_4} (x)=\vec b_2(x)+
\vec  Q_{2}(u_0,\nabla u_0,\overline{u_0},\nabla \overline{u_0})(x),
\\
\tilde c_3=\tilde c_1+R_1\left(\left\{\partial_x^{\alpha} u_0,\partial_x^{\alpha}\overline {u_0}\right\}_{|\alpha|\le s_0/2+1}\right) p_1,
\\
\tilde c_4=\tilde c_2+R_2\left(\left\{\partial_x^{\alpha} u_0,\partial_x^{\alpha}\overline {u_0}\right\}_{|\alpha|\le s_0/2+1}\right) p_2.
\end{array}
$$
Then it suffices to solve the following nonlinear equation for $v$,
$$
\left\{
\begin{array}{l}
\partial_t v=i{\cal L} v+\Psi_{\tilde b_3}v+\vec {\tilde b_4}\cdot \nabla\overline v+
\Psi_{\tilde c_3}v+\Psi_{\tilde c_4}\overline v\\
\quad +\tilde Q_1(v)+\tilde Q_2(v)+\tilde R_1(v)+\tilde R_2(v)\\ 
v(x,0)=J^{s_0} u_0(x).
\end{array}
\right.
$$
This corresponds to solving the integral equation
$$
v(t)=W_1(t) J^{s_0}u_0+ \int_0^t W_1(t-t^{\prime})[\tilde Q_1(v)+\tilde Q_2(v)+\tilde R_1(v)+\tilde R_2(v)](t^{\prime})dt^{\prime},
$$
where $\,w(t)=W_1(t)w_0\,$ is the solution of the linear homogeneous equation
$$
\left\{
\begin{array}{l}
\partial_t w=i{\cal L} w+\Psi_{\tilde b_3}w+\vec {\tilde b_4}\cdot \nabla\overline w+\Psi_{\tilde c_3}w+\Psi_{\tilde c_4}\overline w,\\
w(x,0)=w_0(x).
\end{array}
\right.
$$
The solution of the integral equation is a fixed point of the following map which will turn out to be a
contraction on a suitable function space. Let
$$
\begin{array}{l}
\left[\Phi_{u_0}(w)\right](t)=W_1(t) J^{s_0}u_0\\
\displaystyle+\int_0^t W_1(t-t^{\prime})[\tilde Q_1(v)+\tilde Q_2(v)+\tilde R_1(v)+\tilde R_2(v)](t^{\prime})dt^{\prime},\\
\\
\displaystyle\lambda_1(w)=\sup_{0\le t\le T}\|w(t)\|_{H^{1/2}},\\
\displaystyle\lambda_2(w)=\left(\int_0^T\int_{\BbbR^n}|J^1 w(x,t)|^2\langle x\rangle^{-2N} dxdt\right)^{1/2},\\
\displaystyle\lambda_3(w)=\sup_{0\le t\le T}\|\langle x\rangle ^{2N} 
\partial_t w(t)\|_{H^{-s_0/2+n/2+7/2}},\\
\displaystyle\Lambda (w)=\max \left\{\lambda_j(w)\,:\,j=1,2,3\right\},\\
\displaystyle\lambda_4(w)=\sup_{0\le t\le T}\|\langle x\rangle ^{2N} w(t)\|_{H^{-s_0/2+n/2+7/2}},
\\
\displaystyle X^{a}_T=\left\{w\,:\,\BbbR^n\times \left[0,T\right]\to \BbbC\,:\,\Lambda (w)\le a\right\}.
\end{array}
$$
For suitable $\,a\,$ and sufficiently small $\,T\,$ it will be shown that $\,\Phi_{u_0}\,$ maps the
complete metric space $\,X^a_T\,$ into $\,X^a_T\,$ and is a contraction.

Let $\,w\in X^a_T\,$. The first goal is to show that $\,\Phi_{u_0}(w)\in X^a_T\,$. Notice that $\,v=\Phi_{u_0}(w)\,$
solves the linear equation
$$
\left\{
\begin{array}{l} 
\displaystyle\partial_t v=i{\cal L} v+\Psi_{\tilde b_3}v+\vec {\tilde b_4}\cdot \nabla
\overline v+\Psi_{\tilde c_3}v+\Psi_{\tilde c_4}\overline v\\
\quad +\tilde Q_1(w)+\tilde Q_2(w)+\tilde R_1(w)+\tilde R_2(w)\\
v(x,0)=J^{s_0} u_0(x).\end{array}
\right.
$$
By the linear smoothing effect,
$$
\begin{array}{l}
\lambda_1^2(\Phi_{u_0}(w))+\lambda_2^2(\Phi_{u_0}(w))
\le c\{\|u_0\|^2_{H^{s_0+1/2}}\\
\displaystyle+\int_0^T\int_{\BbbR^n}|\tilde Q_1(w)|^2\langle x\rangle^{2N}dxdt+
\int_0^T\int_{\BbbR^n}|\tilde Q_2(w)|^2\langle x\rangle^{2N}dxdt\\
\displaystyle+\int_0^T\int_{\BbbR^n}|\tilde R_1(w)|^2\langle x\rangle^{2N}dxdt+\int_0^T\int_{\BbbR^n}|\tilde
 R_2(w)|^2\langle x\rangle^{2N}\}\\
\displaystyle=c\{\|u_0\|^2_{H^{s_0+1/2}}+I_1+I_2+II_1+II_2\}.
\end{array}
$$
Concerning the first term,
$$
\begin{array}{rcl}
I_1&=&\displaystyle\int_0^T\int_{\BbbR^n}|\vec {Q_1}(J^{-s_0}w,\nabla J^{-s_0} w,J^{-s_0}\overline w,
\nabla J^{-s_0}\overline w)\\
&&\displaystyle-\vec Q_1(u_0,\nabla u_0,\overline {u_0},\nabla \overline {u_0})]\cdot\nabla w|^2\langle x\rangle^{2N}dxdt\\
&\le&\displaystyle c\sup_{0\le t\le T}\|\langle x\rangle^{2N}[\vec Q_1
(J^{-s_0}w,\nabla J^{-s_0} w,J^{-s_0}\overline w,\nabla J^{-s_0}\overline w)\\
&&\displaystyle-\vec Q_1(u_0,\nabla u_0,\overline {u_0},\nabla \overline {u_0})]\|^2_{L^{\infty}_{x}}
\int_0^T\int_{\BbbR^n}|\nabla w|^2\langle x\rangle^{-2N}dxdt\\
&=&cI_{1a}I_{1b}.
\end{array}
$$
For $\,j\in\{1,\dots,n\}\,$ one has $\,\partial_{x_j}=(\partial_{x_j}J^{-1})J^1\,$. Observe that
$\,\partial_{x_j}J^{-1}\,$ is a 0th order classical pseudo-differential operator and is hence bounded on
$\,L^2(\langle x\rangle^{-2N}dx)\,$ by Lemma \ref{le2.3.3} in Section \ref{sec2}. Therefore,
$$
\displaystyle I_{1b}\le c\int_0^T\int_{\BbbR^n}|J^1w|^2\langle x\rangle^{-2N}dxdt=c\lambda_2(w).
$$
In order to estimate $\,I_{1a}\,$, let $\,j\in\{1,\dots,n\}\,$ and $\,t\in[0,T]\,$. Then
$$
\begin{array}{l}
\vec {Q_1}
(J^{-s_0}w,\nabla J^{-s_0} w,J^{-s_0}\overline w,\nabla J^{-s_0}\overline w)-\vec {Q_1}(u_0,\nabla u_0,\overline {u_0},\nabla \overline {u_0})\\
\displaystyle=\int_0^t\partial_t[\vec {Q_1}(J^{-s_0}w,\nabla J^{-s_0} w,J^{-s_0}\overline w,\nabla J^{-s_0}\overline w)(t^{\prime})]dt^{\prime}.
\end{array}
$$
By the product rule for $\,t$-differentiation, Lemma \ref{le6.1.2} and Sobolev\rq{}s theorem,
$$
I_{1a}\le c T\lambda_3^2(w)S_1(\lambda_1^2(w))
$$
where $\,S_1\,$ is some polynomial of one variable. In the rest of the proof, $\,S_j\,,\,j=2,3\dots\,$, will denote other such polynomials. 
Combining estimates for $\,I_{1a}\,$ and
$\,I_{1b}\,$,
$$
I_{1}\le c T\lambda_3^2(w) \lambda_2^2(w) S_1(\lambda_1^2(w)).
$$
The estimate for $\,I_2\,$ is similar.
Concerning the last two terms,
$$
\begin{array}{l}
II_1=\\
\displaystyle\int_0^T\int_{\BbbR^n}\Big|\big[R_1(\{\partial^{\alpha}_x
J^{-s_0}w,\partial^{\alpha}_x J^{-s_0}\overline w\}_{|\alpha|\le s_0/2+1})-
R_1(\{\partial^{\alpha}_xu_0,\partial^{\alpha}_x \overline {u_0}\}_{|\alpha|\le s_0/2+1})\big]\cdot\\
\cdot \Psi_{p_1}w\Big|^2\langle x\rangle^{2N} dxdt\\
\le cT\lambda_1^2(w)\cdot\\
\cdot \|\langle x\rangle^{2N}[R_1(\{\partial^{\alpha}_x
J^{-s_0}w,\partial^{\alpha}_x J^{-s_0}\overline w\}_{|\alpha|\le s_0/2+1})-
R_1(\{\partial^{\alpha}_xu_0,\partial^{\alpha}_x \overline {u_0}\})_{|\alpha|\le s_0/2+1})]\|^2_{L^{\infty}_x}
\\
=cT\lambda_1^2(w)II_{1a}.
\end{array}
$$
To estimate $\,II_{1a}\,$ let $\,t\in [0,T]\,$. Then
$$
\begin{array}{l}
R_1(\{\partial^{\alpha}_x
J^{-s_0}w,\partial^{\alpha}_x J^{-s_0}\overline w\}_{|\alpha|\le s_0/2+1})(t)-
R_1(\{\partial^{\alpha}_x u_0,\partial^{\alpha}_x \overline {u_0}\}_{|\alpha|\le s_0/2+1})\\
\displaystyle=\int_0^t\partial_t\big[R_1(\{\partial^{\alpha}_xJ^{-s_0}w,\partial^{\alpha}_x J^{-s_0}\overline w\}_{|\alpha|\le s_0/2+1})(t^{\prime})\big]dt^{\prime}.
\end{array}
$$
By the product rule for $\,t$-differentiation, Lemma \ref{le6.1.2} and Sobolev\rq{}s theorem,
$$
II_{1a}\le C T\lambda_3^2(w)S_2(\lambda_1^2(w)).
$$
Hence
$$
II_{1}\le C T^2\lambda_3^2(w) \lambda_1^2 S_2(\lambda_1^2(w)).
$$
The estimate for $II_2\,$ is similar. Combining estimates for $I_1\,,\,I_2\,,\,II_1\,$ and
$\,II_2$, one gets that
$$
(\lambda_1^2+\lambda_2^2)(\Phi_{u_0}(w))\le c\|u_0\|^2_{H^{s_0+1/2}}+cT(1+T)\Lambda^4(w)S_3(\Lambda^2(w)).
$$
Next $\lambda_4\,$ will be estimated and then used in the estimate of $\,\lambda_3$. By Lemma \ref{le6.1.3}(b),
$$
\begin{array}{l}
\lambda_4(\Phi_{u_0}(w))\le
\displaystyle\sup_{0\le t\le T}\|\langle x\rangle^{2N} W_1(t)J^{s_0}u_0\|_{H^{-s_0/2+n/2+7/2}}
\\
\displaystyle+T\sup_{0\le t^{\prime}\le t\le T}\|\langle x\rangle^{2N} W_1(t-t^{\prime})[\tilde Q_1(w)+
\tilde Q_2(w)\\
\displaystyle+\tilde R_1(w)+\tilde R_2(w)](t^{\prime})\|_{H^{-s_0/2+n/2+7/2}}\\
\le c (1+T^N)\|\langle x\rangle^{2N} J^{s_0} u_0\|_{H^{-1}}+ c T (1+T^N)\cdot\\
\displaystyle\cdot \{\sup_{0\le t\le T}\|\langle x\rangle^{2N} \tilde Q_1(w(t))\|_{H^{-1}}
+\sup_{0\le t\le T}\|\langle x\rangle^{2N} \tilde Q_2(w(t))\|_{H^{-1}}\\
\displaystyle+\sup_{0\le t\le T}\|\langle x\rangle^{2N} \tilde R_1(w(t))\|_{H^{-1}}+
\sup_{0\le t\le T}\|\langle x\rangle^{2N} \tilde R_2(w(t))\|_{H^{-1}}\}\\
=c(1+T^N)\|\langle x\rangle^{2N} J^{s_0} u_0\|_{H^{-1}}+c T (1+T^N)\{I_1+I_2+II_1+II_2\}.
\end{array}
$$
By the product rule for $\,x_j$-differentiation and the fundamental theorem of calculus in the
$\,t$-variable,
$$
\begin{array}{c}
\tilde Q_1(w(t))=\\
\displaystyle\sum_j\partial_{x_j}[w(t)\int_0^t\partial_t[Q_{1,j}(J^{-s_0}w,\nabla J^{-s_0} w,J^{-s_0}\overline w,\nabla J^{-s_0}\overline
w)(t^{\prime})]dt^{\prime}\Big]\\
\displaystyle+w(t)\sum_j\int_0^t\partial_{x_j}\partial_t[Q_{1,j}(J^{-s_0}w,\nabla J^{-s_0} w,J^{-s_0}\overline w,\nabla J^{-s_0}\overline w)(t^{\prime})]dt^{\prime}
\end{array}
$$
so that, by Sobolev\rq{}s theorem,
$$
\begin{array}{l}
\displaystyle I_1\le c\sum_j\sup_{0\le t\le T}\|\langle x\rangle^{2N}w(t)\cdot\\
\displaystyle\cdot\int_0^t\partial_t[Q_{1,j}(J^{-s_0}w,\nabla J^{-s_0} w,J^{-s_0}\overline w,\nabla J^{-s_0}\overline
w)(t^{\prime})]dt^{\prime}\|_{L^2}
\\
\displaystyle+c\sum_j\sup_{0\le t\le T}\|\langle x\rangle^{2N}w(t)\cdot\\
\displaystyle\cdot\int_0^t\partial_{x_j}\partial_t[Q_{1,j}(J^{-s_0}w,\nabla J^{-s_0} w,J^{-s_0}\overline w,\nabla J^{-s_0}\overline
w)(t^{\prime})]dt^{\prime}\|_{L^2}\\
\displaystyle\le c T \lambda_3(w)\lambda_1(w) S_4(\lambda_1(w)).
\end{array}
$$
Similarly,
$$
I_2\le C T \lambda_3(w)\lambda_1(w) S_5(\lambda_1(w)).
$$
Concerning $\,II_1\,$, one has
$$
\tilde R_1(w(t))=\Psi_{p_1} w(t) \int_0^t\partial_t
[R_1(\{\partial_x^{\alpha}J^{-s_0}w,
\partial_x^{\alpha}J^{-s_0}\overline w\}_{|\alpha|\le s_0/2+1})(t^{\prime})]dt^{\prime}.
$$
By Sobolev\rq{}s theorem and the $\,L^2\,$ boundedness of $\,\Psi_{p_1}$,
$$
\begin{array}{l}
\displaystyle II_1=\sup_{0\le t\le T}\|\langle x\rangle^{2N}\Psi_{p_1} w(t) \int_0^t\partial_t\left[R_1(\{\partial_x^{\alpha}J^{-s_0}w,\partial_x^{\alpha}J^{-s_0}\overline w\}_{|\alpha|\le s_0/2+1})(t^{\prime})\right]dt^{\prime}\|_{L^2}\\
\displaystyle\le c T \lambda_3(w)\lambda_1(w)S_6(\lambda_1(w)).
\end{array}
$$
Similarly
$$
II_2\le c T \lambda_3(w)\lambda_1(w)S_7(\lambda_1(w)).
$$
Combining estimates for $\,I_1\,,\,I_2\,,\,II_1\,$ and $\,II_2$, it follows that
$$
\lambda_4(\Phi_{u_0}(w))\le c (1+T^N) \|\langle x\rangle^{2N} J^{s_0} u_0\|^2_{H^{-1}}+cT^2(1+T^N)\Lambda^2(w)S_8(\Lambda(w)).
$$
This estimate of
$\,\lambda_4\,$ will be used in that of $\,\lambda_3\,$. Let $\,v=\Phi_{u_0}(w)\,$ and note that
$$
\partial_t v=(i{\cal L}+\Psi_{\tilde b_3}+\Psi_{\tilde c_3}) v+
(\vec {\tilde b_4}\cdot \nabla+
\Psi_{\tilde c_4})\overline v+
\tilde Q_1(w)+\tilde Q_2(w)+\tilde R_1(w)+\tilde R_2(w)
$$
where $\,i{\cal L}+\Psi_{\tilde b_3}+\Psi_{\tilde c_3}\,$
and $\,\vec {\tilde b_4}\cdot \nabla+
\Psi_{\tilde c_4}\,$
 are classical pseudo-differential operators of order 2 and 1, respectively. By
Lemma \ref{le6.1.2},
$$\begin{array}{c}
\displaystyle\lambda_3(v)=\lambda_4(\partial_t v)\le
c\sup_{0\le t\le T}\|\langle x\rangle^{2N} v\|_{H^{-s_0/2+n/2+11/2}}\\
\displaystyle+c\sup_{0\le t\le T}\big\|\langle x\rangle^{2N}\big[\tilde Q_1(w)+\tilde Q_2(w)+\tilde R_1(w)+\tilde R_2(w)\big]\big\|_{H^{-s_0/2+n/2+7/2}}.
\end{array}
$$
The first term can be estimated as $\,\lambda_4(v)\,$ since $\,s_0\,$ is sufficiently large. The second term is dominated by
$\,I_1+I_2+II_1+II_2\,$ in the estimate of $\,\lambda_4(v)\,$. Hence
$$
\lambda(\Phi_{u_0}(w))\le c (1+T^N) \|\langle x\rangle^{2N} J^{s_0} u_0\|^2_{H^{-1}}+cT(1+T^{N+1})\Lambda^2(w)S_9(\Lambda(w)).
$$
$\,T\,$ will later be chosen small, so it can be assumed that $\,T\le 1\,$. 
Combining estimates for $\,\lambda_1\,,\,\lambda_2\,$ and $\,\lambda_3\,$,
$$
\begin{array}{c}
\Lambda(\Phi_{u_0}(w))\le c(\|u_0\|_{H^{s_0+1/2}}+\|\langle x\rangle^{2N}J^{s_0}u_0\|_{H^{-1}})\\
+cT^{1/2}\Lambda^2(w)(1+\Lambda^{\rho}(w))
\end{array}
$$
where
$\,\rho>0\,$. First fix
$$
a>2C(\|u_0\|_{H^{s_0+1/2}}+
\|\langle x\rangle^{2N} J^{s_0} u_0\|_{H^{-1}}).
$$
Next choose
$$
\displaystyle T=\min\{(2c a(1+a^{\rho}))^{-2},1\}.
$$
Then
$\,\Lambda(\Phi_{u_0}(w))\le a\,$ so $\,\Phi_{u_0}\,$ maps $\,X^a_T\,$ into $\,X^a_T\,$. To see that
$\,\Phi_{u_0}\,$ is a contraction on $\,X^a_T\,$ for sufficiently small $\,T\,$, notice that
$$
\begin{array}{l}
\displaystyle (\Phi_{u_0}(w_1)-\Phi_{u_0}(w_2))(t)=
\int_0^t W_1(t-t^{\prime})[(\tilde Q_1(w_1)-\tilde Q_1(w_2))
\\
\displaystyle +(\tilde Q_2(w_1)-\tilde Q_2(w_2))+
(\tilde R_1(w_1)-\tilde R_1(w_2))+(\tilde R_2(w_1)-\tilde R_2(w_2))](t^{\prime})dt^{\prime}.
\end{array}
$$
The estimates used in showing $\,\Phi_{u_0}\,:\,X^a_T\to X^a_T\,$ therefore give
$$
\Lambda(\Phi_{u_0}(w_1)-\Phi_{u_0}(w_2))\le C T^{1/2}\Lambda (w_1-w_2)a(1+a^{\rho}).
$$
Now choose $\,T=\min\{1,2C a(1+a^{\rho}))^{-2}\}\,$. Then
$$
\displaystyle 
\Lambda(\Phi_{u_0}(w_1)-\Phi_{u_0}(w_2))\le
\frac 12\Lambda (w_1-w_2),
$$
so $\,\Phi_{u_0}\,$ is a contraction. By Banach\rq{}s contraction mapping principle there is a unique
fixed point of $\,\Phi_{u_0}\,$ which solves the nonlinear equation.
Now let $\,u\,$ and $\,v\,$ be solutions of (\ref{eq6.2.1}) with initial values $\,u_0\,$ and $\,v_0\,$
respectively. Then the estimates above give
$$
\Lambda (u-v)\le
C(\|u_0-v_0\|_{H^s}+\|\langle x\rangle^{2N} J^{s-3/2}(u_0-v_0)\|_{L^2})+
CT^{1/2}\Lambda(u-v)
$$
if $\,u\,$ and $\,v\,$ are both in an open ball in $\,X^s_T\,$. Here $\,C\,$ depends on the radius of
the ball. Now choose $\,T^{\prime}>0\,$ so small that $\,C(T^{\prime})^{1/2}\le\frac 12\,$. The proof
of Theorem \ref{th6.2.2} is complete since
$$
\|\,\,\cdot\,\,\|_{H^s}+\|\langle x\rangle^{2N} J^{s-3/2}\,\,\cdot\|_{L^2}
$$
is dominated by finitely many seminorms in $\,{\cal S}$ (see (\ref{eq3.1.2})).


\newpage
        



 
\begin{thebibliography}{99}
\bibitem{1}	
Calder\'on, A.P., and Vaillancourt,R.,
{\em On the boundedness of pseudo-differential operators},
J. Math. Soc. Japan
{\bf 23}
(1971),
374--378.
\bibitem{2}
 Chihara, H.,
{\em Local existence for semilinear Schr\"odinger equations},
Math. Japan
{\bf 42}
(1995),
35--51.
\bibitem{3}
 Constantin, P., and Saut, J. C., 
{\em Local smoothing properties of dispersive equations}
J. Amer. Math. Soc.
{\bf 1}
(1989),
 413--446.
 \bibitem{4}
Craig, W., Kappeler, T., and Strauss, W.,
{\em Microlocal dispersive smoothing for the Schr\"odinger equation},
Comm. Pure Appl. Math.
{\bf 48}
(1995),
769--860.
 \bibitem{5}
 Davey, A., and Stewartson, K.,
{\em On three--dimensional packets of surface waves}
Proc. Roy. Soc. London Ser.
{\bf 338}
(1974),
101--110.
 \bibitem{6}
Djordjevic, V. D., and Redekopp, L. G.,
{\em On two-dimensional packets of capillary-gravity waves},
J. Fluid Mech.
{\bf 79}
(1977),
703-714.


\bibitem{7}
Doi, S.,
{\em On the Cauchy problem for Schr\"odinger type equations and the regularity of solutions},
J. Math. Kyoto Univ.
{\bf 34}
(1994),
319--328.


\bibitem{8}
Doi, S.,
{\em Remarks on the Cauchy problem for Schr\"odinger--type equations}, 
Comm. Partial Differential Equations
{\bf 21}
(1996),
163--178.


\bibitem{9}
Doi, S. ,
{\em Smoothing effects for Schrodinger evolution equation and global behavior of geodesic flow}
Math.-Ann. 
{\bf  318}
(2000),
355--389.

\bibitem{10}
  Hayashi, N.,
   {\em Global existence of small analytic solutions to nonlinear Schr\"odinger equations},
 Duke Math. J.
  {\bf 62}
  (1991),
  575--592.

\bibitem{11}
 Hayashi, N., and Ozawa, T.,
{\em Remarks on nonlinear Schr\"odinger equations in one space dimension},
Differential Integral Equations
{\bf 7}
(1994),
453--461.

\bibitem{12}
H\"ormander, L.,
{\em Pseudo-differential operators and non--elliptic boundary problems},
Ann. of Math. (2)
{\bf 83}
(1966),
129--209.

\bibitem{13}
Ichinose, W.,
{\em On $L^2$ well-posedness of the Cauchy problem for
Schr\"odinger type equations on a Riemannian manifold and Maslov
theory},
Duke Math. J.
{\bf 56}
(1988),
549--588.

  \bibitem{14}
  Ishimori, Y.,
{\em  Multi vortex solutions of a two dimensional nonlinear wave 
equation},
Progr. Theor. Phys.
  {\bf 72}
  (1984),
  33--37.

\bibitem{15}
Kato, T.,
{\em Quasi--linear equations of evolution, with applications to partial differential equations},
Lecture Notes in  Math.  
{\bf 448}
(1975),
Springer--Verlag,
  27--70.
  
\bibitem{16}
 Kato, T., 
   {\em On the Cauchy problem for the (generalized) Korteweg-de Vries equation},
 Advances in Math. Supp. Studies, Studies in Applied Math.
   {\bf 8}	
  (1983),	
   93--128.	
	
\bibitem{17}	
  Kenig, C. E., Ponce, G., and Vega, L.,	 
  {\em Oscillatory integrals and regularity of dispersive equations},	
 Indiana University Math. J.	
  {\bf 40}	
  (1991),	
  33--69.
  
 \bibitem{18}	
Kenig, C. E., Ponce, G., and Vega, L.,
{\em Small solutions to nonlinear Schr\"odinger equations},
Ann. Inst. H. Poincar\'e Anal. Non Lin\'eaire
{\bf 10}
(1993),
255--288.

\bibitem{19}
Kenig, C. E., Ponce, G., and Vega, L.,
{\em On the Zakharov and Zakharov--Schulman systems}, 
J. Funct. Anal.
{\bf 127}
(1995),
204--234.

\bibitem{20}
Kenig, C. E., Ponce, G., and Vega, L.,
{\em On the smoothing properties of some dispersive
hyperbolic systems},
 Nonlinear waves(Sapporo, 1995),
 221-229:
GAKUTO Internat. ser. Math. Sci. Appl.,10, Tokyo
(1997).

\bibitem{21}
Kenig, C. E., Ponce, G., and Vega, L.,
{\em On the Cauchy problem for linear Schr\"odinger systems with variable coefficient lower order
terms},
Harmonic analysis and number theory (Montreal, PQ, 1996),
205--22:
CMS Conf. Proc., 21, Amer. Math. Soc., Providence, RI
(1997).

\bibitem{22}
Kenig, C. E., Ponce, G., and Vega, L.,
 {\em Smoothing effects and local existence theory for the generalized nonlinear Schr\"odinger
equations},
Invent. Math. 
{\bf 134}
(1998),
489-545.

\bibitem{23}
Klainerman, S.,
{\em Long time behavior of solutions to nonlinear evolution equations},
Arch. Ration. Mech. and Analysis
{\bf 78}
(1981),
73--98.

\bibitem{24}
Kruzhkov, D. J., and Faminskii, A. V.,
{\em Generalized solutions for the Cauchy problem for the Korteweg-de Vries equation},
Math. USSR Sb
{\bf 48}
(1990),
93--138.

\bibitem{25}
 Kumano--Go, H.,
{\em Pseudo-differential operators},
MIT Press, Cambridge 
(1981).

\bibitem{26}
Mizohata, S.,
{\em On the Cauchy Problem},
Notes and Reports in Mathematics in Science and Engineering,3, Science press and Academic Press
(1985).

\bibitem{27}
Rolvung, C.,
{\em Non-isotropic Schr\"odinger equations}
PhD. dissertation, University of Chicago
(1998).

\bibitem{28}
  Simon, J., and Taflin, E.,
 {\em Wave operators and analytic solutions for systems of nonlinear Klein-Gordon equations 
  and of non-linear Schr\"odinger equations}
 Comm. Math. Phys.
 {\bf 99}
 (1985),
541--562.

\bibitem{29}
 Sj\"olin, P.,
{\em Regularity of solutions to the Schr\"odinger equations},
 Duke Math. J.
  {\bf 55}
 (1987),
 699--715.

\bibitem{30}
Stein, E. M.,
{\em Harmonic analysis: real--variable methods, orthogonality, and oscillatory integrals},
Princeton University Press 
(1993).

\bibitem{31}
 Vega, L. ,
 {\em The Schr\"odinger equation:  pointwise convergence to the initial date},
Proc. Amer. Math. Soc.
 {\bf 102}
  (1988),
 874--878.

\bibitem{32}
Zakharov, V. E., and Schulman, E. I.,
{\em Small solutions to nonlinear Schr\"odinger equations},
Physica
{\bf 1D}
(1980),
185--250.
	
\end{thebibliography}
 \end{document}